\documentclass[11pt,a4paper,twoside]{article}

\usepackage{amsmath, exscale, epsfig, amssymb}

\small
\normalsize
\usepackage{amsfonts}
\usepackage{amsbsy}
\usepackage{fancyhdr}
\usepackage{amscd}
\usepackage{theorem}
\usepackage{epsf}
\usepackage{psfrag}
\usepackage{multicol}
\allowdisplaybreaks

\def\build#1_#2^#3{\mathrel{\mathop{\kern 0pt#1}\limits_{#2}^{#3}}}
\newcommand{\beq}{\begin{eqnarray*}}
\newcommand{\eeq}{\end{eqnarray*}}
\newcommand{\ben}{\begin{enumerate}}
\newcommand{\een}{\end{enumerate}}
\newcommand{\beqs}{\begin{eqnarray*}&\displaystyle}
\newcommand{\eeqs}{&\end{eqnarray*}}

\newcommand{\un}{{\bf 1}}

\newtheorem{theorem}{Theorem}[section]
\newtheorem{lemma}[theorem]{Lemma}
\newtheorem{proposition}[theorem]{Proposition}

{\theorembodyfont{\rmfamily}\newtheorem{remark}{Remark}[section]}
{\theorembodyfont{\rmfamily}\newtheorem{notation}{Notation}[section]}
{\theorembodyfont{\rmfamily}}

\newcommand{\btl}{\mathbb{T}_{{\bf \ell_1 }}} 
\newcommand{\dl}{{\bf d}}

\newcommand{\bE}{\mathbb{E}}
\newcommand{\E}{\mathbb{E}}
\newcommand{\bF}{\mathbb{F}}

\newcommand{\bN}{\mathbb{N}}
\newcommand{\bP}{\mathbb{P}}
\newcommand{\bR}{\mathbb{R}}

\newcommand{\bT}{\mathbb{T}}
\newcommand{\bU}{\mathbb{U}}
\newcommand{\bX}{\mathbb{X}}
\newcommand{\cA}{\mathcal{A}}

\newcommand{\cC}{\mathcal{C}}
\newcommand{\cE}{\mathcal{E}}
\newcommand{\cF}{\mathcal{F}}
\newcommand{\cG}{\mathcal{G}}

\newcommand{\cM}{\mathcal{M}}
\newcommand{\cN}{\mathcal{N}}
\newcommand{\cP}{\mathcal{P}}

\newcommand{\cT}{\mathcal{T}}

\newcommand{\lgeo}{[\![}
\newcommand{\rgeo}{]\!]}
\def\cqfd{ \hfill $\blacksquare$ }

\begin{document}

\title{GROWTH OF L\'EVY TREES.}
\author{Thomas Duquesne\thanks{Universit\'e Paris 11, Math\'ematiques,
    91405 Orsay Cedex, France; supported by NSF Grants DMS-0203066 
and DMS-0405779} \and Matthias 
Winkel\thanks{University of
  Oxford, Department of Statistics, 1 South Parks Road, Oxford OX1 3TG, UK; 
        email winkel@stats.ox.ac.uk; 
supported by Aon and the Institute of Actuaries, le d\'epartement de
math\'ematique de l'Universit\'e d'Orsay and NSF Grant DMS-0405779}}
\maketitle

\begin{abstract}  We construct random locally compact real trees called
 L\'evy trees that are the genealogical trees associated with  
continuous-state branching processes. More precisely, we define a
growing family of discrete Galton-Watson trees with i.i.d. exponential
branch lengths that is consistent under Bernoulli percolation on
leaves; we define the L\'evy tree as the limit of this growing family
with respect to the Gromov-Hausdorff topology on metric spaces. This
elementary approach notably includes supercritical trees and does not
make use of the height process introduced by Le Gall and Le Jan to
code the genealogy of (sub)critical continuous-state branching
processes. 
We construct the mass measure
of L\'evy trees and we give a decomposition along the 
ancestral subtree of a Poisson
sampling directed by the mass measure.

\end{abstract}

\section{Introduction}

Continuous-state branching processes have been introduced by Jirina 
\cite{Ji} and Lamperti \cite{La1, La2, La3}. They 
are the continuous analogues of the
Galton-Watson Markov chains. Recall that the distribution of a
continuous-state branching process is characterized by a real-valued
function $\psi$ defined on $[0, \infty)$ that is of the form 
\begin{equation}
\label{LevKhin1}
 \psi(c)=\alpha c+\beta c^2+\int_{(0,\infty)}
(e^{-cx}-1+cx{\bf 1}_{\{x<1\}})\Pi (dx), 
\end{equation}
where $\alpha \in \bR$, $\beta\geq 0$, and $\Pi$ is the Lévy measure which satisfies 
$$\int_{(0, \infty)}(1\wedge x^2)\Pi (dx)<\infty .$$ 
$\psi$ is called the {\it branching mechanism} 
of the continuous-state branching process. More precisely, 
$Z=(Z_t,t\geq 0)$ is a continuous-state branching 
process with branching mechanism $\psi $ (a CSBP($\psi $) for short) 
iff it is a $[0, \infty ]$-valued 
Feller process whose transition kernel is characterized by 
$$\bE \left[ \exp (-\lambda Z_{s+t}) \left|  Z_s \right. \right]= \exp (-u(t,\lambda)Z_s) \; , $$
where $u$ is the unique non-negative solution of 
$$ \partial_t u(t, \lambda)=-\psi (u(t, \lambda)) \quad {\rm and }\quad
u(0,\lambda)=\lambda \; ,\; t,\lambda \geq 0 . $$
This equation can be rewritten in the following integrated form 
\begin{equation} 
\label{transcsbp}
\int_{u(t,\lambda)}^{\lambda}  \frac{dc}{\psi (c)}\; = \; t \; . 
\end{equation}
We shall mostly restrict our attention to the case where 
$Z_t$ has a finite expectation which is equivalent
to the fact that the right derivative of $\psi $ at $0$ is finite. 
We denote this right derivative by 
$m:=\psi'(0+)$. 
We refer to the case  $m \in [-\infty , 0)$ (resp. $m=0$ and $m\in (0, +\infty)$) as to the 
{\it supercritical case} (resp. {\it critical
  case} and {\it subcritical case}).

We shall often assume that $Z$ has a positive probability of
extinction which is equivalent to the following analytical condition 
\begin{equation}
\label{posextinct}
\int^{\infty}\frac{dc}{\psi(c)}< \infty 
\end{equation}
(see \cite{Bi2}). In that case, we have
$$ \bP \left(\left. \exists t\geq 0 \; :\; Z_t=0 \right|Z_0=a\right) =\exp (-a\gamma) , $$
where $\gamma$ is the largest root of the equation $\psi(c)=0$
(observe that $\gamma >0$ only in the supercritical case: $m<0$). 
If (\ref{posextinct}) is not satisfied then the 
underlying Lévy tree will fail to be separable.

The main goal of this paper is to construct an $(a,\psi)$-L\'evy tree
that can be interpreted as the genealogical tree of a population whose 
size evolves according a CSBP($\psi$) $Z$ with initial state
$Z_0=a$. We proceed by approximating the L\'evy tree by 
Galton-Watson trees with exponential edge lengths. 
More precisely, recall that 
a Galton-Watson tree with exponential edge lengths is the genealogical tree of
an ancestor and its descendants, where individuals have independent and 
identically exponentially distributed lifetimes with a rate
$c\in(0,\infty)$, and produce offspring at the end of their lives  
independently according to an offspring distribution $\xi$ on 
$\{0,2,3,\ldots\}$. Instead of one single tree, we rather consider 
a random number of independent Galton-Watson trees, the random 
number of ancestors being a Poisson 
random variable with parameter $a$. 
We call such a forest a 
\em Galton-Watson forest \em (a GW($\xi,c,a$)-forest for short).

Let  $F$ be a GW($\xi,c,a$)-forest. We perform a 
\em Bernoulli leaf colouring \em  on $F$ i.e.\ we attach independent
Bernoulli marks with parameter $p$ to all leaves and 
interpret a mark $0$ as black, a mark $1$ as red. 
An elementary calculation will show that the subforest 
$F_b$ 
spanned by the black leaves and the root 
is a GW($\xi_b,c_b,a_b$)-forest, where $\xi_b$, $c_b$, $a_b$ 
are calculated explicitly in terms of $\xi$,  $c$, $a$ and $p$ (see Lemma \ref{percola}).

  One of the aims of the paper is to construct a family $(\cF_{\lambda} \, ; \, \lambda\geq 0)$ of  
random trees such that  for any $ \lambda \geq 0$, 
$\cF_{\lambda}$ is a GW($\xi_\lambda,c_\lambda,a_\lambda$)-forest and
that is consistent under Bernoulli leaf colouring: namely, 
for any $0\leq \mu \leq  \lambda$, we want $\cF_{\mu}$ to be the black subtree obtained from  $\cF_{\lambda}$
by a Bernoulli leaf colouring with parameter $1-p=\mu / \lambda $. Theorem \ref{prolong} asserts that the distribution of 
such a leaf colouring consistent family can be parametrized by $(a, \psi)$, where $a\in(0,\infty)$ and  
$\psi$ is the branching mechanism of a continuous-state branching process 
(CSBP($\psi$)) that is of the form (\ref{LevKhin1}); more precisely we have
\begin{equation} 
\label{specif}
a_\lambda=a \psi^{-1}(\lambda),\quad c_\lambda=\psi^\prime(\psi^{-1}(\lambda)),\quad\varphi_\lambda(r)=\sum_{k=0}^\infty\xi_\lambda(k)r^k=r+\frac{\psi((1-r)\psi^{-1}(\lambda))}{\psi^{-1}(\lambda)\psi^\prime(\psi^{-1}(\lambda))}.
\end{equation}

This offspring distribution appears in \cite{Ho} in the Brownian case and 
Theorem 3.2.1 \cite{DuLG} in the
critical and the subcritical cases as the distribution of the ancestral tree 
corresponding to Poisson marks on 
$[0, \infty)$ via the coding of the L\'evy tree by the height process. 
We refer to Remark \ref{concon} for a detailed discussion 
of the connection of our results and the work in 
\cite{DuLG, DuLG2}.

  Conversely, Proposition \ref{growthproc} asserts that to each 
$(a,\psi)$ there corresponds a growing family 
$(\cF_{\lambda}; \lambda \geq 0)$ of GW-trees with edge lengths, 
consistent under 
Bernoulli leaf colouring as explained before and whose distribution is 
specified by (\ref{specif}). This process can be viewed as a forest-valued continuous-time Markov chain  
whose characteristics are specified by Remark \ref{jumpchain}.

The leaf-colouring consistent forest growth processes that we consider 
can be viewed in a more 
general framework of Markovian forest growth processes. Several schemes to
construct such processes preserving Galton-Watson forests and allowing to 
pass to continuous limits are more or less 
explicit in the literature. They are often more easily described by their
co-transition rules. \em Firstly\em, Neveu \cite{Ne} and Salminen \cite{Sal} erase 
branches in general (non-explosive) Galton-Watson trees with exponential edge lengths 
continuously from their tips. Le Jan \cite{LJ91}, Abraham \cite{Ab92} and 
Pitman \cite{Pit02} reverse the procedure to grow stable/Brownian trees and 
forests from appropriate Galton-Watson trees/forests.
\em Secondly\em, Aldous and Pitman \cite{AlPit98} perform percolation on the edges in 
general Galton-Watson trees (without edge lengths) and retain the connected 
component containing the root, as a tree-valued Markov process as the 
percolation probability varies. They call the procedure pruning of a 
Galton-Watson tree. The viewpoint is to gradually reduce the tree by 
consistently increasing the percolation probability.
Geiger and Kaufmann \cite{GeiKau} discount the offspring distribution to reduce
a given Galton-Watson tree in a size-biased way. This can be seen as a special
case of multiplicity-dependent pruning at vertices. We will see here that it is
also related to the \em third \em scheme of reduction by Bernoulli leaf colouring, that we
study in this paper.

   Let us denote by 
$(Z^\lambda_t)_{t\ge 0}$ the population size process 
associated with $\cF_\lambda$. Assume that $m$ is finite. Then, it is
easy to show that for any $t\geq 0$, a.s.  
 \beqs \frac{1}{\psi^{-1}(\lambda)}Z_t^{\lambda}
\build{\longrightarrow}_{\lambda \to\infty}^{\; }  Z_t,
\eeqs
where $Z$ is a CSBP($\psi$) such that $Z_0=a$. Under the additional
assumptions (\ref{posextinct}) we prove in Theorem \ref{mainresult} an 
a.s.\ convergence for the entire genealogy: 
as in the paper by Evans, Pitman and Winter \cite{EvPitWin}, 
we consider genealogical trees as tree-like metric spaces and more
precisely as 
locally compact rooted real trees, whose 
precise definition is given in Section \ref{realtrees}. We introduce 
the set $\bT$ of root preserving isometry classes of such trees 
equipped with the pointed Gromov-Hausdorff metric $\delta $ 
(see (\ref{defgrohau}) Section \ref{gromovtopo} for the definition); we prove 
in Proposition \ref{treepointed} that $(\bT , \delta)$ is a Polish space. 
This is a simple generalization of the compact case proved in \cite{EvPitWin}. 
Then, we see the growing process of trees 
$(\cF_{\lambda}; \lambda \geq 0)$ as a collection of locally
compact rooted real trees  $(\cF_{\lambda}, d_{\lambda} , \rho)$, 
$\lambda \geq 0$, such that for any 
$0\leq \mu  \leq \lambda$ 
$$ \cF_{\mu} \subset \cF_{\lambda} \quad {\rm and} \quad  
d_{\lambda | \cF_{\mu} \times  \cF_{\mu} }= d_{\mu}  $$
(here $\rho $ stands for the common root of the trees). Then Theorem \ref{mainresult} 
asserts that a.s. 
$$ \delta \left( \cF_{\lambda} , \cF \right)
\build{\longrightarrow}_{\lambda \to\infty}^{\; } 0 \; , $$
where $ \cF$ is the completion of $\bigcup  \cF_{\lambda}$. 
The limiting random tree is called the $(\psi , a)$-{\it L\'evy forest}.

This result is related to the work of Pitman and Winkel 
\cite{PitWink} who perform Bernoulli leaf colouring in the special
case of \em binary \em  
Galton-Watson forests. They show, that in this case, the
forest growth process has independent ``increments'', expressed by a 
composition rule. It can be consistently extended to increase to the Brownian 
forest. This passage to the limit is understood by convergence of coding height
processes via a Donsker type theorem. In the 
critical or subcritical case $m\geq 0$, it is also clear 
(see Remark \ref{concon} for a detailed explanation), 
that the distribution of the 
root preserving isometry class of $\cF$ is the same as the distribution 
induced by the corresponding forest coded by the height process introduced by 
Le Gall and Le Jan \cite{LGLJ1} (see also \cite{DuLG}). Let us mention 
that a framework of real trees and the Gromov-Hausdorff metric has been developed 
for probabilistic applications by Evans in \cite{Ev00}, Evans, Pitman,
Winter in \cite{EvPitWin} and Evans, Winter in \cite{EvWin}.

 In Section \ref{excursionmes}, we define  the 
$\psi$-excursion measure $\Theta$ that can 
be seen as the ``distribution'' of a single 
$\psi $-L\'evy tree. More precisely, Proposition \ref{excuracine}
asserts that the isometry classes of the connected 
components of $\cF \backslash \{ \rho \}$ form a Poisson point process 
on $\bT$ with intensity $a \, \Theta$.

   Our definition of the L\'evy forest also allows to construct the mass measure on $\cF$ denoted by 
${\bf m}$ in the following way: let us denote by ${\bf m}_{\lambda}$ the empirical distribution of the leaves of 
$\cF_{\lambda}$; then Theorem \ref{meslebe} asserts that ${\bf m}_{\lambda} /\lambda$
a.s. converges to ${\bf m}$ for the vague topology of the Radon measures
on $\cF$. It also asserts that the topological support of the mass measure is $\cF$ 
and that the isometry class of the tree spanned by the 
root $\rho$ and the points of a point Poisson process on $\cF$ with 
intensity $\lambda {\bf m} $ has the same 
distribution as the isometry class of $\cF_{\lambda}$.

   In the last section, Theorem \ref{poidecanc} provides  
a decomposition of $\cF$ along $\cF_{\lambda}$. In the supercritical case 
$m<0$, it is easy to see that if $\cF$ is infinite, 
then the infinite subtree of $\cF$ is simply the tree $\cF_0$ and the latter 
decomposition provides a decomposition of the L\'evy forest along its infinite 
component which is distributed as a GW($\xi_0, \psi'(\gamma),
a$)-forest. This generalizes a decomposition known for Galton-Watson
trees (see \cite{LyPe}).

  This paper is organized as follows: in Section \ref{disctrees}, we set notation concerning discrete trees and we discuss the Bernoulli leaf colouring of 
discrete Galton-Watson trees. In Section \ref{locallycompactrealtrees} 
we introduce real trees (Subsection \ref{realtrees}) and we define 
the Gromov-Hausdorff topology 
on the isometry classes of locally compact rooted real trees
(Subsection \ref{gromovtopo}); in Subsections 
\ref{leb} and \ref{decpoissonsample}
for technical reasons we shall need to embed
locally 
compact trees in the Banach space $\ell_1 (\bN)$; the way to do that is 
explained in Subsection \ref{ellunemb}. In Section \ref{growth} we define the 
growth process of the L\'evy forest: we first discuss in Subsection 
\ref{levygal} the Bernoulli leaf colouring of Galton-Watson trees with 
exponential edge lengths and in particular we prove 
Theorem \ref{prolong} that specifies the distribution of 
Bernoulli leaf colouring consistent families of Galton-Watson trees; 
Subsection \ref{percolevy} is devoted to the construction of the 
growth process; we briefly discuss the infinitesimal 
dynamics of the growth process and at the 
end of this subsection we also give a special 
probabilistic construction of the increments of the growing process that 
shall be used in the proofs of the next sections. 
Section \ref{mainlevy} is devoted to the study of the L\'evy forest: 
in Subsection \ref{levyforest} we prove the convergence result Theorem 
\ref{mainresult}; in Subsection \ref{leb} we prove Theorem 
\ref{meslebe} that concerns the mass measure; Subsection 
\ref{excursionmes} is devoted to the definition of the excursion measure;
In Subsection \ref{decpoissonsample} we discuss the decomposition of the 
L\'evy forest along the ancestral tree of the points of a Poisson
sample with intensity the mass measure.

\section{Discrete trees}
\label{disctrees}

\subsection{Basic definitions and notations.}
\label{basic}

   Let us set 
\beqs \bU=\bigcup_{n=0}^\infty\bN^{*\,  n}
  \eeqs
where $\bN^*=\{1,2,\ldots\}$ and by convention $\bN^{* \, 0}=\{\emptyset\}$. 
The concatenation of words in $\bU$ is denoted $w=vu=(v_1,\ldots,v_m,u_1,\ldots,u_n)$
for $v=(v_1,\ldots,v_m),u=(u_1,\ldots,u_n)\in \bU$. 
Following Neveu \cite{Ne} we represent an ordered rooted tree as a subset
$t\subset \bU$ satisfying 
\begin{itemize}
  \item $\emptyset\in t$; $\emptyset$ is called the \em ancestor \em of
    $t$.
  \item $j\in\bN,vj\in t\Rightarrow v\in t$; $v$ is called the \em
  parent \em of $vj$.
  
\item for any $v\in t$, there exists an integer 
$ k_v(t)$ such that $vj \in t$ , $1\leq j\leq  k_v(t)$. $k_v(t)$ is the \em number of children \em of $v$. 

\end{itemize}

We denote by $\bT_{{ {\rm discr}}}$ the space of all discrete 
ordered rooted trees. 
On each $t\in\bT_{{ {\rm discr}}}$, we have the \em genealogical  order 
\em given by 
  \beqs v\preceq w\quad \Longleftrightarrow\quad vu=w\quad\mbox{for some $u\in \bU$.}
  \eeqs
Any tree $t\in \bT_{{{\rm discr}}}$ is also totally ordered by the \em
lexicographical order \em on $\bU$
denoted by $\leq $. Note that if $t$ is infinite, then 
$(\bT_{ {\rm discr}}, \le)$ cannot in general 
be embedded in $(\bN,\le)$ in an order-preserving way.

  Let $u\in t$. We say that $u$ is \em a leaf of $t$ \em iff $k_u(t)=0$. 
We denote by ${\rm Lf}(t)$ the set of leaves of $t$. Note that ${\rm Lf}(t)$ is possibly 
empty. We define the shifted tree $t$ at $u$ by 
$$ \theta_u t= \{ v\in \bU \; : \; uv\in t \} \, .$$ 
Then $\theta_u t=\emptyset$ iff $u\in {\rm Lf}(t)$. Let $v\in t$. We denote by 
$\lgeo u,v\rgeo$ the shortest path between $u$ and $v$ and by $u\wedge v$ the
last common ancestor (or branching point) of $u$ and $v$. We set 
$ \rgeo u,v\rgeo := \lgeo u,v\rgeo \backslash \{ u\} $ and we define similarly  
$\lgeo u,v\lgeo$ and $\rgeo u,v\lgeo$.

  We endow $\bT_{{ {\rm discr}}}$ with the $\sigma$-algebra $\cG_{{{\rm discr}}}$ 
generated by the countable
family of subsets 
$\{ t\in \bT_{{ {\rm discr}}}: \; u\in t\}$ , $u\in \bU$. Unless 
otherwise specified, the
random variables that we consider in this paper are defined on the
same probability space $(\Omega , \cA, \bP)$ which is assumed to be
large enough to carry as many independent random variables as we
require. Let $\xi $ be a probability distribution on $\bN$. 
We call {\it Galton-Watson tree} with offspring distribution 
$\xi$ (a GW($\xi $)-tree for short) any 
$\cG_{{{\rm discr}}}$-measurable random variable $\tau $ 
whose distribution is characterized by the two following 
conditions: 
\begin{description}
\item{(i)} $\bP (k_{\emptyset }(\tau )=i) = \xi (i) \; ,\; i\geq 0 $.

\item{(ii)} For every $i\geq 1$ such that $\xi (i) \neq 0$, the shifted trees $\theta_{1} (\tau) , \ldots , \theta_{i} (\tau) $  
under \\ 
$\bP (\cdot \mid k_{\emptyset } (\tau )=i)$ are independent copies of 
$\tau$ under $\bP$.
\end{description}

   We shall sometimes consider finite sequences of discrete trees $f=(t_1,
\ldots ,t_n)$. We call them {\it forests of discrete trees} and we denote their
set by $\bF_{{ {\rm discr}}}$. The elements of the forest are ordered 
by putting first the vertices of the first tree, next the vertices of the
second tree etc. The genealogical order on a forest 
is defined tree by tree. 
A Galton-Watson forest with $n$ elements is just a sequence of $n$
i.i.d. GW-discrete trees.

\subsection{Bernoulli leaf colouring of Galton-Watson trees. }
\label{secdiscrete}

   In this section we discuss Bernoulli colouring of the leaves of a GW-tree
and we compute the distribution of the whole tree conditionally on the genealogy 
of the leaves remaining after the colouring. More precisely, 
let $p\in [0,1]$ and let $\tau $ be a GW($\xi $)-tree. We assume that 
$\tau$ has a.s. leaves which is obviously equivalent to the condition
\begin{equation}
\label{nontrivial}
\xi (0) >0,
\end{equation}
and it will be convenient to assume $\xi(1)=0$.

We colour independently at random
each leaf of $\tau$ in red with probability 
$p$ and in black with probability $1-p$. If there is at least one
black leaf, we also colour in black the subtree spanned by the root and
the black leaves, namely the ancestral tree of the black leaves; then, 
we colour in red the remaining vertices. If there is no black leaf, we
colour all the tree in red.

\begin{figure}[ht]
\psfrag{taugw}{$\tau$}
\psfrag{tausub}{$\tau_{{\rm sub}}$}
\psfrag{taublack}{$\tau_{{\rm b}}$}
\psfrag{tau1}{$\tau_1 $}
\psfrag{tau2}{$\tau_2 $}
\psfrag{tau3}{$\tau_3 $}

\epsfxsize=5cm
\centerline{\epsfbox{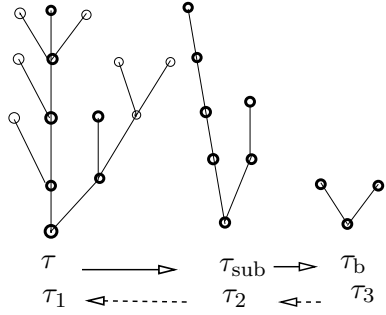}}
\caption{{\small The black vertices are represented by the thick circles 
and the red ones by the thin circles. The dashed arrows represent the
   reconstruction procedure.}}
\label{discperc}
\end{figure}

   Assume that $\tau$ is not completely red. Then, the 
black subtree is isomorphic to a random 
tree in $\bT_{ {\rm discr}}$ denoted by $\tau_{{ {\rm sub}}}$ and
also called the {\it black subtree}.  
The {\it black tree} (which is distinct from the black subtree) 
is obtained as follows: define a graph with set of
vertices $V$ and set of edges $E$ given by
$$ V=\{ \emptyset  \}\cap\{ u\in \tau_{{ {\rm sub}}} \, :\; 
k_u (\tau_{{ {\rm sub}}})
\neq 1 \} $$
and 
$$ E=\left\{ \{ u,v\}\, ;  \, u,v \in \tau_{{ {\rm sub}}}
\; :\; u\neq v \; {\rm and }\; V\cap \rgeo u,v\lgeo =\varnothing \right\}.$$
Here $\rgeo u,v\lgeo $ is the shortest path between $u$ and $v$ in 
$\tau_{{ {\rm sub}}}$. Put on $V$ the order inherited from $\tau_{{
    {\rm sub}}}$. Then $(V,E)$, 
with the distinguished vertex $\emptyset$, is an ordered rooted tree
isomorphic to a unique element $\tau_{{\rm b}}$ in $\bT_{{ {\rm discr}}}$
that is taken as the definition of the {\it black tree} (see Figure 
\ref{discperc}).

The main goal of this 
section is to give the joint distribution of $\tau_{{\rm b}}$, 
$\tau_{{{\rm sub}}}$ and $\tau $ in terms of $\xi$ and $p$. 
More precisely, let us define the colour of each vertex 
$u\in \tau$ as the mark $c_u \in \{ 0,1\}$: $c_u=1$ if $u$ is black and $c_u=0$ if it is
red. The two-colours tree is the $\{0,1\}$-marked tree
$\tilde{\tau}=(\tau ; c_u,u\in\tau )$ distributed as follows:

\begin{itemize}\item Conditionally on $\tau$, the random variables are 
$\{ c_u,u\in {\rm Lf}(\tau )\}$  i.i.d. Bernoulli random variables with expectation $1-p$.

\item For any $v\in \tau $, we set $c_v=1$ if there is 
$u\in {\rm Lf}(\tau)$ such that $v\preceq u $ and $c_u=1$; set $c_v=0$ otherwise.

\end{itemize}

Let $u\in \tau$. We denote by 
$k_u^r (\tilde{\tau})$ the number of red children of $u$ and by 
$k_u^b (\tilde{\tau})$ the number of black ones. Let $l\geq 0$ and $\varepsilon=
(\varepsilon_1, \cdots, \varepsilon_l) \in \{ 0,1\}^l$. Denote by $l_r$ the
number of $0$ in $\varepsilon$ and by $l_b$ the number of $1$. Let 
$f_1, \cdots ,f_l$ be $l$ nonnegative measurable functions on the set of
two-coloured discrete trees equipped with the smallest $\sigma$-field
making the marks measurable. Then 
it is easy to show that 
\begin{equation}
\label{perco1}
 \bE \left[\prod_{i=1}^l f_i(\theta_i \tilde{\tau} ) \; ;\; k_{\emptyset} (\tau)=l \; ;\; 
(c_1, \cdots, c_l)=\varepsilon \right]
 =\xi (l) g(p)^{l_r} (1-g(p))^{l_b} \prod_{i=1}^l 
\bE \left[ f_i( \tilde{\tau} ) 
\left| c_{\emptyset}=\varepsilon_i \right. \right],
\end{equation}
where $g(p)=\bE [p^{\# L(\tau )} ]$. Here  $\theta_j \tilde{\tau}$
stands for the marked tree shifted at the $j$-th children of the
ancestor:
$$ \theta_j \tilde{\tau}= (\theta_j \tau ; c_{ju} \, , \, u\in
\theta_j \tau ). $$
Let us denote by $\varphi$ the generating
function of $\xi$: 
$$   \varphi (s)=\sum_{k\geq 0} \xi (k) s^k \; , \quad s\in [0,1] .$$
By splitting $\tau$ at the root, we prove that $g$ satisfies
\begin{equation}
\label{equag}
g(s)=\varphi(g(s))-\xi(0)(1-s) \; , \quad s\in [0,1]  .
\end{equation}
\noindent
Formula (\ref{perco1}) implies that $ \tilde{\tau}$ is a two-types
Galton-Watson tree whose branching mechanism is described as follows: 

\begin{enumerate}\item[(a)] The tree $ \tilde{\tau}$ is completely red iff  
$c_{\emptyset}= 0$ which happens with probability $g(p)$.  
The tree conditioned to have no black vertices is 
a (completely red) GW($\xi_r$)-tree where $\xi_r$ is given by 
\begin{equation*}
\xi_r (l)= \left\{   
\begin{array}{ll}
\displaystyle 
\xi(l)g(p)^{l-1} &{\rm if} \;\displaystyle 
l\geq 1\; ; \\
\displaystyle \xi(0)p/g(p)
&{\rm if} \displaystyle \; 
l=0 .
\end{array} 
\right. 
\end{equation*}
Then, the generating function $\varphi_r$ of $\xi_r$ is given by  
\begin{equation}
\label{genfunred}
\varphi_r(s)=  1-\frac{\varphi(g(p))-\varphi(g(p)s)}{g(p)}\; ,\; s\in [0,1].
\end{equation}

\item[(b)] Conditionally on having at least one black leaf, 
the two-types offspring
distribution is given by 
\begin{multline*}
\hspace{-0.5cm}
\bP \left((k_{\emptyset}^r (\tilde{\tau}),k_{\emptyset}^b (\tilde{\tau}))=(l_r,l_b)
| c_{\emptyset}=1\right)=  \\
\left\{   
\begin{array}{ll}
\displaystyle 
\xi(l)g(p)^{l_r}(1-g(p))^{l_b-1} \frac{(l_b+l_r)!}{l_b!l_r! } &{\rm if}\; 
\displaystyle 
l_b\geq 1, l_r\geq 0 ; \\
\displaystyle \xi(0)\frac{1-p}{1-g(p)}
&{\rm if} \displaystyle \; 
l_b=l_r=0 .
\end{array} 
\right. 
\end{multline*}

\item[(c)] Conditionally on $\{ k_{\emptyset}^r (\tilde{\tau})=l_r;
k_{\emptyset}^b (\tilde{\tau})=l_b  \}$, $(c_1, \cdots , c_{l_r+l_b})$ is
uniformly distributed among the $(l_b+l_r)!/l_r ! l_b !$ possibilities.

\item[(d)] Conditionally on $\{ (c_1, \cdots , c_{l_r+l_b})=
\varepsilon \}$, $\varepsilon \in \{ 0,1\}^l$, 
the marked trees $\theta_1 \tilde{\tau}, \cdots , \theta_l \tilde{\tau}$ are 
independent and $\theta_i \tilde{\tau}$ has the same distribution as 
$\tilde{\tau}$ under $\bP (\; \cdot \; |  c_{\emptyset}=\varepsilon_i)$. 
\end{enumerate}

   Before giving the joint law of $\tau_{{\rm b}}$ , $\tau_{{\rm sub}}$ and
   $\tau$, 
we need to introduce some notation: we first define the ``black'' offspring 
distribution $\xi_b$ by 
\begin{equation*}
\xi_b(l)= \left\{   
\begin{array}{lll}
\displaystyle 
\frac{\varphi^{(l)}(g(p)) (1-g(p))^{l-1}}{1-\varphi'(g(p))}&{\rm if} 
\;\displaystyle 
l\geq 2 \; ; \\
\displaystyle \quad \quad 0 \quad 
&{\rm if} \displaystyle \; 
l=1 ;
\\
\displaystyle \frac{\xi(0)(1-p)}{(1-g(p))(1-\varphi'(g(p)))}
&{\rm if} \displaystyle \; 
l=0 .
\end{array} 
\right. 
\end{equation*}

\noindent Then, its generating function is given by 
\begin{equation}
\label{genfunblack}
\varphi_b(s)= s+\frac{\varphi(g(p)+s(1-g(p))) -g(p)-s(1-g(p))}{(1-g(p))  
(1-\varphi'(g(p)))}\; ,\quad s\in [0,1].
\end{equation}
\noindent  We also introduce for any $l\geq 1$ the following 
probability distribution on $\bN$: 
 \begin{equation}
\label{addvert}
\nu_{l}(k)=\xi(l+k)\frac{(l+k)!g(p)^k}{k!\varphi^{(l)}(g(p))}\;
,\quad k\geq 0 .
\end{equation}
\noindent 
The joint law of the black tree and the red forest is given by the
following: 

\vspace{3mm}

\noindent  

{\bf Reconstruction procedure for discrete GW-trees. }

\begin{itemize}
\item{{\bf Step 1}}: Let $\tau_1$ be a GW($\xi_b$)-tree. 
For any $u\in \tau_1$ distinct from the root, insert a line-tree 
with a random number $N_u$ of edges at the end of the edge between 
$u$ and its parent, and graft directly 
on the root a line-tree with a random number $N_{\emptyset}$ of edges.   
The $N_{u}$, $u\in \tau_1 $ are distributed as follows: 
conditionally on $\tau_1$, they are i.i.d. random variables  
with a geometric distribution given by
\begin{equation}
\label{attach}
 \bP (N_u =k| u\in \tau_1)= (1-\varphi'(g(p)))\varphi'(g(p))^k \; ,\quad k\geq
 0.
\end{equation}
\noindent
The resulting random element in $\bT_{{ {\rm discr}}}$ is denoted by $\tau_2$ and has the
same distribution as the black subtree.

\item{{\bf Step 2}}: Independently, on each vertex $u\in \tau_2$ such that 
$k_u( \tau_2) =l >0$ graft a random number with distribution $\nu_l (\cdot)$
of red vertices. Insert these
new red vertices uniformly at random among the $l$ black ones. 
Then, graft independently on each newly added red vertex 
an independent GW($\xi_r$)-tree. We obtain 
a two-colours tree denoted by $\widetilde{\tau}_3$. 

\end{itemize}

\noindent We get the following identity: 
\begin{equation}
\label{identdiscr}
\left( \tau_b , \tau_{{ {\rm sub}}}, \widetilde{\tau}  \right) 
\quad {\rm under} \quad 
\bP \left( \; \cdot \; | c_{\emptyset}=1 \right) \; 
\overset{{\rm (dist)}}{=} \; 
\left(\tau_1 , 
\tau_2 , \widetilde{\tau}_3 \right). 
\end{equation}
This identity is a consequence of an elementary computation based on 
(a), (b), (c) and (d), and it is left to the reader. Note that 
(\ref{identdiscr}) implies in particular 
that the black tree $\tau_b$ is distributed 
as a GW($\xi_b $)-tree.

  Denote by $N$ the number of red trees grafted on the black subtree 
of $\widetilde{\tau}$ if $\widetilde{\tau}$ is not completely red and 
set $N=1$ if $\widetilde{\tau}$ is completely red. 
Denote by $ \kappa $ the generating function of $N$: $ \kappa (s):=
\bE [ s^{N} ] . $
By splitting $\widetilde{\tau}$ at the root and by an elementary computation
based on (\ref{perco1}), we show that $ \kappa (s)$ satisfies the following equation
\begin{equation}
\label{numbredtrees}
\varphi \left( \kappa (s)\right)- \kappa (s)= 
\varphi \left( sg(p) \right)- sg(p) - \left( 
\varphi \left( g(p)\right) - g(p) 
\right). 
\end{equation}
We shall use this identity in Section \ref{levyforest}.

\section{The space of locally compact rooted real trees.}
\label{locallycompactrealtrees}
\subsection{Real trees.}
\label{realtrees}

   Real trees form a class of loop-free length spaces, which turn out to be the class of limiting objects 
of many combinatorial and discrete trees, extending the class of trees with 
edge lengths. More precisely we say that a metric space $(T,d,\rho)$ is 
{\it a rooted real tree} if it satisfies the following conditions: 
\begin{itemize} 
  \item  For all $s,t\in T$, there is a 
    unique isometry $f_{s,t}:[0,d(s,t)]\rightarrow T$ such that 
    $f_{s,t}(0)=s$ and $f_{s,t}(d(s,t))=t$; 
  \item If $q$ is a continuous injective map from $[0,1]$ into $T$, we have
    \beqs q([0,1])= f_{q(0), q(1)} ([0, d(q(0), q(1) )])
    \eeqs
  \item $\rho\in T$ is a distinguished point, called the \em root\em.
\end{itemize} 
Let us introduce some notation: we denote by $\lgeo s,t\rgeo $ the trace of $f_{s,t}$:
$\lgeo s,t\rgeo :=f_{s,t}([0, d(s,t)]) $. We also denote by 
$\rgeo s,t\rgeo $, $\lgeo s,t \lgeo$ and 
$\rgeo s,t\lgeo$ the respective images of $(0,d(s, t)] $, $[0,d(s, t)) $ and 
$(0,d(s, t)) $ by $f_{s,t}$. There is a nice characterization 
of real-trees that we use in the next subsection 
which is called the {\it four points condition}: 
let $(X,d)$ be a
complete path-connected metric space; then it is a real tree iff 
\begin{equation}
\label{fourpoint}
d(s_1, s_2) +d(s_3, s_4) \leq 
(d(s_1, s_3) +d(s_2, s_4))\vee (d(s_3, s_2) +d(s_1, s_4)). 
\end{equation}
We refer to \cite{Dress84,DMT96,DT96} for general
results concerning real trees, \cite{Pau88, Pau89} for applications of real 
trees to group theory and to \cite{Ev00, EvPitWin, EvWin},\cite{ DuLG} 
and also \cite{LyoHam} for a probabilistic use of real 
trees.

   In this paper we restrict our attention to 
{\it locally compact real trees}. By the Hopf-Rinow theorem (see for
instance \cite{Gro}, Chapter 1)
the closed balls are 
compact sets. For any 
$s\in T$ we denote by ${\rm n}(s,T)$ the {\it degree} of $s$, namely the 
(possibly infinite) 
number of connected components of $T\setminus\{s\}$. For convenience
of notation, we often denote ${\rm n}(s,T)$ by ${\rm n}(s)$ when
there is no risk of confusion. We denote by 
$${\rm Lf}(T)=\{ s\in T\setminus\{ \rho \}: \; {\rm n}(s,T)=1\} \quad 
{\rm and } 
\quad {\rm Br}(T)=\{ s\in T\setminus\{ \rho \}: \; {\rm n}(s,T)\geq 3\}$$
respectively the set of the {\it leaves} of $T$ and 
the set of {\it branching points} of $T$. We also denote by 
${\rm Sk} (T)$ the {\it internal skeleton} of $T$ that is defined by ${\rm Sk} (T)= T\backslash
{\rm Lf}(T)$. We can easily prove that for any 
sequence $(s_n,n\geq 1)$ dense in $T$, we have 
\begin{equation}
\label{intskel}
{\rm Sk} (T)=\bigcup_{n\geq 1} \lgeo \rho, s_n \lgeo \; . 
\end{equation}
\noindent
Then, the closure of ${\rm Sk} (T)$ is $T$. Note that the trace on
${\rm Sk} (T)$ of the Borel $\sigma$-field is generated by the ``intervals''  $\lgeo s,s'\rgeo$, 
$s,s'\in {\rm  Sk} (T)$. Thus we can define a unique positive Borel 
measure $\ell_T (ds)$ on $T$ such that 
$$ \ell_T (L(T))=0 \quad {\rm and } \quad \ell_T (\lgeo s,s'\rgeo)= d(s,s') .$$
The measure $\ell_T$ is usually called the {\it length measure} of $T$. 
We next prove the following simple lemma. 

\begin{lemma}
\label{branching}
The set of branching points of a locally compact real tree is at most
countable. 
\end{lemma}

\noindent
{\bf Proof}:  Let $(T,d)$ 
be a locally compact real tree. Assume that ${\rm Br}(T)$ is 
uncountably infinite. 
Since $ {\rm Br}(T)\subset {\rm Sk} (T)$, by 
(\ref{intskel}) there is a positive integer $n$ 
such that the set $\lgeo \rho, s_n \rgeo \cap 
{\rm Br}(T)$ is uncountable. Thus we can find 
an injective map ${\rm j}$ from $\bR$ into $\lgeo \rho, s_n \rgeo \cap 
{\rm Br}(T)$. Then with any 
$x\in \bR$, we can associate a connected component 
$C_x$ of $T\setminus{\{ {\rm j}(x)\}}$ such that 
$\lgeo \rho, s_n \rgeo \cap C_x=\varnothing $ since ${\rm j}(x)$ is a
branching point. A simple argument implies that 
$ C_x \cap C_y=\varnothing$ for any $x\neq y$ and $T$ cannot be 
separable, which contradicts the fact that it is locally compact. \cqfd

\vspace{5mm}

In the present paper we define step by step a growing family of trees by
recursively grafting independent random trees on
 nodes and branches of the tree of the previous step. Let us explain in
the deterministic setting one step of this grafting procedure: 
let $(T, d,
\rho)$ be a locally compact real rooted tree, let 
$(T_i, d_i, \rho_i)$, $i \in I$, be 
a family of locally compact real trees and let $(s_i, i\in I)$ be a
collection of vertices of $T$. We specify $T^\prime$ as disjoint union 
$$ T' =T \coprod_{i\in I} T_i\setminus{\{ \rho_i\}}  $$
and we define a distance $d'$ on $T'\times T'$ as follows: 
$d'$ coincides with $d$ on $T\times T$ and if 
$s\in  T_i\setminus{\{ \rho_i\}}$ and $s'\in T'$, then we set 
\begin{equation*}
d'(s,s')= \left\{   
\begin{array}{lll}
\displaystyle 
d_i(s,\rho_i) + d(s_i,s')  &{\rm if} \;\displaystyle 
s'\in T; \\
\displaystyle d_i(s,\rho_i)+d(s_i,s_j)+d_j(s',\rho_j)
&{\rm if} \displaystyle \; 
s'\in T_j\setminus{\{ \rho_j\}}, \; i\neq j; \\
\displaystyle
d_i(s,s')&{\rm if} \displaystyle \; 
s'\in  T_i\setminus{\{ \rho_i\}}.
\end{array} 
\right. 
\end{equation*}
It is easy to prove that $(T',d',\rho'=\rho)$ is a real tree and we use the 
notation 
$$ (T, d',\rho')= T \circledast_{i\in I} (s_i, T_i) $$
to mean that $ (T', d',\rho')$ is obtained from $(T, d, \rho)$ by this
``grafting''  procedure.

\subsection{Gromov-Hausdorff convergence of pointed metric spaces.}
\label{gromovtopo}
  
The purpose of this section is to introduce a nice topology 
on the set $\bT$ of {\it isometry classes} of locally 
compact rooted real trees:
more precisely, we say that two pointed metric spaces 
$(X_1,d_1,\rho_1)$ and $(X_2,d_2,\rho_2)$ 
are equivalent iff there exists an isometry 
$f$ from $X_1$ onto $X_2$ such that $f(\rho_1)=\rho_2$. 
Evans, Pitman and Winter \cite{EvPitWin} 
showed that the set $\bT_{{\rm cpct}}$ of isometry  
classes of 
\em compact \em rooted real trees equipped with the Gromov-Hausdorff 
distance whose definition is recalled below, is a complete and separable 
metric space. Here we define a metric on $\bT$ and 
prove a similar result in the the locally compact 
case. For sake of clarity, we actually prove this result for locally
compact length spaces, the real tree case being a simple consequence
of the four points conditions (\ref{fourpoint}) 
that characterizes real trees.

Let us first recall the definition of the Gromov-Hausdorff distance of 
two pointed compact metric spaces 
$(X_1,d_1,\rho_1)$ and $(X_2,d_2,\rho_2)$: we set 
$$ \delta_{{\rm cpct}} (X_1,X_2)= \inf 
\left\{ d_{{\rm Haus}}(f_1(X_1),f_2(X_2))\vee d (f_1(\rho_1),
f_2(\rho_2)) \right\}$$
where the infimum is taken over all isometric 
embeddings $f_i:X_i \rightarrow E$, $i=1,2$  into a common metric space 
$(E,d)$. Here $d_{{\rm Haus}}$ stands for the Hausdorff distance 
on the set of 
compact sets of $E$. Observe that $ \delta_{{\rm cpct}}$ 
only depends on the isometry classes of the $X_i$'s 
and we can show that it defines a metric on the set of isometry classes 
of all pointed compact metric spaces (see \cite{Gro}).

  There is a useful way to control $ \delta_{{\rm cpct}} (X_1,X_2)$ via 
$\varepsilon $-isometries. Namely, we say that a (possibly not continuous) map
$f:X_1\rightarrow X_2$ is {\it a pointed $\varepsilon $-isometry} if 

\begin{enumerate}
  \item[(i)] $f(\rho_1)=\rho_2$ 
  \item[(ii)] $ {\rm dis} (f):=\sup \{| d_1(x,y)-d_2(f(x),f(y))|\, ;\;  
  x,y\in X_1 \} \; <\; \varepsilon ;$
  \item[(iii)] $f(X_1)$ is an $\varepsilon $-net of $X_2$. 
\end{enumerate}
The quantity ${\rm dis} (f)$ is called the {\it distortion} of $f$. 
The following lemma is a straightforward consequence of the non-pointed case
stated in Corollary 7.3.28 in \cite{BuBu}.  
\begin{lemma}
\label{epsiso} 
Let 
$(X_1,d_1,\rho_1)$ and $(X_2,d_2,\rho_2)$ be two pointed compact metric spaces. Then, 
 
\begin{enumerate}
  \item[{\rm (a)}] If $ \delta_{{\rm cpct}} (X_1,X_2) < \varepsilon $, then there exists a 
$4\varepsilon$-isometry from $X_1$ to $X_2$. 

\item[{\rm (b)}] If there exists a $\varepsilon$-isometry 
from $X_1$ to $X_2$, then $ \delta_{{\rm cpct}} (X_1,X_2) < 4\varepsilon $.
\end{enumerate}
\end{lemma}
\noindent

Let us now recall from \cite{BuBu}, 
Chapter 8,  a way to extend the Gromov-Hausdorff convergence to
non-compact metric spaces. Let  
$(X,d)$ be a metric space. Let $r\geq 0$ and $\rho\in X$. We denote 
by $B_{X}(\rho,r)$ the closed ball centered at $\rho$ with radius $r$. Let 
$(X_n,d_n,\rho_n)$, $n\geq 1$, be a sequence of pointed metric spaces; we say 
that this sequence converges in the pointed Gromov-Hausdorff sense to 
the pointed metric space $(X,d,\rho)$ if for any $r, \varepsilon >0$ there 
exists $n_0=n_0 (r, \epsilon) \geq 1$ such that for every $n\geq n_0$,
there is a map 
$f_n: B_{X_n}(\rho_n,r) \rightarrow X$ satisfying the following conditions:

\begin{enumerate}
\item[{\rm(i')}]  $f_n(\rho_n)=\rho;$
\item[{\rm(ii')}] ${\rm dis} (f_n) <\varepsilon $;
\item[{\rm(iii')}] 
The $\varepsilon$-neighbourhood of $f_n(B_{X_n}(\rho_n,r))$ contains 
$B_{X}(\rho,r-\varepsilon) $.
\end{enumerate}
We use the following notation: 
\begin{equation}
\label{limitformu}
(X_n,d_n,\rho_n) \build{\longrightarrow}_{n\to\infty}^{{\rm G-H}}
(X,d,\rho).
\end{equation}
Let us briefly recall from \cite{BuBu} useful 
properties of pointed Gromov-Hausdorff
convergence. Assume that (\ref{limitformu}) holds. Then, 
\begin{enumerate}
\item[(a)]{(8.1.8 and 8.1.9  \cite{BuBu})} 
If the $X_n$'s are locally compact length spaces and if $X$ 
is complete, then $X$ is a locally compact length space.
   
\item[(b)]{(8.1.2 \cite{BuBu})} If the $X_n$'s are compact and if $X$
  is compact, then 
$$\delta_{{\rm cpct}}(X_n,X )\build{\longrightarrow}_{n\to\infty}^{\;
} 0. $$ 

\item[(c)]{(8.1.3 \cite{BuBu})} If $X$ is a length space, then for any $r>0$
 $$\delta_{{\rm cpct}}(B_{X_n}(\rho ,r), B_{X}(\rho ,r)
 )\build{\longrightarrow}_{n\to\infty}^{\; } 0.$$

\item[(d)]{(8.1.9 \cite{BuBu})} ({\it Pre-compactness}): Let
  $\cC$ be a set of pointed metric spaces. Assume that for any 
$r,\varepsilon >0$, there exists $N(r,\varepsilon)$ such that for every
$(X,d,\rho)\in \cC$ the closed ball $B_X(\rho ,r)$ admits an $\varepsilon$-net 
with at most $N(r,\varepsilon)$ points. Then, any sequence of elements of
$\cC$ contains a converging subsequence in the pointed Gromov-Hausdorff
sense. 
\end{enumerate}

For locally compact length spaces, the pointed 
Gromov-Hausdorff convergence is compatible with the following metric: 
let $(X_1,d_1,\rho_1)$ and $(X_2,d_2,\rho_2)$ be two pointed  
locally compact length spaces;
under our assumptions $(B_{X_i}(\rho_i,r), d_i,\rho_i)$ is a pointed
compact space so it makes sense to define  
\begin{equation}
\label{defgrohau}
\delta (X_1,X_2)= \sum_{k\geq 1} 2^{-k} 
\delta_{{\rm cpct}} 
(B_{X_1}(\rho_1,k),B_{X_2}(\rho_2,k)) .
\end{equation}
Clearly, $\delta$ only depends on the isometry classes of $X_1$ and $X_2$. 
Let us denote by $\bX$ the set of isometry classes of pointed
locally compact length spaces and by $\bX_{{\rm cpct}}$ the set of
pointed compact length spaces. 
\begin{proposition}
\label{polish} Let $(X_n,d_n,\rho_n)$, 
$n\geq 1$ and $(X,d,\rho)$ be representatives
of elements in $\bX$. Then 
$$ (X_n,d_n,\rho_n) \build{\longrightarrow}_{n\to\infty}^{{\rm G-H}} (X,d,\rho)
\quad \Longleftrightarrow \quad \lim_{n\rightarrow \infty} \delta (X_n,X)=0 . $$
Moreover $(\bX, \delta)$ is complete and separable. 
\end{proposition}
{\bf Proof:} The fact that the $\delta $-convergence implies the
pointed Gromov-Hausdorff convergence is easy to deduce from properties (b) and Lemma 
\ref{epsiso}. The
converse is a consequence of (c).

  Next, we prove that $\delta$ is a metric on $\bX$. Since $\delta_{{\rm cpct}}$ satisfies
the triangle inequality, so does $\delta$. Let $(X_1,d_1,\rho_1)$ and $(X_2,d_2,\rho_2)$
be two pointed locally compact length spaces such that 
$\delta (X_1,X_2)=0$. Then, for every $k\geq 1$ there exists an isometry $f_k$ from $B_{X_1}(\rho_1,k)$ onto 
$B_{X_2}(\rho_2,k)$ with $f_k(\rho_1)=\rho_2$. Let $(x_n,n\geq 1)$ be a dense sequence in $X_1$. 
By the Cantor diagonal procedure we can find an increasing sequence of indices
$(k_i,i\geq 1)$ such that for any $n\geq 1$, $(f_{k_i}(x_n), i\geq 1)$
converges in $X_2$. Set $f(x_n)=\lim_{i\rightarrow \infty} f_{k_i}(x_n)$: it 
defines an isometric embedding of 
$(x_n,n\geq 1)$ into $X_2$ such that $f(\rho_1)=\rho_2$ which can be easily
extended to an isometry $f$ from $X_1$ into $X_2$. 

  It remains to prove that $f(X_1)=X_2$. By exchanging the roles of $X_1$ and
  $X_2$, we get  an isometric
embedding $g$ from $X_2$ into $X_1$ such that $g(\rho_2)=\rho_1$. Then, for any 
$k\geq 1$, $f\circ g$ is an isometric map from the compact set
  $B_{X_2}(\rho_2,k)$ into itself. Thus, 
it is a bijective map and we get $f(B_{X_1}(\rho_1,k))=B_{X_2}(\rho_2,k)$ for any $k\geq 1$ which easily proves that 
$f$ is actually onto $X_2$.

  It remains to prove that $\bX$ equipped with the metric $\delta$ is
complete and separable. Since the set of isometry classes of compact metric
spaces equipped with $\delta_{{\rm cpct}}$ is separable, so is $(\bX ,
\delta )$  
for $\bX_{{\rm cpct}}$ is dense in $(\bX ,  \delta )$ by definition of $\delta $.

   We have to show that $(\bX ,\delta )$ is complete. Let $(X_n,d_n,\rho_n)$,
$n\geq 1$, be a Cauchy sequence of representatives of elements of $\bX$. 
To prove that this sequence converges, we only have to prove that it forms a
$\delta$-precompact set. Fix $r,\varepsilon>0$; choose $k> r+1$ and
$n_0\geq 1$ such that for any $n,m\geq n_0$ , $\delta
(X_n,X_m)<2^{-k}\varepsilon /12$. It implies 
\begin{equation}
\label{cauchyseq}
\delta_{{\rm cpct}}(B_{X_n}(\rho_n ,k), B_{X_{n_0}}(\rho_{n_0} ,k)) < \varepsilon /12 \; , \quad n\geq n_0 . 
\end{equation}
By Lemma \ref{epsiso}, there exists a pointed $\varepsilon /3$-isometry 
$$f_n: B_{X_{n_0}}(\rho_{n_0} ,k) \build{\longrightarrow}_{\quad}^{\quad} 
B_{X_{n}}(\rho_{n} ,k) .$$ 

Let $\{x_1, \cdots , x_N \} $ be a $\varepsilon /3$-net of $B_{X_{n_0}}(\rho_{n_0} ,k)$.
Then, for any $n\geq n_0$, the set 
$$\{f_n(x_1),  \;  \cdots   f_n(x_N) \} $$ 
is an $\varepsilon$-net of $B_{X_n}(\rho_n ,k)$ and thus, of $B_{X_n}(\rho_n ,r)$. So 
we can find $N(r, \varepsilon)=N$ such that for any $n\geq 1$ the closed ball $B_{X_n}(\rho_n ,r)$ admits an $\varepsilon$-net 
with at most $N(r,\varepsilon)$ points. The compactness criterion $(d)$
completes the proof. \cqfd

\vspace{3mm}

  Recall that $\bT$ denotes the set of isometry classes of locally compact
rooted real trees. Since the four points condition is obviously a closed condition
for $\delta $, it implies that $\bT$ is a closed subset of $\bX$ and we deduce from Proposition
\ref{polish} the following result. 
\begin{proposition}
\label{treepointed}
$(\bT, \delta)$ is a complete and separable metric space. 
\end{proposition}
Following the proof of Lemma 2.7 in \cite{EvPitWin}, we prove the following
lemma that we shall use in the next section. 
\begin{lemma}
\label{evaconv} Let $((T_n,d_n,\rho))_{n\ge 1}$ be a 
Cauchy sequence of representatives of elements of  
 $(\bT,\delta )$ such that $T_n\subset T_{n+1}$ and $d_{n+1}|_{T_n\times T_n}=d_n$,
  $n\ge 1$. Set for any $a,b\in T_n$, $ n\ge 1$:
  \beqs d(a,b)=d_n(a,b),\quad a,b\in T_n,n\ge 1.
  \eeqs
  This defines a metric on $T_{\infty}:=\bigcup_{n\ge 1}T_n$. Furthermore, all metric
  completions of $(T_{\infty},d,\rho)$ are isometric and form the limit in $(\bT,\delta
  )$ of the sequence $((T_n,d_n,\rho))_{n\ge 1}$. 
\end{lemma} 

\subsection{Galton-Watson real trees with exponential edge lengths.}
\label{lengthtrees}

  Let us consider a discrete tree with positive marks. Namely let 
$t\in \bT_{{ {\rm discr}}}$ and let $m=(m_u; u\in t)$ be a collection of
marks in $[0, \infty]$. We assume that if $m_u=\infty$ then $u$ has no child: $k_u (t)=0$. 
Such a pair $(t,m)$ is called a {\it marked tree} and the set of marked trees is denoted by 
$\bT_{{{\rm mark}}}$. We
denote by $\cG_{{ {\rm mark}}}$ the $\sigma$-algebra generated by the
events $\{ (t,m) : \, u\in t\, , m_u >a \}$, $u\in \bU$ and $a\in \bR$. 
Thinking of the marks as distances between the nodes of $t$, we can 
associate with $(t,m)$ a real tree denoted by $T(t,m)=(T, d, \rho)$ 
as follows: set $\rho =(\emptyset , 0)$ and 
$$ T= \{ \rho \} \bigcup_{ u\in t \,: \, m_u <\infty } 
\left\{ (u, s ) \, , s\in (0,m_u] \right\} 
\bigcup_{ u\in t \,: \, m_u =\infty } 
\left\{  (u, s ) \,, s\in (0,\infty) \right\} 
 $$
and we define $d$ as follows : let $\sigma =(u,s)\in T\setminus
\{\rho\}$, then we set 
$$d(\rho, \sigma )=s+\sum_{v\in \lgeo \emptyset , u \lgeo } m_v $$
where we recall notation $\lgeo \emptyset , u \lgeo=\lgeo \emptyset ,
u \rgeo \setminus \{u\}$. Let $\sigma' =(u',s')\in T\setminus
\{\rho\}$. We define  
\begin{equation*}
d(\sigma, \sigma')= \left\{   
\begin{array}{ll}
\displaystyle 
d(\rho , \sigma)+ d(\rho , \sigma') -2 \sum_{v\in \lgeo \emptyset , u
  \wedge u' \rgeo} m_v
& \quad {\rm if} \displaystyle \quad   u\wedge u' \notin \{ u, u'\}\\
\displaystyle |d(\rho , \sigma)- d(\rho , \sigma') | 
&{\rm otherwise.}
\end{array} 
\right. 
\end{equation*}

It is easy to check that $T(t,m)=(T,d,\rho)$ is a
real tree. Instead of a single tree, consider now a {\it marked forest} $(f,m)$ that is 
a finite sequence 
$(f,m)=( (t_i , m_i); 1\leq i\leq n)$  of marked trees; the set of
marked forests is denoted by $\bF_{{ {\rm mark}}}$. With a 
marked forest $(f,m)$ we associate the real tree $T(f,m)$ defined by
$$ T(t,m)=\{\rho\}  \circledast_{1\leq i\leq n} \left( \rho , T(t_i, m(i)) \right) $$
which obtained by pasting at $\rho$ the trees $T(t_i, m(i))$. 
We also denote by $\overline{T}(f,m)$ the equivalence class of $T(f,m)$ up to 
root preserving isometries. Note that $T(f,m)$ may fail to be locally compact.
For instance if $T(f,m)$ is locally 
compact if for any infinite line of descent: $u_0 \preceq \ldots \preceq u_n \preceq \ldots$ , we have 
\begin{equation}
\label{localdisctree} 
\sum_{n\geq 0} m_{u_n} =+\infty .
\end{equation}
If (\ref{localdisctree}) is satisfied then the real tree 
$T(f,m)$ that is obtained from $f$ and $m$, 
is called a 
{\it a discrete tree with edge lengths}, namely a rooted real tree 
$(T,d,\rho)$ such that 
\begin{equation}
\label{distree}
\forall r >0 \; , \quad \# B_T (\rho , r) \cap {\rm Br} (T) \; < \; \infty \quad {\rm and } \quad n(\sigma , T)<\infty \; , 
\quad \sigma \in T . 
\end{equation}

   Conversely, with each discrete tree with edge lengths $(T,d,\rho)$ we can associate a discrete forest $f\in \bF_{discr}$ and a 
set of marks $m=(m_u,u\in f)$ such that $(T,d,\rho) =T(f,m)$. One way to proceed is the following: 
we call an {\it edge} 
of $T$ the connected components of $T\backslash ({\rm Br}(T) \cup \{ \rho \})$; each edge is isometric to an interval 
of the real line (that possibly has one infinite end); by convention, the {\it left end} of an edge is the closest end to the root; 
observe that $T$ is the closure of the union of its edges by (\ref{distree});  fix an order on each group of edges 
sharing the same left end and then 
label the edges of $T$ by words written with 
integer in the following recursive way: 
\begin{itemize}
\item Each of the ${\rm n} (\rho , T)$ edges of $T$ having the root $\rho $ as a left end are labelled by the empty word 
$\varnothing$. 
\item Take a finite edge whose right end is denoted by $y \in T$. Assume that this edge is labelled by $u\in \bU$ and consider the 
edges whose left end is $y$: the $j$-th edge with respect to the fixed order is then labelled by the word $uj$. 
\end{itemize}
In this way we construct a discrete forest $f$. Consider the edge labelled by the word $u\in \bU$. There are two cases: if the edge is 
infinite, then set $m_u=\infty$; if the edge is finite, then set $m_u=d(\rho , y) -d(\rho , x)$, where $x$ and $y$ stand  for its 
resp. left and right ends.  We clearly have $T(f,m)=(T,d,\rho)$. 
Note that such a marked forest $(f,m)$ is  by no way unique. However, it is uniquely determined if we assume first that 
$f$ is {\it proper} that is $k_u (f)\neq 1$, $u\in f$, and then if we specify some order on the edges of $T$ sharing the same 
left end.

  Let $\xi$ be an offspring distribution and let $c $ be a positive
  real number. Let 
$\tau $ be a GW($\xi$)-tree and  conditionally on $\tau$, let $m=(m_u, u\in \tau)$ be i.i.d. 
exponentially distributed random variables with parameter $c $. 
The random real tree $T(\tau , m)=(\cT, d, \rho)$ is called a 
{\it Galton-Watson real tree} with parameters $(\xi , c)$
(a GW($ \xi , c$)-real tree for short). Define for any  $t\geq 0$ , $ Z_t (\cT)=\#\{v\in\cT : d(0,v)=t\}$. Then,
we can show that $(Z_t (\cT), t\geq 0)$ is a continuous-time Markov branching
process. Moreover, if we denote by $\varphi$ the generating function of $\xi $, then 
$$ \bE \left[ \exp (-\theta Z_t (\cT))\right]= \exp (- v(t, \theta)) , $$
where  $v(t, \theta)$ is the unique non-negative solution of the integral equation 
\begin{equation}
\label{inteqharris}
\int_{e^{-\theta}}^{e^{-v(t,\theta)}} \frac{dr}{\varphi (r) -r} \, =\,
c t  
\end{equation}
(see Chapter III, Section 3, p. 106 \cite{AthNey}). $\cT$ is a discrete tree with edge lengths 
(namely, $\cT$ satisfies (\ref{localdisctree})) iff 
$Z_t (\cT)$ is a.s. finite for all $t\geq 0$, which is equivalent to the following analytical condition
\begin{equation}
\label{nonexpdiscr}
\int^{1-} \frac{dr}{|\varphi (r) -r|} \, =\, \infty.  
\end{equation}
Unless otherwise specified, we assume that {\it all the GW-real trees that we consider 
in this paper satisfy}  (\ref{nonexpdiscr}).

   Define the
height of $\cT$ by 
$h(\cT): =\sup \{ d(\rho , \sigma) \; , \; \sigma \in \cT \} \in [0, \infty]$. Then observe
that $ \bP (h(\cT) \leq t) =\exp (-v(t))$, where for any $t\geq 0$ we set 
$v(t)=\lim_{\theta \rightarrow \infty} v(t,\theta)$. It satisfies
\begin{equation}
\label{heightharris}
\int_{0}^{e^{-v(t)}} \frac{dr}{\varphi (r) -r} \, =\,
c t . 
\end{equation}

We end this subsection by precisely defining the class of random
discrete trees that we shall consider: more specifically,   
let $(\tau_i; i\geq 1)$ be an i.i.d. sequence of GW($\xi$)-trees and conditionally on the 
$\tau_i$'s, let $(m_u(i),  u\in \tau_i , i\geq 1)$ be independent exponentially distributed 
random variables with parameter $c$. Fix a positive real number $a>0$ and denote by $N$ a 
Poisson random variable with expectation $a$ that is assumed to be independent
of the $m(i)$'s and of the $\tau_i$'s.  Set 
$(f,m)=(\tau_i, m(i); 1\leq i \leq N)$. The real tree $T(f,m)= (\cF, d, \rho)$ is called a 
{\it Galton-Watson real forest} with parameters $(\xi, c, a)$ (a 
GW($ \xi , c, a$)-real forest  for short).

\subsection{Isometrical embeddings of real trees in ${\bf \ell_1 } (\bN)$.}
\label{ellunemb}

   For technical reasons we shall sometimes need to 
consider specific representatives of real trees rather that isometry classes. 
Following Aldous's idea (see \cite{Al1}), we may choose to embed locally compact rooted trees in the vector space 
${\bf \ell_1 } (\bN)$ of the summable real-valued sequences equipped with the $||\cdot  ||_1$-norm. Namely, 
$$ {\bf \ell_1 } (\bN) = \left\{ x=(x_n)_{n\geq 0} \in \bR^{\bN} \; : 
\quad  ||x ||_1 :=  \sum_{n\geq 0} |x_n| < \infty \right\} \; . $$
We introduce the space 
$\bT_{{\bf \ell_1 }}$ of the subsets $T \subset  {\bf \ell_1 } (\bN) $ such that $(T, ||\cdot  ||_1 , 0)$ is a locally 
compact rooted real tree. 
Let us denote by $d_{{\rm Haus}}$ the Hausdorff distance on compact
subsets of ${\bf \ell_1 } (\bN) $. Then, for any 
$T, T'\in \btl$, define 
$$ \dl (T, T')=\sum_{k\geq 0} 2^{-k} d_{{\rm Haus}} \left( B_T(0, k) , B_{T'}(0, k)  \right) . $$
Note that 
\begin{equation}
\label{compar}
\delta (T, T') \leq \dl (T, T').
\end{equation}
\begin{proposition}
\label{ellumet}
$(\btl , \dl)$ is a Polish space.
\end{proposition}
\noindent
{\bf Proof:} It is easily proved that $(\btl , \dl)$ is a separable metric space. Let us prove it is complete. Let 
$(T_n, n\geq 0)$ be a Cauchy sequence of elements of $(\btl, \dl)$. Then, for any 
$k\geq 1$, $B_{T_n} (0, k) ,n \geq 0$  is a $d_{Haus}$-Cauchy sequence of closed subsets of ${\bf \ell_1 } (\bN) $.  
Thus, by a well-known property of Hausdorff distances, for any $k\geq 0$ there exists a closed set 
$C_k \subset {\bf \ell_1 } (\bN) $ such that 
$$\lim_{n\rightarrow \infty} d_{{\rm Haus}} \left( B_{T_n} (0, k) , C_k \right)=0 \; , $$
which implies $\lim_{n\rightarrow \infty} \delta \left( B_{T_n} (0, k)
  , C_k \right)$ by (\ref{compar}). By Theorem 
\ref{treepointed}, $(C_k ,  ||\cdot  ||_1  , 0)$ has to be a rooted compact  real tree. Moreover, for any $k'\geq k$ we have 
$C_k \subset C_{k'}$  and Property (c) in Section 
\ref{gromovtopo} implies that 
$$ B_{C_{k'}} (0, k)=C_k \; . $$
Now set $T=\bigcup_{k\geq 0} C_k$ . The previous observations easily implies that $(T, ||\cdot  ||_1 , 0)$ is a
locally compact rooted real tree and that $\lim_{n\rightarrow \infty} \dl \left( T_n, T \right)=0$, which completes the proof. 
\cqfd

\vspace{3mm}

  Let us now briefly explain how to isometrically embed a discrete tree with edge lengths $(T , d, \rho)$ in 
${\bf \ell_1 } (\bN) $. Recall from Section \ref{lengthtrees} that we can find a discrete forest 
$f\in \bF_{discr}$ and marks $m=(m_u, u\in f)$ 
such that $\overline{T}(f,m)= (T , d, \rho)$. Recall also the definition of an edge of $T$ and recall that to each 
vertex $u\in f$ corresponds an edge in $T$. We now order the vertices of $f$ as follows: order the 
roots of $f$ and put them first; then order the vertices at height $1$ and put them after the roots of $f$; 
order the vertices at height $2$ and put them next ... etc. Recall
from the previous section the definition of an edge of $T$. For any $k\geq 0$, denote 
by $I(k)$ 
the edge of $T$ corresponding to the $k$-th vertex of $f$ visited with respect to the linear order above defined and denote by 
$x_k$ the left end of $I(k)$. Clearly $x_k$ belongs to the closure of the set $ \{ \rho\} \bigcup_{j<k} I(j)$. Then, 
let us introduce for any $k\geq 0$ the sequence $e_k\in {\bf \ell_1 } (\bN) $ given by $e_k(n)=0$ if $n\neq k$ and 
$e_k(k)=1$. Let $P=(n_k , k\geq 0)$  be an $\bN$-valued increasing sequence; 
we recursively define the map $f_P$ from $T$ to ${\bf \ell_1 } (\bN) $
as follows. 
\begin{itemize}
\item $f_P (\rho)=0$;
\item For any $k\geq 1$, and any $\sigma \in I(k)$, 
$$f_P (\sigma )= f_P (x_k) + d(x_k, \sigma) e_{n_k} \; . $$
\end{itemize}

It is easy to check that $f_P$ is an isometry. Thus $(T, d , \rho)$ and 
$(f_P(T), ||\cdot  ||_1 , 0)$ are equivalent. Now we prove the following proposition. 
\begin{proposition}
\label{embedding}
Every element of $\bT$ has a representative in $\btl$. 
\end{proposition}
\noindent
{\bf Proof:} We have to prove that any locally compact rooted real tree 
$(T, d , \rho)$ can be embedded isometrically in ${\bf \ell_1 } (\bN) $. It is possible to find  a non-decreasing 
sequence of subsets $K_n $, $n\geq 0$, with no limit points and such that $K_n$ is a $2^{-n}$-net of $T$. We set 
$$ T_n = \bigcup_{\sigma \in K_n} \lgeo \rho , \sigma \rgeo \quad {\rm and } \quad  T_{\infty} = \bigcup_{n\geq 0} T_n \; .$$
Clearly, the $T_n$'s are discrete trees with edge lengths and the closure of $T_{\infty}$ is $T$. We  recursively define a map 
$f$ from $T_{\infty}$ to ${\bf \ell_1 } (\bN) $ in the following way: 
\begin{itemize}
\item Let $P_0$ and 
$P_{n,i}$ , $n\ge 0,i\geq 0$ be disjoint subsets of $\bN$. We consider $f_{P_0}$, the isometrical embedding of $T_0$ 
into ${\bf \ell_1 } (\bN) $ as defined above and we require that $f$ coincides with $f_{P_0}$ on $T_0$. 
\item Assume that $f$ is defined on $T_n$; Denote by $T^{o}_{n,i}$ , $i\in I_n$, the connected components 
of $T_{n+1}\backslash T_n$. Denote by $\rho_{n,i}$ the closest point to the root of the closure 
$T_{n,i}$ of  $T^{o}_{n,i}$. Then, $\rho_{n,i}\in T_n$  and the $(T_{n,i}, d, \rho_{n,i})$ 
are rooted discrete trees with edge lengths. We assume for convenience of notations that the sets of indices 
$I_n$ are subsets of $\bN$. Then, for any 
$\sigma \in T_{n,i}$, we set 
$$ f(\sigma) = f(\rho_{n,i}) + f_{P_{n,i}} (\sigma) \; , $$
where $f_{P_{n,i}} $ stands for the above defined isometrical embedding of $ T_{n,i}$ in ${\bf \ell_1 } (\bN) $. 
\end{itemize}

  Thus, $f$ is an isometrical embedding of $T_{\infty}$ into ${\bf \ell_1 } (\bN) $, which has a unique extension to the 
closure $T$ of $T_{\infty}$. This completes the proof of the proposition. \cqfd

\section{The growth process.}
\label{growth}
\subsection{Bernoulli colouring of the leaves and extensibility of GW-real trees.}
\label{levygal}

  In this section we 
discuss the Bernoulli colouring of the leaves of GW-real trees and forests. In particular, we introduce 
the class of 
L\'evy GW-real trees that is, roughly speaking, the class of  GW-real
trees consistent under Bernoulli colouring. 
More precisely, let $T$ be a discrete tree with edge lengths, that is
a  rooted real tree 
satisfying (\ref{distree}). Let 
$p\in [0, 1]$. Then, colour independently each leaf of $T$ in black with probability $1-p$ and in red with probability $p$. 
Denote by $A$ the set of the black leaves. If $A$ is non-empty, then colour in black the following subtree:
$$ T_{black}=\bigcup_{\sigma \in A}  \lgeo \rho , \sigma \rgeo ; $$
Then, colour in red the remaining part $T \backslash T_{black}$ of the tree. If $A$ is empty, then colour in red the whole tree 
$T$ and set  $T_{black}= \{\rho \}$. As in the discrete case such a colouring is called a {\it p-Bernoulli leaf colouring of T}. 

\begin{remark}
\label{lapalisse}
Observe that if $T$ has leaves and if it 
is not reduced to a point, then the black subtree is reduced to the root iff $T$ is completely red.
\end{remark}

Let $\xi $ be an offspring distribution on $\bN$ such that $\xi (0) >0$. Let us assume that $\xi $ is proper, namely 
$\xi (1)=0$. Fix two positive real numbers $a,c >0$ and denote by $\cT$ (resp. $\cF $) a GW($\xi, c$)-real tree  (resp. 
a GW($\xi, c,a$)-real forest). Let $ p\in (0, 1)$. Denote by $\cT_{black}$ (resp. $\cF_{black}$)
the black subtree of $\cT$ (resp. $\cF $) resulting from a $p$-Bernoulli leaf colouring (here the extra random variables 
used for the Bernoulli colourings are chosen independent of $\cT$ and of $\cF $). 
As in the reconstruction procedure discussed in Section \ref{secdiscrete}, we first compute the distribution of $\cT$ (resp. $\cF $)
conditionally on $\cT_{black}$ (resp. $\cF_{black}$). To that end, recall the notation $g(p), \xi_b, \xi_r, \nu_l $ 
from Section \ref{secdiscrete}. Let  $\cT'$ (resp. $\cF' $) be a GW($\xi_b,(1-\varphi' (g(p))) c $)-real tree  (resp. 
a GW($\xi_b,(1-\varphi' (g(p))) c,(1-g(p))a$)-real forest). Let $\cP =\{ \sigma_i ; i\in I \}$ 
be a Poisson point process on $\cT'$ (resp. on $\cF' $) with intensity
$\varphi' (g(p)) c \, \ell_{\cT'}$ 
(resp. $\varphi' (g(p)) c \, \ell_{\cF'}$). 

\vspace{3mm}

\noindent
{\bf Reconstruction  procedure on GW-real trees.}
\nopagebreak\begin{itemize}\nopagebreak
\item For $\cT'$: on  each vertex $\sigma\in \cP \cup {\rm Br}(\cT')$  graft independently 
a random number $N_{\sigma}$ of independent GW($\xi_r , c$)-real trees; conditionally on  
$\cP \cup {\rm Br}(\cT')$ the $N_{\sigma}$'s are independent and the conditional distribution of $N_{\sigma}$ is $\nu_l$ where 
$l={\rm n}(\sigma , \cT')-1$. Denote by $\cT''$ the resulting tree.

\item For $\cF'$: do the same thing as for $\cT'$ and graft on the root $N_{\rho}$ additional 
independent GW($\xi_r , c$)-real trees, where $N_{\rho}$ stands for an independent Poisson random variable with parameter $ag(p)$. 
Denote by $\cF''$ the resulting tree. 

\end{itemize}

\begin{lemma}
\label{percola}
Assume that (\ref{nonexpdiscr}) 
holds. Then, 
$$ \left( \overline{\cT}_{black} \; ,\;  \overline{\cT}
\right) \; \; {\rm under} \; \; \bP \left( \; \cdot \; | 
\cT_{black} \neq \{ \rho \} \right) \quad 
\overset{(d)}{=} \quad \left( \overline{\cT'},  \overline{\cT''} \right)$$
and 
$$ \left( \overline{\cF}_{black} \; ,\;  \overline{\cF}
\right) \quad 
\overset{(d)}{=} \quad \left( \overline{\cF'},  \overline{\cF''} \right). $$
\end{lemma}

\noindent 
{\bf Proof}: Recall from Section \ref{lengthtrees} the definition of an edge of a discrete tree 
with edge lengths. Let us assume that $\cT= T(\tau, m)$ where $\tau $
is a GW($\xi$)-tree and 
where $m=(m_u, u\in \tau)$
is a collection of independent exponential random variables with parameter $c$. Similarly we can write 
$$ \overline{\cT}_{black} = \overline{T}(\tau_{black} , m_{black})\quad {\rm
  and } \quad \overline{\cT'}= \overline{T}(\tau', m') . $$
Since the leaves of $\cT$ are exactly the leaves of $\tau$,
$\tau_{black}$ is obtained from $\tau$ by a 
$p$-Bernoulli leaf colouring. Thus, by the result of Section
\ref{secdiscrete}, conditionally on $\{\cT_{black} \neq \{ \rho \} \}$, $\tau_{black}$ is distributed as 
$\tau'$, namely as a GW($\xi_b $)-real tree. Moreover, the marks $m'$ are independent exponential random variables 
with parameter $c(1-\varphi' (g(p)))$. We need the following elementary claim whose proof is left to the reader. 
\begin{itemize}
   \item{{\it Claim}}. Let $M$ be an exponential random 
variable with parameter $\alpha >0$; consider  
an independent Poisson process on $[0, \infty)$ with intensity $\beta >0$, which splits the interval 
$[0,M]$ into $N$ subintervals with lengths $L_1, \ldots , L_N$; then, $N$ is a geometric random variable with parameter 
$\beta /(\alpha +\beta )$: 
$$ \bP (N=k+1)=\frac{\alpha}{\beta +\alpha} \, \left(\frac{\beta}{\beta +\alpha}\right)^k \; . $$
Moreover, conditionally on $N$ the $L_i$'s are independent exponentially distributed random variables with parameter
$\alpha +\beta$.  
\end{itemize}

  Now consider one edge $I \subset \cT'$ that corresponds to a vertex $u\in \tau'$ as explained in Section 
\ref{lengthtrees}. Condition on $\tau'$ and use the claim with 
$$ M=m_u \; , \; \alpha =  c(1-\varphi' (g(p))) \quad {\rm and } \quad \beta  =c\varphi' (g(p))  $$
in order to show that the Poisson point process $\cP$ splits $I$ into $N$ subintervals whose lengths are 
independent exponential variables with parameter $c$; Moreover $N$ has a geometric distribution with parameter 
$\varphi' (g(p))$. Now observe that adding the points of $\cP$ in
$\cT'$ corresponds to adding the 
line-trees to $\tau'$ 
as in Step 1 of the reconstruction procedure for discrete trees 
in Section \ref{secdiscrete}. Then, note that we next graft on $\cT'$ independent red GW($\xi_r , c$)-real 
trees according to Step 2 of the reconstruction procedure for discrete trees. Thus, we can write  
$\overline{\cT''}= \overline{T}(\tau'', m'')$ where $\tau''$ is obtained by Steps 1 and 2 of the reconstruction procedure 
for discrete trees in Section \ref{secdiscrete} and where $m''$ is a collection of independent exponential variables 
with parameter $c$. This proves the first identity of the lemma. The second one is a simple consequence of the first
one and its proof is left to the reader. \cqfd

\vspace{3mm}

We now discuss the converse problem to determine the possible offspring distributions that appear as ``black'' distributions; more precisely, we say that a proper offspring distribution $\xi_b$ is 
{\it p-extensible} if we can find a proper offspring distribution $\xi$ such that $\xi_b$ is the ``black'' distribution 
associated with a $p$-Bernoulli leaf colouring of a GW($\xi $)-tree. 

\begin{theorem}
\label{prolong}
Let $\xi_b$ be a proper offspring distribution on $\bN$. Then, the two following assertions are equivalent 
\begin{enumerate}

\item[{\rm(I)}] $\xi_b$ is $p$-extensible for all sufficiently large $p\in (0,1)$.

\item[{\rm(II)}] There exists $\psi$ that is the branching mechanism of a CSBP (thus of the form (\ref{LevKhin1})) such that  
$$ \varphi_b (r)=  r+ \psi (1-r) \; , \quad  r\in [0, 1] \, , $$ 
where $\varphi_b$ stands for the generating function of $\xi_b$. 
\end{enumerate}
\end{theorem}
{\bf Proof:} Let us first prove that (I) implies (II). With any $p\in(0,1)$ we can associate the $p$-extension $\xi$ of 
$\xi_b$ ($\xi$ depends on $p$ but we skip it for convenience of notation). Recall (\ref{genfunblack}) 
and set $v_p= g(p) /(1-g(p))$ where $g$ is defined by (\ref{equag}). Observe that (\ref{genfunblack}) implies that $\varphi_b$ 
is $C^{\infty}$ on $(-v_p , 1)$ and continuous on $[-v_p , 1]$. Moreover 
\begin{equation}
\label{complmono1}
\forall v \in (-v_p , 1)\; , \; \forall n \geq 2 \; : \quad \varphi_b^{(n)} (v) \geq 0 \,  .  
\end{equation}
We first prove the following equation
\begin{equation}
\label{vpetp}
p=1- \frac{\varphi_b(0)}{v_p + \varphi_b (-v_p)} .
\end{equation}
To that end, first note that 
\begin{equation}
\label{black0}
\varphi_b(0)=\xi_b (0)=(1+v_p)(1- \varphi_b ' (-v_p))(\varphi (g(p)) -g(p))
\end{equation}
Then, observe that 
\begin{equation}
\label{offsp0}
\varphi (0)=\xi (0)=\frac{v_p + \varphi_b (-v_p)}{(1+v_p)(1- \varphi_b' (-v_p))}. 
\end{equation}
Deduce from (\ref{equag}) that 
$$ 1-p = \frac{\varphi (g(p)) -g(p)}{\xi (0)}$$
and use (\ref{black0})and (\ref{offsp0}) to prove (\ref{vpetp}). 

\vspace{4mm}

Let us now define $v_{max}\in (0, \infty]$ by
$$ v_{max}=\sup \{ v\geq 0 \; :\; \varphi_b^{(n)} (u)   \geq 0 \; , \; u\in (-v, 1)\, , \; n \geq 2 \} . $$
Suppose that $ v_{max} <\infty $. First observe that by (\ref{vpetp}) we can find an increasing sequence 
$p_n \in (0, 1) \rightarrow 1$ such that 
$$ \lim_{n\rightarrow \infty } \varphi_b (-v_{p_n}) = +\infty .$$
Since $\varphi_b$ is convex on $(-v_{max} , 1]$, it implies that 
$$ \lim_{v\rightarrow v_{max} } \varphi_b (-v) = +\infty \quad {\rm and } \quad  
\lim_{v\rightarrow v_{max} } \varphi_b' (-v) = -\infty .$$
But the second limit is impossible for $\varphi_b'$ is a convex non-decreasing function on 
$(-v_{max} , 1)$. Thus, we must have $ v_{max} =\infty $ and (\ref{complmono1}) implies that 
$$  \forall v \in (-\infty , 1)\; , \; \forall n \geq 2 \; : \quad \varphi_b^{(n)} (v) \geq 0 \, . $$
Set $\psi (u)= \varphi_b (1-u) -1 + u $ , $u\in [0, +\infty )$. The previous observation implies that $\psi$ has the following 
properties

\begin{enumerate}
\item[(a)] $\psi (0)=0$ and $\psi' (1)=1$; 
\item[(b)] $\psi'' $ is completely monotone on $[0, +\infty )$.
\end{enumerate}

Bernstein's Theorem and a standard integration argument adapted from
the proof of Theorem 2, Chapter XIII.7 in \cite{Fe} imply that 
$\psi$ is of the form (\ref{LevKhin1}). 

\vspace{3mm}

The fact that (II) implies (I) is an easy consequence of the following computation (which is left to reader). If 
$\varphi_b (r) = r+ \psi (1-r)$ , $r\in [0, 1]$  and if $p\in (0,1)$, then the offspring distribution $\xi$, whose 
generating function $\varphi$ is given by 
\begin{equation}
 \label{varphiequa}
\varphi (r) =r + \frac{ \psi \left( (1-r) \psi^{-1} \left( \psi (1) /1-p \right)  
\right) }{\psi^{-1} \left( \psi (1) /1-p \right)  \psi' \left( \psi^{-1} \left( \psi (1) /1-p \right)  \right) }
\end{equation}
is a $p$-extension of $\xi_b$. \cqfd

\begin{remark}
\label{broadassum}
Observe that Theorem \ref{prolong} is true for offspring distributions that do not satisfy (\ref{nonexpdiscr}). If $\varphi_b $
is of the form given by Theorem \ref{prolong} (II), then it is easy to prove 
that if even $$ \varphi_b' (1)=\sum_{k\geq 0}k\xi_b (k) < \infty \; , $$
then, $\psi' (0+)$ is finite. 
\end{remark}

The main objects that we discuss in this paper are families of GW-real forests that are consistent under Bernoulli leaf colouring. 
More precisely, let $(\cF_{\lambda}; \lambda \in [0, \infty))$ be a collection of random locally compact rooted trees such that 
for any $\lambda \geq 0$ ,  $\cF_{\lambda}$ is a
GW($\xi_{\lambda},c_{\lambda} , a_{\lambda})$-real 
forest, such that for any $\lambda >0$, $\xi_{\lambda}$ 
is a proper offspring distribution satisfying $\xi_{\lambda} (0) >0 $
and such that $c_{\lambda}$ and $a_{\lambda}$ are non-negative real numbers. We say that 
 $(\cF_{\lambda}; \lambda \in [0, \infty))$ is {\it Bernoulli leaf colouring consistent} if for any 
$0\leq \mu\leq \lambda$, $\cF_{\mu}\subset\cF_\lambda$ and $\cF_{\mu}$ 
is obtained from $\cF_{\lambda}$ as the ``black'' tree resulting from a 
$p$-Bernoulli leaf colouring with $1-p=\mu/\lambda$. According to Lemma \ref{percola} it implies that 
$\xi_{\lambda}$ is the $(1-\mu/\lambda)$-extension of
$\xi_{\mu}$. Therefore, $\xi_{\mu}$ is $p$-extensible 
for any sufficiently 
large $p$ and $\xi_{\mu}$ has to be of the form given by Theorem \ref{prolong} (II). 
Accordingly, up to a linear time change of the family $(\cF_{\lambda};
\lambda \in [0, \infty))$, there 
is a unique function 
$\psi$ satisfying $(\ref{LevKhin1})$ such that for any $\lambda \geq 0$: 
\begin{equation}
\label{xilambda}
\xi_{\lambda}(k)= \frac{\psi^{-1}(\lambda)^{k-1}\left| \psi^{(k)}(\psi^{-1}(\lambda))\right|}{k! 
\psi'(\psi^{-1}(\lambda))} \quad {\rm if } \quad k\neq 1 
\end{equation}
and $\xi_{\lambda}(1)=0$. The generating function $\varphi_{\lambda }$ of $\xi_{\lambda }$ is then given by 
\begin{equation}
\label{varphilambda}
 \varphi_{\lambda }(s)= s+ \frac{\psi\left( (1-s)\psi^{-1}(\lambda)\right)}{\psi^{-1}(\lambda)\psi'(\psi^{-1}(\lambda))}. 
\end{equation}

\begin{remark}
\label{gammapsi}
Recall that $\gamma $ is the largest root of $\psi$. Thus $\gamma >0$ iff 
$m=\psi'(0+)<0$. Observe that if $\gamma >0$, then  
$$\varphi_{0}(s)= s+ \frac{\psi\left( (1-s)\gamma \right)}{\gamma \psi'(\gamma)}\;  $$ 
and $\xi_0 (0)=0$. A GW($\xi_0$)-tree is infinite with no leaf. 
\end{remark}

By definition, the black distribution associated with $\xi_{\lambda}$ via a $(1-\mu/\lambda)$-Bernoulli leaf colouring is 
$\xi_{\mu}$. It is also easy to compute the function $g$ that solve (\ref{equag}). Namely, 
\begin{equation}
\label{gequa}
 g(s):= 1-\frac{\psi^{-1}\left((1-s)\lambda
 \right)}{\psi^{-1}(\lambda)} \; .
\end{equation}
Thus, the probability for a GW($\xi_{\lambda })$-tree to be completely red is 
\begin{equation}
\label{redproba}
g(p)= 1-\frac{\psi^{-1}(\mu)}{\psi^{-1}(\lambda)}. 
\end{equation}
The red distribution associated with $\xi_{\lambda }$ via a $(1-\mu/\lambda)$-Bernoulli leaf colouring is denoted by 
$\xi_{\mu, \lambda}:=\xi_r$ and is given by 
\begin{equation}
\label{ximulambda}
\xi_{\mu, \lambda}(k)= \frac{\left| \psi^{(k)}(\psi^{-1}(\lambda))\right|}{
   \psi'(\psi^{-1}(\lambda))}
\frac{\left( \psi^{-1}(\lambda)-\psi^{-1}(\mu) \right)^{k-1}}{k!}
\; ,\quad k\geq 2 ,
\end{equation}
$$ \xi_{\mu, \lambda}(1)=0 \quad {\rm and } \quad \xi_{\mu, \lambda}(0)=\frac{\lambda -\mu}{
 \left( \psi^{-1}(\lambda)-\psi^{-1}(\mu) \right) \psi'(\psi^{-1}(\lambda))}. $$
The generating function of $\xi_{\mu , \lambda}$ is denoted by $\varphi_{\mu , \lambda}$ 
and is given by  
\begin{equation}
\label{redgenfun}
 \varphi_{\mu , \lambda}(s)= s + \frac{ \psi \left(  \psi^{-1}(\lambda)
- s(\psi^{-1}(\lambda)-\psi^{-1}(\mu))\right)-\mu }{
(\psi^{-1}(\lambda)-\psi^{-1}(\mu))\psi'(\psi^{-1}(\lambda))}. 
\end{equation}

\begin{remark}
\label{mulambda}
Observe that $\varphi'_{\mu , \lambda}(1)=1- \psi'(\psi^{-1}(\mu)) /
\psi'(\psi^{-1}(\lambda)) <1 $. Thus, for $\mu<\lambda$, 
$\xi_{\mu , \lambda}$ is a subcritical offspring distribution and therefore, any 
GW($\xi_{\mu , \lambda}$)-real tree is a.s. finite. 
\end{remark}

  For any $l\geq 1$ we denote by $\nu^{\mu,\lambda}_l$ the distribution given by (\ref{addvert}) with 
$\varphi=\varphi_{\lambda}$, $\xi =\xi_{\lambda}$ and $g(p)$ as in (\ref{redproba}). For any $l\geq 2$, 
$\nu^{\mu,\lambda}_l$ is given by 
\begin{equation}
\label{nulevy}
\nu^{\mu, \lambda}_{l}(k)= \frac{\left| \psi^{(l+k)}(\psi^{-1}(\lambda))\right|}{\left|
    \psi^{(l)}(\psi^{-1}(\mu)) \right|}
\frac{\left( \psi^{-1}(\lambda)-\psi^{-1}(\mu) \right)^k}{k!}
\; ,\quad k\geq 0 
\end{equation}
and for $l=1$ 
\begin{equation}
\label{nulevy1}
\nu^{\mu, \lambda}_{1}(k)= \frac{\left| \psi^{(1+k)}(\psi^{-1}(\lambda))\right|}{
   \psi'(\psi^{-1}(\lambda))- \psi'(\psi^{-1}(\mu)) }
\frac{\left( \psi^{-1}(\lambda)-\psi^{-1}(\mu) \right)^k}{k!}
\; ,\quad k\geq 1 ,
\end{equation}
with $\nu^{\mu, \lambda}_{1}(0)=0$. Now observe that the parameter 
$\varphi'(g(p))$ of the geometric distribution in (\ref{attach}) is given by 
\begin{equation}
\label{rateattach}
\varphi'(g(p))= 1-\frac{\psi'(\psi^{-1}(\mu))}{\psi'(\psi^{-1}(\lambda))} \; .
\end{equation}
Then, according to Lemma \ref{percola} we have
$$ a_{\mu}/ a_{\lambda}= \psi^{-1}(\mu)/\psi^{-1}(\lambda) \quad {\rm and }\quad  
c_{\mu}/ c_{\lambda}= \psi' (\psi^{-1}(\mu))/\psi' (\psi^{-1}(\lambda)) . $$
{\bf We choose the following  normalization:} 
\begin{equation}
\label{normali}
 a_{\lambda}= a \, \psi^{-1}(\lambda)
 \quad {\rm and }\quad  c_{\lambda}= \psi' (\psi^{-1}(\lambda)), 
\end{equation}
where $a>0$. Such a Bernoulli leaf colouring consistent family $(\cF_{\lambda}; \lambda \in [0, \infty))$ whose distribution 
is specified by (\ref{xilambda}) and 
(\ref{normali}) is called an 
$(a, \psi)$-{\it L\'evy growth process}. 
\begin{remark}
\label{psishift} For any $\mu \geq 0$, we set 
$$ \psi_{\mu} (x)= \psi ( x+\psi^{-1} (\mu ))-\mu . $$
Then, $\psi'_\mu(0+) = \psi'(\psi^{-1} (\mu ))$. If $\mu  >0$, then $\psi'_{\mu}(0+) $ is finite. 
It is also easy to check that 
$$ \psi^{-1}_{\mu} ( x)=\psi^{-1}(x+ \mu)- \psi^{-1}(\mu) ;$$
Thus for any $\mu \leq \lambda$
$$\left( \psi_{\lambda- \mu}\right)_{\mu}= \psi_{\lambda} . $$
Note that $\xi_{\lambda}$ and $\xi_{\mu, \lambda}$ actually depend on $\psi$: 
$\xi_{\lambda ; \psi}=\xi_{\lambda}$ , 
$ \xi_{\mu , \lambda ; \psi}=\xi_{\mu , \lambda}$. Then, it is easy to check that
for any $\mu_0 \leq \mu \leq \lambda$:
\begin{equation}
\label{xishift}
 \xi_{\mu, \lambda}=\xi_{\lambda -\mu ; \psi_{\mu}} \quad {\rm and} \quad 
\xi_{\mu , \lambda ; \psi_{\mu_0}}= \xi_{\mu +\mu_0 , \lambda +\mu_0 } .
\end{equation}
\end{remark}

\begin{notation}
\label{notGWreal} Fix $0\leq \mu \leq \lambda$ and $a>0$. We shall use the following notation. We denote by
\begin{itemize}\item $\Delta_{\lambda} (d\overline{T})$ the distribution on $\bT$ of the isometry class of 
a GW($\xi_{\lambda}, \psi' (\psi^{-1}(\lambda)) $)-real tree,
\item $\Delta_{\lambda}^a (d\overline{T})$ the distribution on $\bT$ of the isometry class of 
a GW($\xi_{\lambda}, \psi' (\psi^{-1}(\lambda)),  a \psi^{-1}(\lambda) $)-real forest,
\item $\Delta_{\mu ,\lambda} (d\overline{T})$ the distribution on $\bT$ of the isometry class of 
a GW($\xi_{\mu ,\lambda}, \psi' (\psi^{-1}(\lambda)) $)-real tree,
\item $\Delta_{\mu, \lambda}^a (d\overline{T})$ the distribution on $\bT$ of the isometry class of 
a GW($\xi_{\mu,\lambda}, \psi' (\psi^{-1}(\lambda)), \\ 
 a (\psi^{-1}(\lambda)-\psi^{-1}(\mu) $)-real forest. 
\end{itemize}
 According to the previous remark, we get 
\begin{equation}
\label{deltashift}
\Delta_{\mu,\lambda}^a= \Delta_{\lambda -\mu ; \psi_{\mu}}^a \quad {\rm and } \quad 
\Delta_{\mu,\lambda}= \Delta_{\lambda -\mu ; \psi_{\mu}}
\end{equation}
with an obvious notation. Observe also that $ \Delta_{\lambda}^0= \delta_{\overline{\{\rho\}}}$ that is the Dirac mass at the
isometry class $\overline{\{\rho\}}$ of the point tree. Thus  $\Delta_{\lambda}^0 \neq  \Delta_{\lambda}$. 
\end{notation}

\subsection{Construction of the growth process.}
\label{percolevy}

  In this section we discuss how to grow a tree in order to obtain Bernoulli colouring consistent families of GW-real trees and 
related tree-valued processes. The definition given in this subsection 
is slightly more general for we want to start the growth process at any 
discrete real tree with edge lengths. Let $(T, d, \rho)$ be such a tree. Fix $0\leq \mu\leq \lambda$ and $a>0$. 
Let $\psi$ be of the form (\ref{LevKhin1}) such that $\psi'(0+)$ is finite.
We first define a random tree denoted by $Q_{\mu, \lambda}^a (T)$ via the following grafting procedure: 

\begin{itemize}
\item 
{\bf The grafting procedure on} $T$:
Let $\cP$ be a Poisson point process on $T$ with intensity 
$$  \left(\psi'(\psi^{-1}(\lambda))-\psi'(\psi^{-1}(\mu))\right) \, 
\ell_T \; .$$
Graft a random number $N_{\sigma}$ of independent 
GW($\xi_{\mu, \lambda}, \psi' (\psi^{-1} (\lambda))$)-real trees on 
each vertex $\sigma \in \cP \cup{\rm Br}(T)$; 
here $N_{\sigma}$ has distribution 
$\nu^{\mu, \lambda}_l$, where $l={\rm n}(\sigma , T)-1$. The resulting tree is denoted by $Q_{\mu , \lambda }(T)$. 
Then, 
graft on $\rho$ a random number 
$N_{\rho}$ of independent GW-real trees with the same distribution, where $N_{\rho}$ is a Poisson random variable with parameter
$a(\psi^{-1}(\lambda)-\psi^{-1}(\mu))$. Denote by $Q_{\mu, \lambda}^a (T)$ the resulting tree.
\end{itemize}

\begin{remark}
\label{Qzero}
Observe that $ Q^0_{\mu , \lambda } (T)= Q_{\mu , \lambda  } (T)$ and note that if $T$ reduces to its root $\rho$ 
then $Q_{\mu , \lambda  } (T)=\{\rho\} $. 
\end{remark}

Consider a GW($\xi_{\lambda}, \psi' (\psi^{-1} (\lambda))$)-real tree (resp. a 
GW($\xi_{\lambda}, \psi' (\psi^{-1} (\lambda)), a\psi^{-1}(\lambda)$)-real forest) denoted by $\cT(\lambda)$ (resp. by
 $\cF(\lambda)$). Denote by $\cT_{\mu}(\lambda)$ (resp. by
 $\cF_{\mu}(\lambda)$) the black subtree obtained by a $(1-\mu/\lambda)$-Bernoulli leaf colouring of $\cT(\lambda)$ (resp. of 
$\cF(\lambda)$). Let $\cT'$ (resp. $\cF'$) be a GW($\xi_{\mu}, \psi' (\psi^{-1} (\mu))$)-real tree (resp. a 
GW($\xi_{\mu}, \psi' (\psi^{-1} (\mu)), a\psi^{-1}(\mu)$)-real
forest). The grafting procedure corresponds to the reconstruction
procedure explained at the beginning of the previous section and Lemma \ref{percola} implies that

\begin{equation}
\label{percoarbre}
 \left( \overline{\cT}_{\mu}(\lambda ) \; ,\;  \overline{\cT}(\lambda )
\right) \; {\rm under} \; \bP \left( \; \cdot \; | 
\cT_{\mu}(\lambda ) \neq \{ \rho \} \right) \quad 
\overset{(d)}{=} \quad \left( \overline{\cT'},  \overline{Q}_{\mu, \lambda}(\cT') \right)
\end{equation}
and 
\begin{equation}
\label{percoforet}
\left( \overline{\cF}_{\mu}(\lambda ) \; ,\;  \overline{\cF}(\lambda )
\right)  \quad 
\overset{(d)}{=} \quad \left( \overline{\cF'},  \overline{Q}^a_{\mu, \lambda}(\cF')\right). 
\end{equation}
\noindent 
Here $\overline{Q}_{\mu, \lambda}(\cT')$ and $\overline{Q}^a_{\mu, \lambda}(\cF')$ stand for the isometry 
classes of resp. $Q_{\mu, \lambda}(\cT')$ and $Q^a_{\mu, \lambda}(\cF')$ (the extra random variables 
used to define the grafting procedures on $\cT'$ and $\cF'$ are chosen independent of these trees).

The grafting procedure enjoys a Markov-like property in the following sense: 
fix $a\geq 0$ and let $0\leq \lambda_1 <\lambda_2<\lambda_3$. 
Set $\cF=Q^a_{\lambda_1, \lambda_3 }(T)$. Let 
$U_{\sigma}$ , $\sigma\in {\rm Lf}(\cF)\backslash {\rm Lf}(T) $ be $[0, 1]$-uniform independent random variables conditionally 
on $\cF$. We define 
$$ \cF_b= T \bigcup \left\{  \lgeo \rho , \sigma \rgeo  \; ; \; \sigma \in  {\rm Lf}(\cF)\backslash {\rm Lf}(T) 
\; : \; U_{\sigma} \leq \frac{\lambda_2 -\lambda_1}{\lambda_3- \lambda_1}  \right\} \; .$$
$\cF_b$ is thus the black tree resulting from a $(1-(\lambda_2
-\lambda_1)/(\lambda_3 -\lambda_1))$-Bernoulli 
colouring 
of the leaves of $\cF$ that are not in $T$. 

\begin{proposition}
\label{chaining}
For any discrete tree with edge lengths $T$, any $a\geq 0$ and any $0\leq \lambda_1 <\lambda_2<\lambda_3$, we have
$$\left( \overline{Q}^a_{\lambda_1 , \lambda_2}(T) \; ,\;  \overline{Q}^a_{\lambda_2 , \lambda_3} 
\left( Q^a_{\lambda_1 , \lambda_2}(T)\right)
\right)  \quad 
\overset{(d)}{=} \quad \left(   \overline{\cF}_b  \; , \;\overline{\cF} \right)$$
(here the extra random variables used to define $Q^a_{\lambda_2 , \lambda_3} $ are chosen independent of 
$Q^a_{\lambda_1 , \lambda_2} (T)$).
\end{proposition}
{\bf Proof:} By Remark \ref{psishift}, 
we only have to prove 
\begin{equation}
\label{reduced}
\left( \overline{Q}^a_{0 , \mu}(T) \; ,\;  \overline{Q}^a_{\mu , \lambda} 
\left( Q^a_{0 , \mu}(T)\right)
\right)  \quad 
\overset{(d)}{=} \quad \left(   \overline{\cF}_b  \; , \;\overline{\cF} \right) 
\end{equation}
by replacing $\psi$ by $\psi_{\lambda_1}$ and by taking $\mu= \lambda_2-\lambda_1$ and $\lambda =\lambda_3 - \lambda_1$ 
in (\ref{reduced}).

   Let us denote by $\cP$ the Poisson point process on $T$ involved in the grafting procedure defining $\cF$. For any
$\sigma \in \cP\cup {\rm Br} (T)\cup \{ \rho \}$, we denote by $\cT^i_{\sigma}$ , $1\leq i\leq N_{\sigma}$, the trees grafted 
on $\sigma$. Denote by $\cT^i_b (\sigma)$ the tree $\cF_b \cap \cT^i_{\sigma}$ and
set 
$$ J_r (\sigma ) = \left\{ i\in \{ 1, \ldots, N_{\sigma} \} \; : \; \cT^i_b (\sigma)= \{ \sigma \}  \right\} 
\quad {\rm and} \quad J_b (\sigma )=\{ 1, \ldots, N_{\sigma} \} \backslash J_r (\sigma ). $$
Then, observe that performing a $(1-\mu/\lambda)$-Bernoulli leaf colouring on 
${\rm Lf}(\cF)\backslash {\rm Lf}(T)$ is the same as performing independent $(1-\mu/\lambda)$-Bernoulli leaf 
colourings on the $\cT^i_{\sigma}$'s. Accordingly, conditionally on $J_r (\sigma )$ and on 
$J_b (\sigma )$ the pairs of trees $(\cT^i_{\sigma}, \cT^i_b (\sigma) )$, $i\in J_b (\sigma )$, and the trees 
$\cT^i_{\sigma}$ , $i\in J_r (\sigma )$, are independent; moreover, by (\ref{percoarbre}), conditionally on  
$J_r (\sigma )$ and on $J_b (\sigma )$ the isometry classes of  $(\cT^i_{\sigma}, \cT^i_b (\sigma) )$,  
$i\in J_b (\sigma )$, are independent copies of $(\overline{Q}_{\mu , \lambda} (\cT), \overline{\cT})$ where $\cT$ is a 
GW($\xi_{\mu}, \psi' (\psi^{-1} (\mu))$)-real tree. To simplify notation we assume that for 
any $\sigma \in \cP\cup {\rm Br} (T)\cup \{ \rho \}$ and any $i\in J_b (\sigma )$ 
\begin{equation}
\label{psassum}
Q_{\mu , \lambda} (\cT^i_b (\sigma))= \cT^i_{\sigma}. 
\end{equation}

Now recall that conditionally on $N_{\sigma}$, the events that $
\cT^i_{\sigma}$ is completely red, $i\in \{1, \ldots , N_{\sigma} \}$, are independent events 
and have probability 
$1-\psi^{-1} (\mu)/\psi^{-1} (\lambda)$. Since $N_{\sigma}$ has distribution 
$\nu^{0, \lambda}_l$, with $l={\rm n} (\sigma , T)-1$, we get 
\begin{multline}
\label{jointJ}
\bP \left( \# J_b (\sigma)=k_b;\# J_r (\sigma)=k_r \arrowvert \cP 
\right)= \\ 
\frac{(k_b+k_r)!}{k_b! k_r!} \left(1-\frac{\psi^{-1} (\mu)}{\psi^{-1} (\lambda)} \right)^{k_r} 
\left(\frac{\psi^{-1} (\mu)}{\psi^{-1} (\lambda)} \right)^{k_b} \nu^{0,\lambda}_l (k_b+k_r).
\end{multline}
Now set $\cP_1=\{ \sigma \in \cP \; :\;  \# J_b (\sigma) \geq 1\} $ and $\cP_2= \cP \backslash \cP_1$. It is easy to deduce from the 
latter observations that $\cP_1$ and $\cP_2$ are independent Poisson point processes with respective intensities 
$$ \psi^{-1}(\mu) \, \ell_T \quad {\rm and} \quad (\psi^{-1}(\lambda)-\psi^{-1}(\mu)) \, \ell_T \; . $$
A long but straightforward computation based on (\ref{jointJ}) (which is left to the reader) implies that conditionally on 
$\cP_1$ and $\cP_2$ the following assertions are true: 
\begin{enumerate}
\item[(a)]  If $\sigma \in \cP_1$, then $\# J_b (\sigma)$ has distribution $\nu^{0,\mu}_1$ and 
conditionally on $\# J_b (\sigma)=j$, $\# J_r (\sigma)$ has distribution $\nu^{\mu, \lambda}_{j+1}$; 

\item[(b)]  If $\sigma \in \cP_2$, then $\# J_r (\sigma)$ has distribution $\nu^{\mu, \lambda}_1$; 

\item[(c)] If $\sigma \in {\rm Br}(T)$, then $\# J_b (\sigma)$ has distribution $\nu^{0,\mu}_l$, 
where $l={\rm n}(\sigma , T)-1$. Moreover, conditionally on $\# J_b (\sigma)=j$, $\# J_r (\sigma)$ has distribution $\nu^{\mu, \lambda}_{j+l}$; 

\item[(d)]  $\# J_b (\rho)$ and $\# J_r (\rho)$ are independent Poisson random variables with respective parameters 
 $\psi^{-1}(\mu)a$ and $(\psi^{-1}(\lambda)-\psi^{-1}(\mu))a$. 

\end{enumerate} 

Next, observe that 
\begin{equation}
\label{effbe}
\cF_b = T \circledast_{\substack{ \sigma \in \cP_1 \cup {\rm Br}(T)\cup \{ \rho \} \\ i\in J_b (\sigma)}} 
\; \left( \sigma , \cT^i_b (\sigma) \right). 
\end{equation}
According to the distribution of $\cP_1$ and of $J_b (\sigma)$, $\sigma \in \cP_1 \cup {\rm Br}(T)\cup \{ \rho \}$, 
(\ref{effbe}) implies that $\cF_b$ is obtained from $T$ by the grafting procedure corresponding to the 
``grafting operator'' $Q^a_{0,\mu}$ and, more precisely, that $\overline{\cF_b}$ has the same distribution as the isometry class 
of $Q^a_{0,\mu}(T)$. To simplify notation we assume that 
\begin{equation}
\label{assueffbe}
\cF_b=Q^a_{0,\mu}(T). 
\end{equation}
We now graft trees on $\cF_b$ according to the ``grafting operator'' $Q^a_{\mu, \lambda}$. Observe that this procedure can be 
split in the three following steps: 
\begin{enumerate}
\item[(i)] Graft trees according the ``grafting operator'' $Q_{\mu, \lambda}$ independently on each 
$\cT^i_b (\sigma)$, $ i\in J_b (\sigma)$ , $\sigma \in \cP_1 \cup {\rm Br}(T)\cup \{ \rho \}$. Note that by  (\ref{psassum})
the resulting trees have the same distribution as the $\cT^i_{\sigma}$'s . 

\item[(ii)] Choose additional grafting points on $T$ according to a Poisson point process with the same distribution as 
$\cP_2$. We denote this set of points by $\cP_2'$.

\item[(iii)] Graft a random number of independent 
GW($\xi_{\mu , \lambda} , \psi' (\psi^{-1}(\lambda ))$)-real trees at each $ \sigma \in \cP_1 \cup \cP_2' \cup {\rm Br}(T)\cap \{ \rho \}$, the random number of trees grafted on $\sigma$ having distribution $\nu^{\mu, \lambda}_l$, with $$l={\rm n}(\sigma , \cF_b)-1= {\rm n}(\sigma , T)-1 + \# J_b (\sigma) .$$  
If $\# J_b (\sigma)=j$, then by the grafting procedure, $\# J_r (\sigma)$ is distributed 
$\nu^{\mu,\lambda }_{l+j}$, 
and the resulting trees have the same distributions as $\cT^i_{\sigma}$ , $i\in J_b (\sigma)$, 
$ \sigma \in \cP_1 \cup \cP_2 \cup {\rm Br}(T)\cup \{ \rho \}$.
\end{enumerate}

This implies that the isometry class of $Q^a_{\mu, \lambda}(\cF_b)$
has the same distribution as 
$\overline{\cF}$ and it completes the proof of (\ref{reduced}) by (\ref{effbe}). \cqfd

\vspace{3mm}

Fix $\lambda>0$, set $\cF (\lambda)=Q^a_{0, \lambda} (T)$, and for all $\mu \in [0, \lambda]$ set 
$$ \cF_{\mu}(\lambda) =T \, \bigcup \left\{  \lgeo \rho , \sigma \rgeo \; ; \; \sigma \in 
{\rm Lf}(\cF (\lambda))\backslash {\rm Lf}(T) \; : \; U_{\sigma } \leq \mu /\lambda \right\} , $$
where the $U_{\sigma}$'s are i.i.d. $[0, 1]$-uniform variables conditionally on $\cF (\lambda)$. The following 
proposition discusses how to construct an $(a, \psi)$-{\it growth process} starting from a discrete tree with 
edge lengths $(T, d, \rho)$. 
\begin{proposition}
\label{growthproc} Assume that $m=\psi' (0+)$ is finite. Then, there exists a family of random rooted locally compact real trees 
$(\cF_{\lambda }, d_{\lambda}, \rho)$, $\lambda \in [0, \infty)$ such that a.s. 
\begin{enumerate}
\item[{\rm(i)}] For any $0 \leq \mu \leq \lambda $
$$ \cF_{\mu } \subset \cF_{\lambda } \quad {\rm and} \quad  
d_{\mu} = d_{\lambda \, |\cF_{\mu }\times \cF_{\mu }}\; . $$

\item[{\rm(ii)}] The map $\lambda \longrightarrow \overline{\cF}_{\lambda }$ is cadlag in $(\bT, \delta)$ and 
$$ \left( \overline{\cF}_{\mu  }  \; , \; 0\leq \mu \leq \lambda \right)
\overset{(d)}{=}\left( \overline{\cF}_{\mu  } (\lambda) \; , \; 0\leq \mu \leq \lambda \right) \; . $$
\end{enumerate}
\end{proposition}

\noindent
{\bf Proof:} Let $(\lambda_n; n\geq 0)$ be an increasing sequence that goes to $\infty$ and such that $\lambda_0=0$.   
Set $\cF_0=T$ and define the sequence $(\cF_{\lambda_n}; n\geq 1)$ by  
$\cF_{\lambda_{n+1}}= Q^a_{\lambda_n , \lambda_{n+1}}(\cF_{\lambda_n} ), n\geq 0 $, where 
the extra random variables used in the grafting procedure at step $n$
are chosen to be independent of $\cF_{\lambda_n}$. 
Associate a random variable 
$V_{\sigma}$ with any 
$\sigma \in  \cup {\rm Lf}(\cF_{\lambda_n}) \backslash {\rm Lf}(T)$ such that conditionally on the sequence $(\cF_{\lambda_n}; n\geq 1)$, 
the $V_{\sigma}$'s are i.i.d. uniformly 
distributed in $[0, 1]$. Then, for any $\lambda\in [\lambda_n  ,\lambda_{n+1} )$ we define the growth process
as follows:
$$ \cF_{\lambda }=\cF_{ \lambda_n } 
\cup \, \bigcup \left\{ \lgeo \rho , \sigma \rgeo \; , \; 
\sigma \in  {\rm Lf} (\cF_{\lambda_{n+1}}) 
\backslash  {\rm Lf} (\cF_{\lambda_{n}})  \; : \; V_{\sigma}\leq 
\frac{\lambda - \lambda_{n}}{\lambda_{n+1} - \lambda_{n}} \right\} $$ 
and  
$$ d_{\lambda} = d_{\lambda_{n+1} \;  |\cF_{ \lambda}\times \cF_{\lambda}} \; , $$
Thus, point (i) clearly holds and it implies that 
$\lambda \longrightarrow \overline{\cF}_{\lambda }$ is cadlag in $(\bT, \delta)$. 
Fix $n\geq 0$ and take $\lambda=\lambda_{n+1}$. Then use Proposition \ref{chaining}
successively with $\mu =\lambda_0, \ldots , \lambda_n$ to prove that 
the joint distribution of $ \overline{\cF}(\lambda_{n+1})$ and 
$$  \left( \overline{\cF}_{\lambda_{k}  }(\lambda_{n+1}) \, ; \; U_{\sigma }  \, , \, \sigma \in 
{\rm Lf} (\cF_{\lambda_{k+1}}(\lambda_{n+1})) 
\backslash  {\rm Lf} (\cF_{\lambda_{k}}(\lambda_{n+1}))
  \, , \;  0\leq k\leq n \right) $$
is the same as the joint distribution of $\overline{\cF}_{\lambda_{n+1}}$ and 
$$\left( \left( \overline{\cF}_{\lambda_{k} } \, ; \;\frac{\lambda_{k}+ (\lambda_{k+1} - \lambda_{k})V_{\sigma }}{\lambda_{n+1}} \, , \, \sigma \in 
{\rm Lf} (\cF_{\lambda_{k+1}}) 
\backslash  {\rm Lf} (\cF_{\lambda_{k}}) \right)  \, , \;  0\leq k\leq n \right) .$$
Thus, for any $n \geq 0$: 
$$ \left( \overline{\cF}_{\mu  }  \; , \; 0\leq \mu \leq \lambda_{n+1} \right)
\overset{(d)}{=}\left( \overline{\cF}_{\mu  } (\lambda) \; , \; 0\leq \mu \leq \lambda_{n+1} \right) \; , $$
which implies the second part of the proposition by an easy argument. \cqfd

\begin{remark}
\label{ellungrowth}
Following the construction given in the proof of Proposition \ref{embedding}, we can embed the growth process in 
$\ell^1 (\bN)$ and we obtain a non-decreasing cadlag process in $(\btl,
\dl)$. 
\end{remark}

\begin{remark}
\label{jumpchain}

Observe that the distribution of $\overline{Q}^a_{\mu, \lambda}(T)$
only 
depends on the isometry class of 
$(T, d, \rho)$ so it makes sense to denote by ${\bf P}_{\mu, \lambda} 
(\overline{T}, d\overline{\cT})$ the 
distribution on $\bT$ of $\overline{Q}^a_{\mu, \lambda}(T)$. Proposition \ref{chaining} and (\ref{percoforet})
imply that the isometry classes $(\overline{\cF}_{\lambda}; \lambda \geq 0)$ of  
a $(a, \psi)$-L\'evy growth process as defined in the end of Section
\ref{levygal} 
is a $\bT$-valued inhomogeneous Markov process 
with transition kernel ${\bf P}_{\mu, \lambda} (\overline{T},
d\overline{\cT}) $ (in the Brownian case 
$\psi (\lambda)=\lambda^2 /2$, Pitman and Winkel in \cite{PitWink} proved 
that this process has independent growth increments expressed by a composition rule).
Observe, however, that 
$\overline{Q}^a_{\mu, \lambda}(T)$ is only defined for discrete trees with edge lengths.

   More specifically, it is clear from the construction that the growth 
process $(\cF_{\lambda})_{\lambda\geq 0}$ is a pure jump process
obtained by adding 
single branches. More precisely, we get  
the following jump-chain with 
holding times construction of the process of $(\cF_\lambda)_{\lambda\ge 0}$
started at a \em compact \em discrete tree with edge lengths $(T, d, \rho)$. 
The equivalence classes of 
$(\cF_{\lambda})_{\lambda\geq 0}$ have the same distribution as the 
equivalence classes of the non-decreasing family of real trees 
$({\widetilde \cF}_{\lambda})_{\lambda\geq 0}$ that 
has a discrete set of jump times
$(\Lambda_n)_{n\ge 1}$ at which branches of lengths $(L_n)_{n\ge 1}$ are 
added, at locations $(\Sigma_n)_{n\ge 1}$
and such that the process 
  $(\Lambda_n,\Sigma_n,L_n,{\widetilde \cF}_{\Lambda_n})_{n\ge 0}$ 
is a Markov chain with
transition kernel
\begin{eqnarray*}
\bP \left( \left.
\Lambda_{n+1} \in d\lambda \, ;\, 
\Sigma_{n+1}  \in d\sigma \, ; \, L_{n+1}=dy \, ; \, {\widetilde \cF}_{\Lambda_{n+1}} \in 
 dT'  \,  \right| \Lambda_{n}=\mu  \, ;  \, 
{\widetilde \cF}_{\Lambda_{n}}=T    \right)  \\
 =\psi' (\psi^{-1}(\lambda))\,  \exp
 \left( -\psi' (\psi^{-1}(\lambda))\, y - 
\int_\mu^\lambda ds < \! M_{s,T} \!>  \right) \,  \\
\times  \quad  d\lambda \; M_{\lambda , T}  (d\sigma )  \; dy  \;   
\delta_{\{T*(\sigma,\lgeo 0,y\rgeo)\}}  (dT')  
\end{eqnarray*}
where 
\begin{eqnarray*}
\psi' (\psi^{-1}(\lambda)) \; M_{\mu,T}(d\sigma ) =  
\psi''(\psi^{-1}(\mu))\ell_T   (d\sigma)  \; &+& \\
\sum_{v \in Br(T)\backslash\{\rho\}} &&
\frac{|\psi^{({\rm n}(\sigma,T))}(\psi^{-1}(\mu))|}{|\psi^{({\rm n}(\sigma,T)-1)}
(\psi^{-1}(\mu))|}\delta_v (d\sigma) +a\delta_{\rho} (d\sigma)
\end{eqnarray*}
and $< \! M_{s,T} \!>$ stands for the total mass of $M_{s,T} $. 
Since the result is not important in the sequel, we skip the proof that 
is an easy consequence of the grafting procedure. 
\end{remark}

For technical purposes, we end the subsection by  
providing an alternative definition of the grafting procedure 
that is less direct but that is used in the proofs of the results of 
the next section. For convenience of notation we set 
\begin{equation}
\label{qandc}
q=\psi'(\psi^{-1}(\lambda)) \quad {\rm and} \quad
c=\psi^{-1}(\lambda)-\psi^{-1}(\mu).
\end{equation}
Then, for any non-negative integer $l \geq 2$, we define 
the distribution $\eta_{\mu ,l} (dx)$ on $[0, \infty)$ by 
\begin{equation*}
\eta_{\mu ,l} (dx)=\frac{2\beta\, {\bf 1}_{\{ l=2\}} }{|\psi^{(2)} (\psi^{-1}(\mu))|} 
\delta_0 (dx) 
+ \frac{ x^l e^{-x\psi^{-1}(\mu)}}{|\psi^{(l)} (\psi^{-1}(\mu))|}
\Pi(dx) .
\end{equation*}
It is easy to check that $\eta_{\mu ,l} (dx)$ is a probability measure. Let  
$$ \cP_1 = \{ (\sigma_i^{(1)} , x_i) \; , \; i\in I^{(1)}\}  \quad {\rm and } \quad  
\cP_2 = \{ (\sigma_i^{(2)} , y_i) \; , \; i\in I^{(2)}\} $$
be two independent Poisson point processes on $\bT\times [0, \infty)$  with respective intensities
$$ \ell_T (d\sigma)\otimes x e^{-\psi^{-1}(\mu)x} \Pi (dx) \quad {\rm and } \quad 2\beta \; 
\ell_T (d\sigma) \otimes dy . $$
We shall use the following notation: define for any $k \in \{1, 2\}$, 
$$S^{(k)}_{\mu}=\left\{ \sigma_i^{(k)} \; , \; i\in I^{(k)} \right\} $$
and set $ S_{\mu}= S^{(1)}_{\mu} \cup S^{(2)}_{\mu} \cup {\rm Br}(T) \cup \left\{ \rho \right\}  $ and 
$$  S'_{\mu}= S^{(1)}_{\mu}\cup   {\rm Br}(T) \cup \left\{ \rho \right\} . $$
We then introduce the collection of random variables $\cA_{\mu}= \{ a_{\sigma
  }(\mu), \; \sigma \in S'_{\mu} \}$ that are distributed as follows: 
\begin{itemize}\item $a_{\sigma_i^{(1)}} (\mu)=x_i$ , $i\in  I^{(1)}$;  
\item $(a_{\sigma  }(\mu), \; \sigma \in  {\rm Br}(T) \cup \{ \rho\} )$ is a set of independent real-valued random variables  
independent of 
$\cP_1$ and $\cP_2$. Moreover, $a_{\rho }(\mu)=a$ and 
for any $\sigma \in {\rm Br}(T)\setminus{\{ \rho\} }$ , $a_{\sigma }(\mu)$ is distributed according to 
$\eta_{\mu ,l} (dx)$ where $l={\rm n}(\sigma , T)-1$. 
\end{itemize}

\vspace{3mm}

   We next define a collection of random trees 
$ \{ (F_{\sigma }(\lambda), d_{\sigma , \lambda }, \rho_{\sigma,\lambda} ) \; , \; \sigma \in S_{\mu} \} $ independent conditionally on 
$\cP_1$ , $\cP_2$ and $\cA_{\mu}$, and whose conditional distribution 
is given as follows.

\begin{itemize}
\item If $\sigma \in S'_{\mu}$, then $(F_{\sigma }(\lambda), d_{\sigma  , \lambda }, \rho_{\sigma,\lambda})$
is distributed as a GW($ \xi_{\mu, \lambda} , q, c \, a_{\sigma })$-real forest {\it with the convention that}
$F_{\sigma }(\lambda)=\{ \rho_{\sigma,\lambda}\}$  {\it if} $a_{\sigma }=0$. 

\item If $\sigma=\sigma_j^{(2)}$ , $ j\in I^{(2)}$, then $(F_{\sigma }(\lambda), d_{\sigma  , \lambda }, \rho_{\sigma,\lambda})$ is a single 
GW($ \xi_{\mu, \lambda} , q)$-real tree if $y_j \leq \psi^{-1}(\lambda)-\psi^{-1}(\mu)$
 and it is simply 
the point tree $\{ \rho_{\sigma,\lambda}\}$  otherwise. 
\end{itemize}

\vspace{3mm}

\noindent 
We set 
$$  S_{\mu, \lambda }= \left\{ \sigma \in S_{\mu} \; :\; F_{\sigma }(\lambda) \neq \{\rho_{\sigma,\lambda}\} \right\}  $$
and 
$$ T' =
T  \, \circledast_{\sigma \in S_{\mu , \lambda} \setminus{\{ \rho_{\sigma,\lambda} \}}} 
(\sigma , F_{\sigma } (\lambda)  )  \quad {\rm and} \quad 
F'= T' \, \ast (\rho , F_{\rho}(\lambda) ).$$
\noindent
The following lemma implies that 
\begin{equation}
\label{identifii}
\overline{T}'\overset{(d)}{=} \overline{Q}_{\mu , \lambda} (T) \quad {\rm and} \quad 
\overline{F}' \overset{(d)}{=}\overline{Q}^a_{\mu , \lambda} (T) \; . 
\end{equation}

\begin{lemma}
\label{connection} 
Assume that $\psi'(0+)$ is finite. 
\begin{enumerate}
\item[{\rm(I)}] Let $E$ be a connected component of $T\backslash{({\rm Br}(T)\cup \{ \rho\})}$ (an edge of $T$). 
Then, $E\cap S_{\mu, \lambda }$ is 
a Poisson point process with intensity 
$$ \left( \psi'( \psi^{-1}(\lambda))
-\psi'( \psi^{-1}(\mu)) \right) {\bf 1}_{E} (\sigma ) \ell_T (d \sigma) .$$

\item[{\rm(II)}] Conditionally on $S_{\mu, \lambda }$, the random real forests  
$(F_{\sigma } (\lambda), d_{\sigma , \lambda}, \rho_{\sigma,\lambda})$ , $\sigma \in S_{\mu, \lambda }$, are independent. Moreover, for any 
$\sigma \in S_{\mu, \lambda }\backslash\{\rho\}$, the forest $F_{\sigma } (\lambda)$ consists of a random number $N_{\sigma }(\lambda)$  of independent 
GW($ \xi_{\mu, \lambda} , \psi' (\psi^{-1} (\lambda)) $)-real rooted trees, whose conditional distribution is given by 
$$ \bP \left( N_{\sigma }(\lambda)=k \; |  \; S_{\mu, \lambda } \right) =\nu^{\mu, \lambda}_l (k)  \; , \; k\geq 1 ,$$
where $l={\rm n}(\sigma , T)-1$.  
\end{enumerate}
\end{lemma}
{\bf Proof:} Set $S^{1,2}_{\mu}=S^{1}_{\mu} \cup S^{2}_{\mu}$ and let $M$ be the measure on $\bT$ given by 
$$ M (d\overline{\cT})=  
2\beta (\psi^{-1}(\lambda)-\psi^{-1}(\mu)) \Delta_{\mu, \lambda} (d\overline{\cT}) +
\int_{(0, \infty)}\Pi (dx) \, x e^{-\psi^{-1}(\mu)x} \Delta^x_{\mu, \lambda}  (d\overline{\cT}) .$$
An easy computation implies that
\begin{eqnarray*} 
M (\overline{\cT} \neq \{ \rho \} )&=&  2\beta (\psi^{-1}(\lambda)-\psi^{-1}(\mu)) + \int_{(0, \infty)}\Pi (dx) x e^{-\psi^{-1}(\mu)x} 
\left( 1-e^{-(\psi^{-1}(\lambda)-\psi^{-1}(\mu) )x} \right) \\
&=& \psi'( \psi^{-1}(\lambda))
-\psi'( \psi^{-1}(\mu)) .
\end{eqnarray*}
If we set $\widetilde{M}= 
M(\; \cdot \;  | \,\overline{\cT} \neq \{ \rho \} )$, then standard 
results on Poisson point processes imply that 
$$ \left\{ (\sigma ,F_{\sigma }(\lambda) )   \; ,\;  \sigma \in S^{1,2}_{\mu} \; :\;   F_{\sigma }(\lambda)\neq \{ \rho \} 
\right\}$$
is a Poisson point process with intensity 
$$ \left( \psi'( \psi^{-1}(\lambda))
-\psi'( \psi^{-1}(\mu)) \right) \ell_T (d\sigma) \otimes \widetilde{M}(d\overline{\cT}) .$$ 
This implies the first point of the lemma.

  Now, observe that if $\overline{\cT}$ has distribution  
$\widetilde{M}$, then $\overline{\cT}$ is obtained by pasting at the root $N$ independent copies of 
GW($ \xi_r , q)$-real rooted trees, where $\xi_r=\xi_{\mu,\lambda}$, $q$ is given by (\ref{qandc}), and the distribution of $N$ is given first by  
$\widetilde{M}(N=0)=0$ and for any $k\geq 1$ by 
\begin{eqnarray*}
\left( \psi'( \psi^{-1}(\lambda))
-\psi'( \psi^{-1}(\mu)) \right)\widetilde{M}(N=k) &=& 2\beta \, c \, {\bf 1}_{\{ k=1 \}} \;  +
\\ & \, &  
\int_{(0, \infty)}\Pi (dx) x e^{-\psi^{-1}(\mu)x} e^{-cx}(cx)^k/ k! \\
&=& (-1)^{k+1} \psi^{(k+1)} (\psi^{-1}(\lambda)) \, c^k/k! \; \; . 
\end{eqnarray*}
Accordingly, $\widetilde{M} (N=k)=\nu^{\mu, \lambda}_1(k)$ , $k\geq 0$, which implies the second part of the lemma in the 
$\sigma \in S^{1,2}_{\mu}\cap S_{\mu , \lambda }$ case.

   It remains to consider $\sigma \in {\rm Br}(T)\cup \{ \rho\}$: in that case the forest 
$F_{\sigma }(\lambda)$ is composed of $N_{\sigma }$ independent random GW($ \xi_r , q)$-real rooted trees, where 
$N_{\sigma }$ is a mixture of Poisson random variables whose distribution is given  
for any $k\geq 0$ by :
\begin{eqnarray*}
\bP \left( N_{\sigma }= k \; | \; S_{\mu, \lambda } \right) &=& \bE \left[ e^{-ca_{\sigma }}\frac{(ca_{\sigma })^k}{k!}\right]  \\
  &=& 
\frac{2\beta{\bf 1}_{\{ l=2\}}}{|\psi^{(2)} (\psi^{-1}(\mu))|} +
\frac{1}{|\psi^{(l)}(\psi^{-1}(\mu))|}\int_{(0, \infty)} \Pi (dx)x^{k+l} c^k e^{-x(c+\psi^{-1}(\mu))}/k! \\
&=& \nu^{\mu, \lambda}_{l}(k) 
\end{eqnarray*}
(here again $l={\rm n}(\sigma , T)-1$). This completes the proof of the lemma. \cqfd

\begin{remark}
\label{loifsigmulam}
Deduce from the definition of the $F_{\sigma }(\lambda)$'s that the sets of random variables 
$$ \cP_1 (\lambda)=\{ (\sigma ,\overline{F}_{\sigma }(\lambda) ) \; , \; \sigma \in 
S^1_{\mu} \cap  S_{\mu, \lambda }
 \} \; , \; 
 \cP_2 (\lambda)=\{ (\sigma ,\overline{F}_{\sigma }(\lambda) ) \; , \; \sigma \in 
S^2_{\mu} \cap  S_{\mu, \lambda }
 \}$$
and $\cP_3 (\lambda)=\{ (\sigma ,\overline{F}_{\sigma }(\lambda) ) \; , \; \sigma \in 
 {\rm Br}(T)\cup \{ \rho\}\} $ are independent. Their 
distributions are given as follows:
\begin{enumerate}
\item[(i)] $ \cP_1 (\lambda)$ is a Poisson point process on $T\times \bT$ with intensity measure
$$ \ell_{T} (d\sigma) \otimes \int_{(0, \infty)} \Pi (dr)re^{-r\psi^{-1}(\mu)} \Delta^r_{\mu, \lambda} 
(d\overline{\cT} \, \cap \{ \overline{\cT}\neq\overline{\{\rho\}}\}) $$
(Recall that $\overline{\{\rho\}}$ stands for the isometry class of the point tree). 

\item[(ii)] $ \cP_2 (\lambda)$ is a Poisson point process on $T\times \bT$ with intensity measure
$$ 2\beta \ell_{T} (d\sigma) \otimes (\psi^{-1}(\lambda)-\psi^{-1}(\mu)) \, \Delta_{\mu, \lambda}(d\overline{\cT})$$

\item[(iii)] For every $\sigma\in  {\rm Br}(T)$ , $\overline{F}_{\sigma }(\lambda)$ has distribution 
$$ \int_{[0, \infty)} \eta_{\mu , l} (dr) \, \Delta^r_{\mu, \lambda} (d\overline{\cT})\, ,$$
where $l={\rm n} (\sigma ,F_{\sigma }(\lambda) )-1$ and $\overline{F}_{\rho}(\lambda)$  
has distribution $\Delta^a_{\lambda,\mu } (d \overline{\cT})$. 
\end{enumerate}
\end{remark}

\begin{remark}
\label{numbtrees}
Fix $r>0$. Denote the total number of trees added on $B_T(\rho,r)$ by 
$N_{\mu , \lambda} (r)$. Then, note that  
$$ N_{\mu , \lambda} (r)= \sum_{\sigma \in  S_{\mu, \lambda }\cap B_T(\rho , r)}( {\rm n} (\rho_{\sigma,\lambda} ,F_{\sigma }(\lambda) )-1) . $$
Then, conditionally on $\cP_1$, $\cP_2$ and $\cA_{\mu}$, $N_{\mu , \lambda} (r)$
is distributed as a Poisson random 
variable with parameter $c A_{\mu} (r) $ where  
$$  A_{\mu} (r)= 2\beta \ell_T(B_T(\rho , r)) + 
\sum_{ \sigma \in S'_{\mu}\cap B_T(\rho , r)} a_{\sigma } (\mu ).  $$
Thus,
\begin{equation}
\label{equnombre}
 \E \left[ s^{ N_{\mu , \lambda} (r)}\right]= \E \left[ \exp (-cA_{\mu} (r) (1-s) )\right]. 
\end{equation}
\cqfd \end{remark}

\section{The L\'evy forest.}
\label{mainlevy}

\subsection{Construction of the L\'evy forest.}
\label{levyforest}

  In this section, we study the increasing limit of $\cF_\lambda$ as $\lambda\rightarrow\infty$, 
and properties of the limit. Let us consider an $(a, \psi)$-growth process started at the discrete 
tree with edge lengths $T$ denoted by 
$(\cF_{\lambda }, d_{\lambda}, \rho)$, $\lambda \in [0, \infty)$. Set $\cF_{\infty}=\bigcup \cF_{\lambda}$ and define 
a metric $d$ on $\cF_{\infty}$ by $d(\sigma , \sigma')
=d_{\lambda} (\sigma , \sigma')$ if 
$\sigma , \sigma'\in \cF_{\lambda}$. We denote by $(\cF, d)$ the completion of 
$(\cF_{\infty}, d)$. 

\begin{theorem}
\label{mainresult} Assume that (\ref{posextinct}) holds. Almost surely, $(\cF, d, \rho)$ is a locally compact rooted real tree and 
$$ \delta (\overline{\cF}_{\lambda} , \overline{\cF}) 
 \build{\longrightarrow}_{\lambda \to\infty}^{\; } \; 0 \; .$$
\end{theorem}

\begin{remark}
\label{nonloccomp}
If (\ref{posextinct}) does not holds, then the popupation may become extinct but in an infinite time and therefore the underlying genealogical 
tree cannot be locally compact. 
\end{remark}
\noindent 
\noindent {\bf Proof :} Thanks to Lemma \ref{evaconv}, it is sufficient to prove that for any $r\in (0, \infty)$ a.s. 
the collection of closed balls  $( B_{\cF_{\lambda}} (\rho , r) ; \lambda \geq 0)$ is Cauchy when $\lambda $ 
goes to infinity with respect to the Hausdorff distance $d_{{\rm Haus}} $ on compact sets of $(\cF, d)$. 
Set 
$$ \Xi_{\mu, \lambda} (r):=d_{{\rm Haus}} \left( 
 B_{\cF_{\mu}} (\rho , r)\, , \, B_{\cF_{\lambda}} (\rho , r)
\right). $$
 Since $\Xi_{\mu, \lambda} (r)$ is non-decreasing in $\lambda$ and non-increasing in $\mu$, we only have to prove that
for any $t>0$
\begin{equation}
\label{hauscauch}
\lim_{\mu \rightarrow \infty} \; \sup_{\lambda \geq \mu} \; 
\bP \left( 
\Xi_{\mu, \lambda} (r) \leq t \right)\; = \; 1 \; . 
\end{equation} 
We first need to introduce some notation: let $(\cT_i^o , i\in I)$ be the connected components of the open set 
$\cF_{\lambda} \backslash T $ in $\cF_{\lambda}$. Denote by $\sigma_i$
the vertex 
of $T$ on which $\cT_i^o$ is grafted  
and set $\cT_i =\cT_i^o \cup \{ \sigma_i  \}$. Then, the $(\cT_i, d, \sigma_i)$'s are 
compact rooted real trees and 
$$ \cF_{\lambda } =T \circledast_{i\in I} (\sigma_i ,\cT_i ) . $$
Let $\mu \in [0, \lambda]$ and let $i\in I$. Set $\cT_i'=\cT_i \cap \cF_{\mu }$ and denote by 
$(\cT_{i,j}^o , j\in J(i))$ the connected components of $\cT_i \backslash \cT_i'$. Denote by 
$\sigma_{i,j}$ the vertex of $ \cT_i'$ on which $\cT_{i,j}^o$ is grafted and set 
$\cT_{i,j} =\cT_{i,j}^o \cup \{ \sigma_{i,j}  \}$. Clearly, the $(\cT_{i,j}, d, \sigma_{i,j})$'s are 
compact rooted real trees. Observe that 
$$ \cF_{\lambda} \backslash \cF_{\mu } = \bigcup_{\substack{ i\in I, \\ j\in J(i)}} \cT_{i,j}^o 
\; \; {\rm and } \; \;  \cF_{\mu } =T \circledast_{i\in I} (\sigma_i ,\cT_i' ).$$ 
Thus 
$$\cF_{\lambda }=\cF_{\mu }\circledast_{\substack{ i\in I, \\ j\in J(i)}}(\sigma_{i,j} ,\cT_{i,j} ) .$$ 
To simplify notations, we set 
$$h_{i,j}:= h(\cT_{i,j})=\sup \{ d( \sigma_{i,j},\sigma) \, , \, \sigma \in \cT_{i,j}\}$$ 
and  $I(r):= \{ i \in I \, : \, d(\rho , \sigma_i) \leq r   \}$. Then, the previous observations imply
\begin{equation}
\label{ineqdelta}
\Xi_{\mu, \lambda}(r) \leq \max \{ h_{i,j} \; , \; \ i\in I (r),\,   j\in J(i) \}. 
\end{equation}
Now deduce from Proposition \ref{chaining} 
$$ \left( \overline{Q}_{0 ,\mu}^a (T), \, 
\overline{Q}_{\mu, \lambda}^a (Q_{0 ,\mu}^a (T) ) \right) \overset{(d)}{=}\left( \overline{\cF}_{\mu }, 
\overline{\cF}_{\lambda }\right). $$
So if we set 
$$ N_{0, \lambda} (r)= \# I(r) \; , \; \; N_i =\# J(i) \;\;  {\rm and} \; \;  N_{\mu, \lambda} (r)= \sum_{i\in I(r)} N_i \; , $$
then we get the following: 
\begin{enumerate}
\item[(a)] Conditionally on $N_i$ , $i\in I$, the trees $ \overline{\cT}_{i,j}$, $ \ i\in I (r),\,   j\in J(i)$, are independent 
GW($\xi_{\mu, \lambda}, \psi'(\psi^{-1} (\lambda))$)-real trees. 

\item[(b)] Conditionally on $N_{0, \lambda} (r)$, $(\overline{\cT}_i',\overline{\cT}_i)$, 
$i \in I(r)$ are  i.i.d. pairs of trees distributed as 
$( \overline{\cT}',\overline{\cT})$ where $\cT$ is a  GW($\xi_{0, \lambda}, \psi'(\psi^{-1} (\lambda))$)-real tree and 
$\cT'$ is obtained from $\cT$ as the black subtree resulting from a $(1-\mu /\lambda)$-Bernoulli leaf colouring. 
\end{enumerate}
Thus, if we set $K_{\mu, \lambda}(t)=\bP( \max \{ h_{i,j}, i\in I (r), j\in J(i) \} \leq t)$, we deduce from (a)
\begin{equation}
\label{maxheight}
K_{\mu, \lambda}(t)= \bE \left[ \exp (- v_{\mu, \lambda}(t) N_{\mu, \lambda} (r))\right] , 
\end{equation}
where $\exp (- v_{\mu, \lambda}(t)):=\bP(h(\cT'') \leq t) $ and $\cT''$ is a 
GW($\xi_{\mu, \lambda}, \psi'(\psi^{-1} (\lambda))$)-real tree. Then, (\ref{heightharris}) applied to 
$\varphi_{\mu, \lambda}$ and a simple change of variable imply that $ v_{\mu, \lambda}(t)$ satisfies the following equation
\begin{equation}
\label{integequa}
\int_{(\psi^{-1} (\lambda)-\psi^{-1} (\mu))(1-e^{- v_{\mu, \lambda}(t)})}^{\psi^{-1} (\lambda)-\psi^{-1} (\mu)} 
\; \frac{dx}{\psi( \psi^{-1} (\mu) +x) -\mu } =t .
\end{equation}
We now need to compute the distribution of $ N_{\mu, \lambda} (r)$ and accordingly the distribution of the $N_i$, $i\in I(r)$.   
If $( \cT',\cT)$ are as in (b), then denote by $M$ the number of red trees grafted on $\cT'$. Note that 
$M$ is possibly equal to $1$ if 
$\cT'$ is reduced to the point tree $\{\rho\}$, that is if $\cT$ is completely red. Set $\kappa (s)=\bE[s^{M}]$. According 
to (\ref{numbredtrees}), $\kappa$ satisfies
$$ \varphi_{0,\lambda} \left( \kappa (s)\right)- \kappa (s)= 
\varphi_{0,\lambda} ( sg(p) )- sg(p) - \left( 
\varphi_{0,\lambda}( g(p)) - g(p) \right), $$
where we recall that $1-p=\mu/ \lambda$ and 
$$g(p)=\bE[p^{\# {\rm Lf} (\cT)}]= 
1-\frac{\psi^{-1}(\mu)-\gamma}{\psi^{-1}(\lambda)-\gamma} \; .$$
A straightforward computation implies: 
\begin{equation}
\label{kappacomp}
(\psi^{-1}(\lambda)-\gamma)(1-\kappa (s))\; = 
\psi^{-1} \left[ \psi \left( 
(\psi^{-1} (\lambda)-\psi^{-1} (\mu))(1-s)+ \psi^{-1} (\mu)
\right)  -\mu \right] -\gamma.
\end{equation}
Then, by (a) and (b) we get:
$$ \bE \left[ \exp (- v_{\mu, \lambda}(t) N_{\mu, \lambda} (r))\right] =\bE \left[ 
\kappa (e^{- v_{\mu, \lambda}(t)})^{ N_{0, \lambda} (r)} \right] .$$
Recall the notation of Remark \ref{numbtrees}: we take here $\mu=0$ and therefore we set  
$$ A_0 (r) =  2\beta \ell_{T} (B_{T}(\rho , r)) \; + 
\sum_{ \sigma \in S'_{0}\cap B_{T}(\rho , r)} a_{0 } (\sigma ) . $$
Thus, by (\ref{equnombre})
$$K_{\mu, \lambda}(t)=\bE \left[ \exp( - A_0 (r)\, 
(\psi^{-1}(\lambda)-\gamma) (1-\kappa (e^{- v_{\mu, \lambda}(t)}))  )\right].$$
Deduce from (\ref{integequa}) that 
$$ \lim_{\lambda \rightarrow \infty} (\psi^{-1} (\lambda)-\psi^{-1} (\mu))(1-e^{- v_{\mu, \lambda}(t)})=w_{\mu} (t) \; < \; \infty $$ 
which satisfies 
$$ \int_{w_{\mu} (t)}^{\infty}
\frac{dx}{\psi( \psi^{-1} (\mu) +x) -\mu } =t  \; .$$
Notice that here we use (\ref{posextinct}).Thus 
$$  \lim_{\lambda \rightarrow \infty}K_{\mu, \lambda}(t)=
\bE \left[ \exp \left( -A_0  (r) \psi^{-1} \left( \psi \left( 
w_{\mu} (t)  + \psi^{-1} (\mu)  \right)  -\mu \right) -\gamma \right) \right].$$
Finally observe that 
$$ t=\int_{w_{\mu} (t)}^{\infty}
\frac{dx}{\psi( \psi^{-1} (\mu) +x) -\mu }= \int_{ \psi(  \psi^{-1} (\mu) + w_{\mu} (t) )- \mu  }^{\infty}
\frac{dy}{y \psi'( \psi^{-1} (\mu  +y) )}, $$
which implies 
$$\lim_{\mu \rightarrow \infty} \psi(  \psi^{-1} (\mu) + w_{\mu} (t) )- \mu =0 $$
by dominated convergence. Thus 
$$\lim_{\mu \rightarrow \infty}\lim_{\lambda \rightarrow \infty}K_{\mu, \lambda}(t)=1 . $$
It proves (\ref{hauscauch}), which completes the proof of the theorem. \cqfd

\begin{remark}
\label{ellunforestlim}
Assume that the $(a, \psi)$-L\'evy growth process 
$(\cF_{\lambda }, || \cdot ||_1 , 0)$, $\lambda \in [0, \infty)$ is $\btl$-valued. The proof 
actually implies that a.s. 
$$ \dl (\cF_{\lambda }, \cF) \build{\longrightarrow}_{\lambda \to\infty}^{\; } \; 0 \; .$$
\end{remark}

\begin{notation}
\label{levynot}
\begin{itemize} 
\item The random locally compact rooted real tree obtained as a limit of an $(a, \psi)$-L\'evy growth process starting at 
$T$ is called an {\it $(a, \psi)$-L\'evy forest starting at $T$} and we shall sometimes denote such a random tree by the symbol
$Q^a_{0, \infty} (T)$. We also denote by $\overline{Q}^a_{0, \infty} (T)$ its isometry class. 

\item We call {\it $(a, \psi)$-L\'evy forest} the random tree $Q^a_{0, \infty} (\cF_0)$, where 
$\cF_0$ is a GW($\xi_0, \psi' (\gamma) , a\gamma$)-real forest that is independent of the random variables used to define 
the growth process. We denote by $P^a (d\overline{\cT})$ the distribution on $\bT$ of 
$\overline{Q}^a_{0, \infty} (\cF_0)$. 

\item Let $\mu \geq 0$. Observe that $\psi_{\mu}$ satisfies the assumptions of Theorem 
\ref{mainresult}. We denote the limit of the $(a, \psi_{\mu})$-growth process started at $T$ by the symbol 
$Q^a_{\mu, \infty} (T)$. 

\item Observe that $0$ is the only root of 
$\psi_{\mu}(x)=0$. So an $(a, \psi_{\mu} )$-L\'evy forest is the limit of an $(a, \psi_{\mu} )$-growth process
started a the tree reduced to a point. We denote the distribution of the isometry class of an $(a, \psi_{\mu} )$-L\'evy forest by 
$P^a_{\mu} (d\overline{\cT})$. If 
$\gamma >0$, then $P^a_0 \neq P^a$. 

\item We shall also consider the following random trees. Let $\cT_{\mu }$ be a GW($\xi_{\mu}, \psi' (\psi^{-1}(\mu))  $)-real 
tree and let $\cT_{\mu , \lambda }$ be GW($\xi_{\mu , \lambda }, \psi' (\psi^{-1}(\lambda))  $)-real 
tree. We denote by $P_{\mu}(d\overline{\cT})$ the distribution on $\bT$ of $\overline{Q}^{a=0}_{\mu, \infty} (\cT_{\mu})$ and 
we denote by  $P_{\mu, \lambda}(d\overline{\cT})$ the distribution on $\bT$ of 
$\overline{Q}^{a=0}_{\lambda, \infty} (\cT_{\mu, \lambda})$. 
Now observe that $P^0_{\mu}=\delta_{\overline{\{\rho\}}}$ and thus $P^0_{\mu}\neq P_{\mu}$ 
(recall that  $\overline{\{\rho\}}$ stands for the isometry class of the tree reduced to a point). 
\end{itemize}
\end{notation}

Let us end this subsection by two useful observations: first note that
Proposition \ref{chaining} combined with Theorem \ref{mainresult} with
$\psi_{\mu}$ imply that for any 
discrete tree with edge lengths $T$, we have
\begin{equation}
\label{chaininginf}
\overline{Q}^a_{\lambda, \infty} \left( Q^a_{\mu , \lambda} (T)\right) \quad 
\overset{(d)}{=} \quad \overline{Q}^a_{\mu, \infty}(T) \;  
\end{equation}
(here the extra random variables used to define $Q^a_{\lambda, \infty}$ are chosen independent of $ Q^a_{\mu , \lambda} (T)$). 
Then, recall notation $\Delta_{\mu, \lambda}^a$ from the previous section. Apply Theorem 
\ref{mainresult} with $\psi_{\mu}$ to get  
\begin{equation}
\label{weakdelta}
\Delta_{\mu, \lambda}^a
\build{\longrightarrow}_{\lambda \to\infty}^{\; } \;  P^a_{\mu} \; 
\end{equation}
weakly in the space of probability measures on $\bT$.

\subsection{The mass measure.}
\label{leb}

   Let $a\geq 0$ and let $(\cF_{\lambda}; \lambda \geq 0)$ be an $(a, \psi)$-L\'evy growth process. 
{\it We assume that the $\cF_{\lambda}$ are embedded in $\ell_1(\bN)$} and we denote by $\cF$ the 
limit of this growth process in $\btl$. We also denote 
by ${\bf m}_{\lambda}$ the empirical distribution of the leaves of ${\rm Lf}(\cF_{\lambda})$: 
\begin{equation}
\label{finbirth}
{\bf m}_{\lambda}= \sum_{\sigma \in {\rm Lf}(\cF_{\lambda}) \backslash {\rm Lf}(T)} \delta_{\sigma} \; . 
\end{equation}
\begin{theorem}
\label{meslebe}
There exists a random measure ${\bf m}$ on $\ell_1 (\bN)$ such that 
\begin{enumerate}
\item[{\rm(i)}] Almost surely the convergence 
$$ \lambda^{-1}{\bf m}_{\lambda} \; \build{\longrightarrow}_{\lambda \to\infty}^{\; } \; {\bf m} $$
holds for the vague topology of Radon measures on  $\ell_1 (\bN)$; 

\item[{\rm(ii)}]  Almost surely the topological support of ${\bf m} $ is $\cF$; 

\item[{\rm(iii)}] Let $\cP=\{ (\sigma_j , U_j) \;  , \; j\in J  )\}$ be a Cox process on $\ell_1 (\bN) \times [0, \infty )$ with 
random intensity ${\bf m} (d\sigma)\otimes du $. For any $\lambda \geq 0$ denote by 
$\cF_{\lambda}'$  the subtree of $\cF$ spanned by 
$0$ and the set of vertices $\{ \sigma_j  \;  ;  \; j\in J \; , \; U_j \leq \lambda \}$: 
$$ \cF_{\lambda }'= \bigcup \{ \lgeo 0 , \sigma_j \rgeo  \; ; \;  j\in J \; , \; U_j \leq \lambda \} .$$ 
Then, 
$$  (\overline{\cF}_{\lambda}; \lambda \geq 0) \overset{(d)}{=}(\overline{\cF}_{\lambda}'; \lambda \geq 0) .$$
\end{enumerate}
\end{theorem} 

\begin{remark}The measure ${\bf m}$ is concentrated on the leaves of $\cF$ since by definition 
${\bf m}(\cF_\infty)=0$ and since $\cF\backslash{\rm Lf}(\cF)\subset\cF_\infty$.
\end{remark}

\noindent{\bf Proof:} Let us prove (i). By standard density arguments, it is sufficient to prove that for any non-negative 
continuous function $f$ on $\ell_1 (\bN)$ with compact support there exists a non-negative finite random variable 
${\bf m} (f)$ such that 
we a.s. have 
\begin{equation}
\label{singlef}
\lambda^{-1}<{\bf m}_{\lambda}, f> \; \build{\longrightarrow}_{\lambda \to\infty}^{\; } \; {\bf m} (f). 
\end{equation}
Fix $\mu \geq 0 $. We denote by $\cT_i^o$ , $ i\in I(\mu )$, the connected components of $\cF \backslash \cF_{\mu}$  and 
we denote by $\sigma_i$ the vertex of $\cF_{\mu}$ on which $\cT_i^o$ is grafted and we set 
$$ \cT_i = \{ \sigma_i\} \cup \cT_i^o  \quad {\rm and } \quad   \cT_i  (\lambda)=  \cT_i \cap \cF_{\lambda}  \; , \; \lambda \geq \mu 
\; .$$
Set for any $\lambda \geq \mu$
$$ {\bf m}_{\lambda}^{ \cT_i }= \sum_{\sigma \in {\rm Lf}(\cT_i (\lambda))} \delta_{\sigma} , $$
with the conventions that if $\cT_i (\lambda)= \{ \sigma_i\}$ 
then ${\rm Lf}(\cT_i (\lambda))= \emptyset$
and $ {\bf m}_{\lambda}^{ \cT_i }=0$. Then, for any $\lambda_2 \geq \lambda_1 \geq \mu $
$$ \lambda_2^{-1}<{\bf m}_{\lambda_2}, f>- \lambda_1^{-1}<{\bf m}_{\lambda_1}, f> = T_1 + T_2 + T_3, $$
where 
$$ T_1= \left( \lambda_2^{-1}-\lambda_1^{-1} \right) <{\bf m}_{\mu}, f> \; , $$
$$ T_2 = \sum_{i\in I(\mu)} < \lambda_2^{-1}{\bf m}^{ \cT_i }_{\lambda_2} -\lambda_1^{-1}
{\bf m}^{ \cT_i }_{\lambda_1} \; , \;  f-f(\sigma_i)> \; ,  $$ 
$$ T_3 = \sum_{i\in I(\mu)} f(\sigma_i) \left(  < \lambda_2^{-1}{\bf m}^{ \cT_i }_{\lambda_2}> -<\lambda_1^{-1}
{\bf m}^{ \cT_i }_{\lambda_1}> \right) . $$
We set 
$$ M (\lambda):= \lambda^{-1}< {\bf m}_{\lambda} , \un_{B(0, r)}>  = 
 \lambda^{-1} \# \{ \sigma \in {\rm Lf} (\cF_{\lambda} ) \; : \; ||\sigma ||_1 \leq r\} , $$
where $r$ is such that $f(\sigma )=0$ if $ ||\sigma ||_1 >r $. We
 also define 
for any $\lambda \geq \mu$
$$ M_{\mu , f} (\lambda):=  \sum_{i\in I(\mu)} f(\sigma_i)  \lambda^{-1} < {\bf m}^{ \cT_i }_{\lambda}>  \; . $$
\begin{lemma}
\label{martgmass}
There exist two finite random variables  $M_{\mu , f} (\infty)$ and  $M (\infty)$ such that a.s. 
$$   M_{\mu , f} (\lambda)\build{\longrightarrow}_{\lambda \to\infty}^{\; } M_{\mu , f} (\infty) \quad  {\rm and } \quad 
 M (\lambda)\build{\longrightarrow}_{\lambda \to\infty}^{\; } M (\infty) .$$
\end{lemma}
{\bf Proof of the lemma}: For any $\lambda \geq \mu$, denote by $\cG_{\lambda}$ the sigma-field generated by $\cF_{\mu}$, 
the random variables $(\sigma_i, \overline{\cT}_i (\lambda'); \lambda'
\geq \lambda )$, $i\in I(\mu)$, and the $\bP$-null sets. Set also 
$$ I(\mu , \lambda)= \{ i\in  I(\mu ) \; : \; \cT_i (\lambda)\neq  \{ \sigma_i\} \} .  $$
Clearly for any 
$\lambda \geq \lambda' \geq \mu$, we have $\cG_{\lambda} \subset \cG_{\lambda'}$ and $ M_{\mu , f} (\lambda)$ is 
$\cG_{\lambda}$-measurable. 
Moreover the random variable $< {\bf m}^{\cT_i }_{\lambda}>$ only depends on $\cG_{\lambda}$ via 
$\overline{\cT}_i (\lambda)$. Then, observe 
that for any $\lambda \geq \lambda' $ conditionally on $\cF_{\mu}$ and on $ I(\mu , \lambda)$, the 
trees $( \overline{\cT}_i ( \lambda') , \overline{\cT}_i ( \lambda))$ , $i \in I(\mu , \lambda)$, are independent and distributed as 
 $( \overline{\cT}_b , \overline{\cT})$  where $\cT$ is a GW($\xi_{\mu, \lambda}, \psi' (\psi^{-1}(\lambda))$)-real tree 
and where $\cT_b$ is the black 
subtree of $\cT$ resulting from a $(1-(\lambda' -\mu)/(\lambda -\mu))$-Bernoulli leaf colouring. 
Therefore,  
conditional on $\overline{\cT}$, $\# {\rm Lf} (\cT_b)$ has a binomial distribution 
with parameters $\# {\rm Lf} (\cT)$ and  $(\lambda' -\mu)/(\lambda -\mu) $. Accordingly 
$$ \bE \left[ \left. \# {\rm Lf} (\cT_b) \right| \overline{\cT}  \right] = 
\frac{\lambda' -\mu}{\lambda -\mu} \# {\rm Lf} (\cT )  . $$

Then, deduce from the latter observations that

\begin{eqnarray*}
\bE \left[ \left.  M_{\mu , f} (\lambda')
\right|  \cG_{\lambda}  \right] &=&   \sum_{i\in I(\mu)}   f(\sigma_i)
\lambda^{-1} 
\bE \left[ \left.
< \! {\bf m}^{ \cT_i }_{\lambda'} \!>  \right| \overline{\cT}_i (\lambda)  \right] \\
&=&   \frac{\lambda' -\mu}{\lambda -\mu}  \, . \, 
\frac{\lambda}{\lambda'} \, M_{\mu , f} (\lambda) .
\end{eqnarray*}

Thus, $M=(\frac{\lambda }{\lambda -\mu  } M_{\mu , f} (\lambda) ; \lambda \geq \mu  )$ is a 
non-negative backward martingale with respect to $ (\cG_{\lambda}; \lambda \geq \mu  )$. 
A similar result holds for 
$ ( \frac{\lambda }{\lambda -\mu }M (\lambda); \lambda \geq
\mu)$. Therefore, these two backward martingales converge to two
limits in $[0, \infty]$ denoted by resp. $M_{\mu , f} (\infty)$ and $M (\infty)$. Since 
$\lambda /(\lambda -\mu) $ converges to $1$ when $\lambda $ goes to infinity, it implies the two convergences of the lemma. It remains 
to show that these two limiting random variables are a.s. finite.

  To that end, observe that 
\begin{equation}
\label{inega}
 M_{\mu , f} (\lambda) \leq \frac{||f ||_{\infty}}{\lambda} \sum_{i\in I(\mu)} 
\un_{[0,r]} (||\sigma_i ||_1) \; \# {\rm Lf} (\cT_i (\lambda)) .
\end{equation}
Then, recall that conditionally on $\cF_{\mu}$ and $ I(\mu, \lambda)$, the trees $\overline{\cT}_i (\lambda)$ , 
$i\in  I(\mu, \lambda)$ are independent with the same distribution as $\overline{\cT}$.
Fix $\theta >0$. Use Remark \ref{psishift}, take $s=e^{-\theta / \lambda }$ and replace $\psi$ by $\psi_{\mu}$ in (\ref{gequa}),
to get 
$$
g (e^{-\theta / \lambda })= \bE \left[  e^{-\frac{\theta}{\lambda} \# {\rm Lf}(\cT)}  \right] 
= 1-\frac{\psi_{\mu}^{-1} ((1-e^{-\theta/\lambda})(\lambda -\mu))}{\psi_{\mu}^{-1} (\lambda -\mu)} \; .
$$
Set $N_{\mu, \lambda} (r)= \# \{ i\in I(\mu, \lambda) \; : \;
||\sigma_i||_1 \leq r \}$. 
Then the previous observation implies
$$\bE \left[\exp\left(-\theta  \, \lambda^{-1} 
\sum_{ i\in I(\mu, \lambda)}  
\un_{[0,r]} (||\sigma_i ||_1) \; 
\# {\rm Lf}(\cT_i (\lambda ))\right)\right]
= \bE \left[ \left(g (e^{-\frac{\theta}{\lambda}})\right)^{N_{\mu, \lambda} (r)}  \right] \; .$$
Now use Remark \ref{numbtrees} to get 
$$ \bE \left[\left(g (e^{-\frac{\theta}{\lambda}})\right)^{N_{\mu, \lambda} (r)}  \right] =
\bE \left[ \exp ( -A_{\mu} (r)\psi_{\mu}^{-1} 
((1-e^{-\frac{\theta}{\lambda}})(\lambda -\mu))  ) \right] .$$
Thus, 
$$ \lim_{\lambda \rightarrow \infty} \bE \left[\left(g (e^{-\theta /\lambda})\right)^{N_{\mu, \lambda} (r)}  \right] = 
\bE \left[ e^{-A_{\mu} (r)\psi_{\mu}^{-1} (\theta )  } \right] .$$ 
Then, by (\ref{inega})
\begin{eqnarray*} 
 \bE \left[ e^{-\theta M_{\mu, f} (\infty)}  \right] &=& \lim_{\lambda \rightarrow \infty}
\bE \left[ \exp \left( -\theta M_{\mu, f} (\lambda )\right) \right]   \\
& \geq & \lim_{\lambda \rightarrow \infty} 
\bE \left[ \exp \left( -
\frac{\theta ||f ||_{\infty}  }{\lambda} \sum_{i\in I(\mu)} \un_{[0,r]} (||\sigma_i ||_1) \# {\rm Lf} (\cT_i (\lambda))
\right) \right] \\
 &=&\bE \left[ e^{-A_{\mu} (r)\psi_{\mu}^{-1} (\theta ||f ||_{\infty}  )  } \right] . 
\end{eqnarray*}
Since the right member of the last  inequality tends to $1$ when $\theta$ goes to $0$, so does the first member, which implies that 
$ M_{\mu, f} (\infty)$ is a.s. finite. A similar argument works for $ M (\infty)$. This completes the proof of the lemma. 
\cqfd

\vspace{3mm}

  Let us fix $\Omega'\subset \Omega$ such that $\bP(\Omega')=1$ and such that the following limits hold
$$ \lim_{\lambda \rightarrow \infty} \dl (\cF_{\lambda}, \cF)=0 \; , \quad 
 \lim_{\lambda \rightarrow \infty} M_{\mu, f} (\lambda )= M_{\mu, f} (\infty) \quad {\rm and } \quad 
\lim_{\lambda \rightarrow \infty} M  (\lambda )= M (\infty) . $$
We fix $\omega \in \Omega'$. Let $\epsilon >0$. For any $\eta >0$ we denote the modulus of uniform continuity of $f$ by 
$w(f, \eta):=\sup \{ |f(\sigma )-f(\sigma' ) |  ; \, ||\sigma -\sigma'||_1 \leq \eta  \}$.  
\begin{enumerate}

\item[(a)] We choose $\eta $ such that  
$$ 3  M (\infty) \, w(f, \eta) \; \leq \; \epsilon . $$

\item[(b)] We choose $\mu $ large enough such that 
$$\dl (\cF_{\mu}, \cF) \; \leq \; \eta $$

\item[(c)] Then, we choose $\lambda$ large enough such that for any $\lambda_1 , \lambda_2 \geq \lambda$ 
$$ \left| M_{\mu, f} (\lambda_1 )- M_{\mu, f} (\lambda_2 )  \right| <\epsilon  \quad , \quad 
 M(\lambda_1 )+ M (\lambda_2 ) \leq  3  M (\infty)$$
and 
$$ |T_1| =\left|  \lambda_2^{-1}- \lambda_1^{-1} \right| <{\bf m}_{\mu},f>  \; \leq \; \epsilon   .$$
\end{enumerate}

Then, by (b) we have
 $$ < \lambda_2^{-1}{\bf m}^{ \cT_i }_{\lambda_2} -\lambda_1^{-1}
{\bf m}^{ \cT_i }_{\lambda_1} \; , \;  f-f(\sigma_i)>  \; \leq  \;  w(f, \eta) \left( 
 \lambda_2^{-1}<{\bf m}^{ \cT_i }_{\lambda_2}> + \lambda_1^{-1} <{\bf m}^{ \cT_i }_{\lambda_1}>
\right) . $$
Thus, by (a) and (c)
\begin{eqnarray*} 
 |T_2| & \leq & \left(   M(\lambda_1 )+ M (\lambda_2 )  \right)  w(f, \eta) \\
& \leq & 3  M (\infty) \,  w(f, \eta) \; \leq \; \epsilon . 
\end{eqnarray*} 
Now observe that $T_3=  M_{\mu, f} (\lambda_1)- M_{\mu, f} 
(\lambda_2)$. By (c) we get $|T_3|\leq \epsilon$. Thus we have proved that 
for any $\omega \in  \Omega'$ and any $\epsilon >0$, we can find a sufficiently large $\lambda$ 
such that 
$$ \sup_{\lambda_1 , \lambda_2 \geq \lambda} \left|  \lambda_2^{-1}
<{\bf m}_{\lambda_2}, f>-\lambda_1^{-1}<{\bf m}_{\lambda_1}, f> \right| \; \leq \; 3 \epsilon , $$
which implies (\ref{singlef}) and then (i) of Theorem \ref{meslebe}.

\vspace{3mm}

  Let us prove (ii). To that end, set for any  $\theta \in [0, \infty)$ and any $\lambda>\mu$ 
$$ h_{\mu , \lambda} (\theta )= \frac{\psi_{\mu}^{-1}(\theta +\lambda -\mu )-\psi_{\mu}^{-1}(\theta )}{
\psi_{\mu}^{-1}(\lambda -\mu )} . $$
We need the following lemma 

\begin{lemma}
\label{loimui}
Conditionally on $\cF_{\mu}$ and on $I(\mu , \lambda)$, the random variables ${\bf m} (\cT_i)$ , 
$ i\in I(\mu , \lambda)$ are i.i.d.  and the Laplace transform of their conditional distribution is  
$h_{\mu , \lambda} $. 
\end{lemma}
{\bf Proof of the lemma}: Since for any $ i\in I(\mu , \lambda)$, ${\bf m} (\{ \sigma_i\})=0$, it is easy to check that 
\begin{equation}
\label{llimi}
 (\lambda')^{-1}<{\bf m}_{\lambda'}^{\cT_i}> \build{\longrightarrow}_{\lambda' \to\infty}^{\; } 
{\bf m} (\cT_i).
\end{equation}
Now observe that almost surely for $\lambda' \geq   \lambda  \geq \mu$, conditionally on $\cF_{\mu}$ and on $I(\mu , \lambda)$, 
the trees $\overline{\cT_i}(\lambda')$ , $i \in I(\mu , \lambda)$, are independent and distributed as 
$\overline{\cT}$ where $\cT$ stands for a GW($ \xi_{\mu , \lambda '} , \psi' (\psi^{-1} (\lambda')) $)-real tree conditioned 
on not being completely red after a $(1-\frac{\lambda-\mu }{\lambda'-\mu})$-Bernoulli leaf colouring. 
Denote by $\cT_b$ the black subtree resulting from such a colouring. An elementary computation based on Remark 
\ref{psishift} and (\ref{gequa}) implies that 
$$ \bE \left[ \left. s^{\# {\rm Lf (\cT) }}  \right| \cT_b \neq \{ \rho\}  \right]= 
\frac{\psi_{\mu}^{-1} \left(  (1-s)\lambda ' +s\lambda -\mu  \right) - 
\psi_{\mu}^{-1} \left(  (1-s)(\lambda '-\mu) \right)}{ \psi_{\mu}^{-1} \left( \lambda -\mu \right) } . $$
Take  $s=\exp (-\theta / \lambda' )$ and then observe that the right member converges to  
$h_{\mu , \lambda} (\theta )$ when $\lambda '$ goes to infinity. This completes the proof of the lemma. \cqfd

\vspace{3mm}

\noindent
{\bf End of the proof of the theorem:} Since $\psi_{\mu}^{-1} $ is concave, we get 
$$ h_{\mu , \lambda} (\theta ) \; \leq  \;  \frac{1}{\psi_{\mu}^{-1}(\lambda-\mu)} .\frac{\lambda -\mu }{\psi_{\mu}' 
(\psi_{\mu}^{-1}(\theta ))} 
\build{\longrightarrow}_{\theta \to\infty}^{\; } 0 . $$
Thus for any $\lambda>\mu$ ,  ${\bf m} (\cT_i) >0$ , $ i\in I(\mu , \lambda)$ a.s. It implies that a.s. ${\bf m} (\cT_i) >0$
for every  $i \in I(\mu)$. 
Then a.s. for every $\mu \geq 0 $ the topological support of ${\bf m}$ has a non-trivial 
intersection with each of the connected components of $\cF \backslash \cF_{\mu}$, which implies (ii). 

\vspace{3mm}

  Let us prove (iii). Since the process $ (\cF_{\lambda}'; \lambda \geq 0)$ is obviously Bernoulli leaf colouring consistent, we only 
have to prove that for a fixed $\mu>0$, we have
\begin{equation}
\label{onemarg}
\overline{\cF}_{\mu}' \overset{(d)}{=}\overline{\cF}_{\mu} . 
\end{equation}
Conditionally on $ (\cF_{\lambda}; \lambda \geq 0)$, let $V_{\sigma }$ , $ \sigma \in \bigcup_{\lambda \geq 0} 
{\rm Lf}(\cF_{\lambda})$ be 
i.i.d. $[0, 1]$-uniform random variables. Set for any $\lambda \geq \mu$
$$ \cN_{\mu,\lambda} = \sum_{\sigma \in {\rm Lf} (\cF_{ \lambda})} \un_{[0, \mu / \lambda]} (V_{\sigma}) \, 
\delta_{(\sigma , \lambda  V_{\sigma})} .  $$
Denote by  $\cM (\ell_1 (\bN))$ the set of  Radon measures of $\ell_1 (\bN)$ and equip it with a metric compatible with the 
vague topology. Let $K$ be a measurable non-negative function on $\bT_{\ell_1} \times \cM (\ell_1 (\bN))$ and let 
$f$ be a non-negative continuous function on $\ell_1 (\bN) \times [0, \mu ]$ with compact support. Set 
$$\cE_{\lambda}= \bE \left[ K(\cF , {\bf m})
e^{-<\cN_{\mu , \lambda} , f>} \right]  .$$
First observe that 
$$\cE_{\lambda}= 
\bE \left[ K(\cF , {\bf m}) 
\exp \left( \sum_{\sigma \in {\rm Lf}(\cF_{\lambda}) } \log \left(  1- \frac{1}{\lambda} 
\int_0^{\mu }du \left( 1- e^{-f(\sigma , u)}\right)\right) \right) \right].$$ 
Note that all the product and sums involved in the latter expression are finite since $f$ has 
compact support. Let $r>0$ be such that $f(\sigma , u)=0$ for all $\sigma $ such that 
$||\sigma||_1 \geq r$ and all $u\in [0, \mu]$. We now use the elementary inequality 
$$ 0 \leq -\log (1-x) -x \leq \frac{x^2}{2(1-x)} \; , \quad x\in [0, 1)  $$
to get 
\begin{eqnarray*}
\left| \sum_{\sigma \in {\rm Lf}(\cF_{\lambda}) } \log \left(  1- \frac{1}{\lambda} 
\int_0^{\mu }du \left( 1- e^{-f(u, \sigma)}\right)\right)
+\int_0^{\mu }du \int \frac{1}{\lambda} { \bf m}_{\lambda} (d\sigma) \left( 1- e^{-f(u, \sigma)}\right)  \right| \\
\leq \frac{1}{2\lambda} \frac{\mu^2}{1-\mu / \lambda}  \frac{1}{\lambda} { \bf m}_{\lambda} 
\left( B_{\ell_1 (\bN) } (0, r)   \right) \; . 
\end{eqnarray*}
The first point of the Theorem then implies that 
$$\lim_{\lambda \rightarrow \infty} \cE_{\lambda}= \bE \left[ K(\cF , {\bf m}) 
\exp \left( - \int_0^{\mu }du \int { \bf m} (d\sigma ) \left( 1- e^{-f(u, \sigma)}\right) \right) \right]  . $$
This implies that the following joint convergence 
\begin{equation}
\label{loiconverg}
\left( \cF , {\bf m} , \cN_{\mu , \lambda} 
\right) \build{\longrightarrow}_{\lambda \to\infty}^{\; } \left( \cF , {\bf m} , \cN_{\mu , \infty} \right)
\end{equation}
holds 
in distribution on $\bT_{\ell_1} \times \cM (\ell_1 (\bN))\times\cM (\ell_1 (\bN)\times[0,\mu])^2$; here $ \cN_{\mu
  , \infty}$ stands for 
a Cox process on 
$\ell_1 (\bN)\times [0, \mu]$ with 
random intensity $ {\bf m}(d\sigma ) \otimes \un_{[0, \mu ]} (x)dx $. Using Skorohod's 
representation theorem, we assume that (\ref{loiconverg}) holds a.s. (for convenience we keep denoting the random 
variables in the same way). For any $\lambda \in [0, \infty)\cup \{
\infty \}$, we denote by $\cP_{\mu , \lambda}$ the set of 
$\sigma \in \ell_1 (\bN)$ for which there exists $U\in [0, \mu ]$ such that $(\sigma , U)$ is an atom of 
  $ \cN_{\mu , \lambda}$. We also introduce the subtree $\cF_{\mu , \lambda}$ of $\cF$ spanned by $0$ and 
the points of $\cP_{\mu , \lambda}$:
$$ \cF_{\mu , \lambda}: = \bigcup_{ \sigma \in \cP_{\mu , \lambda}} \lgeo 0, \sigma \rgeo  \; .$$
Clearly 
\begin{equation}
\label{identpreuve}
 \overline{\cF}_{\mu , \infty} \overset{(d)}{=}\overline{\cF}_{\mu}'  . 
\end{equation}
Next deduce from Lemma \ref{percola} and from the definition of $\cN_{\mu , \lambda}$ that 
the distribution of $\overline{\cF}_{\mu , \lambda} $ does not depend on $\lambda $ and is equal to 
$\Delta_{\mu}^a$. Observe now that for any $r>0$ such that 
$||\sigma||_1 \neq r$ if $\sigma \in \cP_{\mu , \infty}$,  (\ref{loiconverg}) implies 

$$d_{Haus} \left( \cP_{\mu , \lambda} \cap B(0, r) , \cP_{\mu , \infty} \cap B(0, r)  \right)
 \build{\longrightarrow}_{\lambda \to\infty}^{\; } 0 .$$
(Recall that $d_{Haus}$ stands for the Hausdorff distance on the compact
 sets of $\ell_1 (\bN)$). Next, set for any $r>0$ and any $\lambda \in [\mu , \infty) 
\cup \{\infty \}$ 
$$\cF_{\mu , \lambda}(r) =\bigcup \{ \lgeo 0, \sigma \rgeo \; ; \;  \sigma \in \cP_{\mu , \lambda} \cap B(0, r) \}. $$
or any $\sigma , \sigma ' \in \cF$ 
$$ d_{Haus} \left( \lgeo 0, \sigma \rgeo , \lgeo 0, \sigma' \rgeo
\right) \leq || \sigma - \sigma ' ||_1  . $$
Thus, we get  
$$ d_{Haus} \left(\cF_{\mu , \lambda}(r) , \cF_{\mu , \infty}(r)  \right) \leq 
d_{Haus} \left( \cP_{\mu , \lambda} \cap B(0, r) , \cP_{\mu , \infty} \cap B(0, r)  \right) . $$
Then for  any $r>0$ such that 
$||\sigma||_1 \neq r$ if $\sigma \in \cP_{\mu , \infty}$, 
\begin{equation}
\label{hauspoint}
d_{Haus} \left(\cF_{\mu , \lambda}(r) , \cF_{\mu , \infty}(r)  \right)
 \build{\longrightarrow}_{\lambda \to\infty}^{\; } 0 .
\end{equation}
Let $r>0$ be such that 
$||\sigma||_1 \neq r$ if $\sigma \in \cP_{\mu , \infty}$. Since $
\cP_{\mu , \infty}$ has no limit point, we can find $\eta \in (0,1)$ such
that 
\begin{equation}
\label{couronne}
\cP_{\mu , \infty} \cap (B(0, r+\eta ) \backslash B(0,r-\eta)) =\emptyset.
\end{equation}
For the same reason, there is only a finite number of connected
components $C_1, \ldots, C_k$ of $\cF \backslash
B(0,r)$ containing at least one point of $\cP_{\mu , \infty}$. For any
$1\leq i\leq k$, denote by $\sigma_i$ the point of $\cF$ on which
$C_i$ is grafted (observe that $||\sigma_i||_1=r$). Then, 
\begin{equation}
\label{expli1}
\cF_{\mu , \infty} \cap  B(0, r)= \cF_{\mu , \infty}(r)
\bigcup_{i=1}^k \lgeo 0,\sigma_i \rgeo
\end{equation}
Set $R=\max_{1\leq i\leq k} \min \{ ||\sigma ||_1 \; , \; \sigma \in
\cP_{\mu , \infty}\cap C_i  \}$. Observe that $R>r+\eta $ and that for
any $r'>R$, we have 
\begin{eqnarray}
\label{expli2}
\cF_{\mu , \infty} (r')\cap  B(0, r) &=& \cF_{\mu , \infty}(r)
\bigcup_{i=1}^k \lgeo 0,\sigma_i \rgeo \\
&=& \cF_{\mu , \infty} \cap  B(0, r).
\end{eqnarray}
Now for any $\lambda >\mu$ such that 
$$ d_{Haus} \left( \cP_{\mu , \lambda} \cap B(0, R+1) , \cP_{\mu , \infty} \cap B(0, R+1)  \right)
< \eta /2 ,$$
the connected components of  $\cF \backslash
B(0,r)$ containing at least one point of $\cP_{\mu , \lambda}$ are
exactly $C_1, \ldots, C_k$. Thus, 
\begin{eqnarray}
\label{expli3}
\cF_{\mu , \lambda} (R+1)\cap  B(0, r) &=& \cF_{\mu , \lambda}(r)
\bigcup_{i=1}^k \lgeo 0,\sigma_i \rgeo \\
&=& \cF_{\mu , \lambda} \cap  B(0, r) .
\end{eqnarray}
Then by (\ref{expli1}) (\ref{expli2}) and (\ref{expli3}), 
$$ 
d_{Haus} \left(\cF_{\mu , \lambda}\cap  B(0, r) , \cF_{\mu , \infty} \cap  B(0, r) \right)
\leq d_{Haus} \left(\cF_{\mu , \lambda}(R+1) , \cF_{\mu , \infty}
  (R+1) \right) \build{\longrightarrow}_{\lambda \to\infty}^{\; } 0 .$$

This combined with (\ref{hauspoint}) implies that a.s. 
$$ \delta  \left(\overline{\cF}_{\mu , \lambda} , \overline{\cF}_{\mu , \infty}  \right)
 \build{\longrightarrow}_{\lambda \to\infty}^{\; } 0  .$$
Since the distribution of the $\overline{\cF}_{\mu , \lambda}$ is constant and equal to 
$\Delta_{\mu}^a$ , it implies that $\overline{\cF}_{\mu , \infty}$ is also distributed according 
to $\Delta_{\mu}^a$, which proves (\ref{identpreuve}) 
and which completes the proof of (iii) and the proof of the Theorem by 
(\ref{onemarg}).  

\cqfd

\begin{remark}[Connection with previous works in \cite{LGLJ1, LGLJ2, DuLG, DuLG2}] 
\label{concon} 
L\'evy forests have first been defined in the subcritical or critical case via the coding 
by a process $H=(H_t,t\geq 0)$ introduced by Le Gall and Le Jan in \cite{LGLJ1} called the 
$\psi$-height process. 
This process is obtained
from a L\'evy process $X=(X_t,t\geq 0)$ with Laplace exponent $\psi$, by the
following approximation procedure: for every $t\geq 0$, the following
limit in probability exists 
$$H_t=\lim_{\varepsilon\rightarrow 0} \frac{1}{\varepsilon}\int_0^t 
ds\,{\bf 1}_{\{X_s\leq I^s_t+\varepsilon\}} , $$
where we have set $I^s_t :=\inf_{s\leq r\leq t}X_r$ (this
approximation is a consequence of Lemma 1.1.3 in \cite{DuLG}). Set $T_a=\inf \{t\geq 0 \, : \; X_t =-a\}$.
Then, the process $(H_t,0\leq t\leq T_a)$ represents the ``contour'' of the
tree $(\cF,d,\rho)$ in the following sense. For any $s,s'\in [0, T_a]$, set  
$$ d(s,s')=H_s+H_{s'} -2\inf_{s\wedge s' \leq u\leq s\vee s'} H_u $$
and introduce the equivalence relation $s\sim s'$ 
iff $d(s,s')=0$. Then Theorem 2.1 in \cite{DuLG2} asserts that  
$$ \left( \cF \, ,\;  d \, ,\; \rho \right):=\left( [0, T_a]/\sim \, ,\;  d \, ,\;  \widetilde{0} \right)$$
is a compact random real rooted tree; for any $s\in [0, T_a] $, denote
by $\widetilde{s}$ the $\sim $-isometry class of $s$. 
Let us explain why $\overline{\cF}$, defined in this way, is an 
$(a, \psi)$-L\'evy forest.  
Let $\cP=\{ (t_i,r_i) , i\in J\}$
be a Poisson point process on $[0, \infty)^2$ with intensity the Lebesgue measure.  For any $\lambda \geq 0$ 
we set 
$$ \cF(\lambda )=\bigcup \{ \lgeo \rho , \widetilde{t_i} \rgeo \; ; \; i\in J \; :\; r_i \leq \lambda  \; ; \; t_i \leq T_a \} . $$ 
Obviously the family of real trees $(\cF (\lambda); \lambda \geq 0)$ is consistent under Bernoulli leaf colouring
and Theorem 3.2.1 \cite{DuLG} asserts that 
$(\cF (\lambda),d,\rho)$ is a GW$(\xi_{\lambda}, \psi'(\psi^{-1}(\lambda)), a\psi^{-1}(\lambda))$ -real rooted forest 
Thus, $(\cF (\lambda); \lambda \geq 0)$ is an 
$(a, \psi)$-growth process. Besides, it is clear from the construction that a.s. 
$$ \lim_{\lambda \rightarrow \infty} \delta (\cF (\lambda), \cF ) =0 . $$

   Moreover, if we take 
$(\cF (\lambda); \lambda \geq 0)$ in Theorem \ref{meslebe}, the mass 
distribution is clearly the image of the Lebesgue measure on the line 
by the canonical projection associated with $\sim$. 
We refer to \cite{DuLG2} for discussion of various geometric properties of L\'evy forests. 
\end{remark}

\begin{remark}
\label{echec} The construction of the mass measure on a L\'evy tree given in \cite{DuLG2} 
only relies on the metric structure of the L\'evy tree and not on a particular coding (see the remark before Theorem 4.4 in 
\cite{DuLG}). 
We failed to give a proof that is well-suited to our approach that $ { \bf m}$ is actually a deterministic functional
of its topological support $\cF$. 
\end{remark}

\subsection{Excursion measure of L\'evy trees.}
\label{excursionmes}

  Fix $a>0$  and consider an $(a, \psi)$-L\'evy forest $\cF$. Denote 
by $\cT_i^o$, $i \in J$, the connected components of 
$ \cF \backslash \{ \rho\}$ and set for any $i\in J$, $\cT_i  =\{
\rho\} \cup \cT_i^o$. 
The main goal of this section
is to define a Borel measure $\Theta (d\overline{\cT})$ on $\bT$ such
that the 
following proposition holds:
\begin{proposition} 
\label{excuracine}
The point measure 
$$ \cN (d\overline{T}):= \sum_{i\in J} \delta_{\overline{\cT}_i}  (d\overline{T}) $$
is a Poisson point measure on $\bT$ with intensity $a\, \Theta (d\overline{\cT})$. 
\end{proposition}

   Before proving this proposition, recall the notation $P_{\mu} (d\overline{\cT})$ and  
$P_{\mu ,\lambda } (d\overline{\cT})$ from Subsection \ref{levyforest}. We first establish 

\begin{equation}
\label{posipsi}
\mbox{\rm\bf \hspace{-0.9cm}Claim:} \qquad \qquad \qquad \qquad 
P_{\lambda }= \left( 1-\frac{\psi^{-1}(\mu)}{\psi^{-1}(\lambda)}\right) P_{\mu ,\lambda } + 
\frac{\psi^{-1}(\mu)}{\psi^{-1}(\lambda)} P_{\mu }.\qquad \qquad \qquad \qquad 
\end{equation}
{\bf Proof of the claim:} Let $\cT_{\mu}$ and $\cT_{\mu, \lambda}$ be as in the 
last point of Notation \ref{levynot}. Perform a $(1-\mu/\lambda )$-Bernoulli leaf colouring on $\cT_{\lambda}$. Recall that 
the probability that $\cT_{\lambda}$ is completely red  is $1-\psi^{-1}(\mu)/\psi^{-1}(\lambda)$. Moreover, 
conditionally on this event,  $\overline{\cT}_{\lambda}$ is distributed as 
$\overline{\cT}_{\mu, \lambda  }$ and conditionally on the
complementary event, $\overline{\cT}_{\lambda}$ is distributed as $\overline{Q}_{\mu, \lambda  } (\cT_{\mu}) $. Then, flip a coin with 
probability  $1-\psi^{-1}(\mu)/\psi^{-1}(\lambda)$ to be head. If it is head, then set $\cT'= \cT_{\mu, \lambda}$; otherwise set 
$\cT'= Q_{\mu, \lambda  } (\cT_{\mu}) $. The previous observations imply that $\overline{\cT}'$ and 
$\overline{\cT}_{\lambda}$ have the same distribution. Accordingly, 
$\overline{Q}_{\lambda , \infty } (\cT_{\lambda}) $ and
$\overline{Q}_{\lambda , \infty } (\cT') $ have the same
distribution. Use now (\ref{chaininginf}) with $T=\cT_{\mu}$ to get 
$$ \overline{Q}^a_{\lambda, \infty} \left( Q^a_{\mu , \lambda} (\cT_{\mu})\right) \quad 
\overset{(d)}{=} \quad \overline{Q}^a_{\mu, \infty}(\cT_{\mu}) \;  . $$
This, combined with the previous observation imply the claim. \cqfd

\vspace{3mm}

   Let $\lambda_0=0 < \lambda_1 <  \lambda_2 , \ldots$ be any increasing sequence going to infinity. We define the 
excursion measure by 
$$  \Theta (d\overline{\cT})= \gamma P_{0 } (d\overline{\cT}) +\sum_{n\geq 0} 
(\psi^{-1}(\lambda_{n+1})  -\psi^{-1}(\lambda_n)) P_{\lambda_n,\lambda_{n+1} } (d\overline{\cT}). $$
Recall that $\gamma=\psi^{-1}(\lambda_0)$.
Let us first prove that $  \Theta (d\overline{\cT})$ 
does not depend on $(\lambda_n ; n\geq 0)$ and more precisely for any 
non-negative measurable function $K$ on $\bT$, let us prove that 
\begin{equation}
\label{limtheta}
<\Theta , K> =\lim_{\lambda \rightarrow \infty} \uparrow <\psi^{-1}(\lambda) P_{\lambda} , K> . 
\end{equation}
{\bf Proof of (\ref{limtheta}):} (\ref{posipsi}) implies that $\lambda \rightarrow  <\psi^{-1}(\lambda) P_{\lambda} , K>$
is non-decreasing. Thus, the limit in (\ref{limtheta}) exists in $[0, \infty]$. Denote this limit by $L(K)$ and observe that
$$ \psi^{-1}(\lambda_n) <P_{\lambda_n} , K> =  \gamma P_{0 } (d\overline{\cT}) +\sum_{k= 0}^{n-1}  
(\psi^{-1}(\lambda_{k+1})  -\psi^{-1}(\lambda_k)) <P_{\lambda_k,\lambda_{k+1} } , K >.$$
Thus, by letting $n$ go to infinity, we get $<\Theta , K>=L(K)$, which proves (\ref{limtheta}). \cqfd

\vspace{3mm}
\noindent 
{\bf Proof of the proposition:} Fix $\lambda >0$ and define 
$$ J_{\lambda }= \{i\in J \; :\; \cT_i \cap \cF_{\lambda} \neq \{ \rho \} \} . $$
Set $\cT_i (\lambda)= \cT_i \cap \cF_{\lambda} $ for any $i \in J_{\lambda}$. From the 
construction of the growth process, we deduce that 
$\# J_{\lambda}$ is a Poisson random variable with parameter $a\psi^{-1}(\lambda)$ and that conditionally on $J_{\lambda}$, 
the $\overline{\cT}_i  (\lambda)$ , $i\in J_{\lambda}$, are i.i.d. and distributed as the isometry class of a 
GW($\xi_{\lambda}, \psi' \psi^{-1}(\lambda)$)-real tree.  Now observe that for any $i\in J_{\lambda}$ 
the tree $\cT_i$ is obtained as the limit of a growth process 
started at $\cT_i (\lambda)$. Then $\overline{\cT}_i  $ and $\overline{Q}_{\lambda , \infty} (\cT_i (\lambda))$  have the 
same distribution. Thus conditionally on $J_{\lambda}$ the $\overline{\cT}_i$ ,  
$i\in J_{\lambda}$, are independent and distributed according to
$P_{\lambda}$ and for any non-negative 
measurable function $K$ on $\bT$, we have 
\begin{eqnarray*}
\bE \left[ \exp \left( -\sum_{i\in J_{\lambda}}  K(\overline{\cT}_i)   \right)\right] &=& 
\bE \left[  \left( \int  P_{\lambda} (d\overline{\cT}) e^{-K(\overline{\cT})  } \right)^{\# J_{\lambda}}  \right] \\
&=& \bE \left[ \exp \left( - a\psi^{-1}(\lambda) \int  P_{\lambda} (d\overline{\cT}) \left( 1-e^{-K(\overline{\cT})  } \right)
\right)\right] \; .
\end{eqnarray*}
Now, observe that 
$$<\cN, K>=\lim_{\lambda \rightarrow \infty} \uparrow  \sum_{i\in J_{\lambda}}  K(\overline{\cT}_i) \, ,  $$
which completes the proof by (\ref{limtheta}) and by the dominated convergence theorem. \cqfd

\vspace{3mm}

   Recall that the height $h(T)$ of a rooted real tree $(T, d, \rho)$ is the (possibly infinite)
real number  $ \sup \{ d(\rho , \sigma ) \; , \; \sigma \in T\}$. Observe that $h(T)$ is invariant up to isometry so it makes sense 
to define $h(\overline{T})$ as the height of any representative of
$\overline{T}$. 
It is easy to check here that a.s.
\begin{equation}
\label{limicroisheight}
 \lim_{\lambda \rightarrow \infty} \uparrow h(\overline{\cF}_{\lambda} )=  h(\overline{\cF} ) \; .
\end{equation}
Recall from (\ref{heightharris}) that the probability that the height of a single 
GW($\xi_{\lambda}, \psi'(\psi^{-1}(\lambda))$)-real tree is greater that $x$ is $e^{-v_{\lambda} (x)}$ where 
$v_{\lambda} (x)$ satifies 
$$ \psi^{-1} (\lambda)  \int_0^{e^{-v_{\lambda} (x)}} \frac{ du}{\psi ((1-u)\psi^{-1} (\lambda))}=x . $$
Then, (\ref{limicroisheight}) and a simple computation imply that 
$$ \bP \left( h(\overline{\cF} )\le x \right)= \exp (-av(x)) \; , $$
where $v$ satisfies the equation 
\begin{equation}
\label{equav}
\int_{v(x)}^{\infty}  \frac{ du}{\psi (u)} =x \; .
\end{equation}
Now observe that $ h(\overline{\cF})=\sup \{ h(\overline{\cT}_i ) \; , \; i\in J\}$. Thus Proposition 
\ref{excuracine} implies that 
\begin{equation}
\label{thetav}
\Theta \left(   h(\overline{\cT}) >x  \right) =v(x). 
\end{equation}

\begin{remark}
\label{approracine}
Observe that Proposition \ref{excuracine} and (\ref{thetav}) imply that a.s. 
$$a= \lim_{\epsilon \rightarrow 0} \frac{1}{v(\epsilon)} \# \{ i\in J \; : \;  h (\overline{\cT}_i ) \} . $$
\end{remark}

\begin{notation}
\label{thetanota} 
The measure $\Theta $ is called the $\psi$-{\it excursion measure}. The terminology comes from the fact that in the critical or
subcritical case when the L\'evy forest is coded by a $\psi$-height process as explained in Remark \ref{concon}, 
$\Theta$ is the distribution of the tree coded by one excursion above $0$ of the height process. In the last section we shall use 
the notation $\Theta_{\lambda}$ for the $\psi_{\lambda}$-excursion measure. 
\end{notation}

\subsection{Decomposition of the L\'evy forest along the ancestral tree of a Poisson sample}
\label{decpoissonsample}

  Fix $\mu_0 \geq 0$ and $a\geq 0$. Consider a $\btl$-valued $(a,
  \psi)$-growth 
process $(\cF_{\lambda}; \lambda \geq 0)$ 
denote by $\cF$ the 
limit of this growth process. Recall that 
$\overline{\cF}_{\mu_0}$ is distributed as the isometry class of 
the ancestral subtree of a Poisson sampling on $\cF$ with intensity
  $\mu_0  \, . \,  {\bf m}$.  The aim of this subsection is to compute
  the distribution of 
$\cF$ conditionally on $\cF_{\mu_0}$, as the reconstruction 
procedure does in the discrete case. To avoid technicalities and to make
  easier 
the statement of this decomposition we also assume that 
\begin{equation}
\label{technicassu}
 \cF_{\mu_0}= Q^a_{0, \mu_0} (\cF_0)  , 
\end{equation}
where the extra random variables used to define $ Q^a_{0, \mu_0}$ are
chosen independent of $\cF_0$. Before stating the main result, we need to
introduce some notation: Denote by ${\rm Gr}$ the set of points on which the
connected components of $\cF \backslash \cF_{\mu_0}$ are grafted; for any 
$\sigma \in {\rm Gr}$, 
denote by  $F_i^o (\sigma )$, $i\in J(\sigma)$, the  connected components of 
$\cF \backslash \cF_{\mu_0}$ that are grafted on $\sigma$ and set 
$$F_{\sigma}=\{ \sigma \}\cup \{ F_i^o (\sigma ) , i\in J(\sigma) \}
.$$ 
Observe that 
$F_{\sigma}$ is a closed and connected set. Next, 
let us introduce the sets of points 
$$ S^1_{\mu_0} =\{ \sigma \in {\rm Gr}\backslash ({\rm Br}(\cF_{\mu_0}) \cup \{ \rho \})  \; : \;    \# J(\sigma ) \geq 2  \} $$
and 
$$ S^2_{\mu_0} =\{ \sigma \in {\rm Gr}\backslash ({\rm Br}(\cF_{\mu_0}) \cup \{ \rho \})  \; : \;    \# J(\sigma ) =1  \} .$$
Note that some branching points of $\cF_{\mu_0}$ may not be in ${\rm Gr}$. We then set 
$$  S_{\mu_0} =S^1_{\mu_0} \cup S^2_{\mu_0} \cup {\rm Br}(\cF_{\mu_0})\cup \{ \rho \} , $$
and if $\sigma \in  ({\rm Br}(\cF_{\mu_0}) \cup \{ \rho \}) \backslash {\rm Gr} $, then  we set $F_{\sigma }= \{ \sigma\}$.

   Recall from the end of Section \ref{excursionmes} the notation $\Theta_{\mu_0}$ for the $\psi_{\mu_0}$-excursion 
measure and also recall from Section \ref{levyforest} the notation $P^r_{\mu_0 }$. Let us denote by $v_{\mu_0}$ the 
function defined on $[0, \infty)$ that satisfies 
\begin{equation}
\label{intequamu_0}
\int_{v_{\mu_0} (t)}^{\infty} \frac{du}{\psi_{\mu_0} (u)} =t. 
\end{equation}

\begin{theorem}
\label{poidecanc}
Almost surely for every $\sigma \in S_{\mu_0 }$, the following limit exists and is finite 
$$ a(\sigma):= \lim_{\epsilon \rightarrow 0} \frac{1}{v_{\mu_0} (\epsilon)} \# \{ i\in J(\sigma) \; : \;  h(F_i^o (\sigma )) > 
\epsilon \} . $$
Moreover, conditionally on $\cF_{\mu_0}$ the collections of random variables 
$$ \cP_1 = \{ (\sigma ,a(\sigma)  , \overline{F}_{\sigma}) \; , \; \sigma \in S^1_{\mu_0} \} \quad {\rm } \quad 
\cP_2 = \{ (\sigma , \overline{F}_{\sigma}) \; , \; \sigma \in S^2_{\mu_0} \}$$
and $\cP_3=  \{ (a(\sigma) , \overline{F}_{\sigma})$, $\sigma \in   {\rm Br}(\cF_{\mu_0}) \cup \{ \rho \} \}$, are 
independent. Their conditional distributions are given by the following: 

\begin{enumerate}

\item[{\rm(i)}] $\cP_1 $ is a Poisson point process on $\cF_{\mu_0} \times [0, \infty) \times \bT$ with intensity measure
$$ \ell_{\cF_{\mu_0}} (d\sigma) \otimes e^{-r \psi^{-1}(\mu_0)}
 \Pi (dr) \otimes P^r_{\mu_0} (d\overline{\cT}) ;$$

\item[{\rm(ii)}] $\cP_2 $ is a Poisson point process on $\cF_{\mu_0} \times \bT$ with intensity measure
$$ 2 \beta \ell_{\cF_{\mu_0}} (d\sigma) \otimes \Theta_{\mu_0} (d\overline{\cT}) ;$$

\item[{\rm(iii)}] The $ (a(\sigma) , \overline{F}_{\sigma})$, $\sigma \in   {\rm Br}(\cF_{\mu_0}) \cup \{ \rho \}$ are independent 
random variables; Moreover, for each  $\sigma \in   {\rm
  Br}(\cF_{\mu_0}) $ , the $ [0, \infty) \times \bT$-valued random
variables $ (a(\sigma) , \overline{F}_{\sigma})$ are  distributed according to 
$$ \eta_{\mu_0 , l} (dr) \otimes P^r_{\mu_0}  (d\overline{\cT})$$
where $l={\rm n}(\sigma , \cF_{\mu_0})-1$; $a(\rho)=a$ and 
$\overline{F_{\rho}}$ is distributed according to $P_{\mu_0}^a$. 
\end{enumerate}
\end{theorem}

\begin{remark}
\label{casdeux} 
Recall that $P^0_{\mu_0} =\delta_{\overline{\{ \rho\}}}$. If $l=2$ in (iii), then since  
$$ \eta_{\mu_0 , 2} (\{ 0\})=\frac{2\beta}{\psi^{(2)} (\psi^{-1} (\mu_0))} , $$
$F_{\sigma}$ reduces to a point with probability $\eta_{\mu_0 , 2} (\{ 0\}) >0$ as soon as $\beta >0$. 
\end{remark}
{\bf Proof:} Set for any $\lambda \geq \mu_0 $ and any $\sigma \in S_{\mu_0}$,  $F_{\sigma} (\lambda)= F_{\sigma}\cap 
\cF_{\lambda}$ and 
$$  S_{\mu_0 , \lambda} =\{ \sigma \in S_{\mu_0} \; : \; F_{\sigma} (\lambda) \neq \{ \sigma\} \} .$$
Deduce from Theorem \ref{mainresult} that a.s. for any 
$\sigma \in S_{\mu_0}$ 
\begin{equation}
\label{loclim}
\dl (F_{\sigma} (\lambda), F_{\sigma})
\; \build{\longrightarrow}_{\lambda \to\infty}^{\; } \;0 \; .
\end{equation}
Recall notation $\Delta_{\mu_0 , \lambda}$ and $\Delta^r_{\mu_0 , \lambda}$ from the end of Section \ref{levygal}. For convenience of notation, let us set 
$$ M_1 (d\overline{\cT})= \int_{(0, \infty)} \Pi (dr) re^{-r\psi^{-1}(\mu_0)}P^r_{\mu_0} (d\overline{\cT}) $$
and 
$$  M_2 (d\overline{\cT})= 2\beta \Theta_{\mu_0 } (d\overline{\cT}) .$$
\begin{lemma}
\label{techdelt}
For any non-negative continuous function $R$ on $\bT$, any $\mu_0 \geq
0$ 
and any $\epsilon >0$, we have 
\begin{enumerate}
\item[{\rm(a)}] $\displaystyle 2\beta\left( \psi^{-1}(\lambda) -\psi^{-1}(\mu_0) \right) 
\int_{\{ h(\overline{\cT}) > \epsilon \}} 
\Delta_{\mu_0 , \lambda}( d\overline{\cT}) \left( 1- e^{-R (\overline{\cT}) }\right)$

$$\build{\longrightarrow}_{\lambda \to\infty}^{\; }
\int_{\{ h(\overline{\cT}) > \epsilon \}}  M_2 (d\overline{\cT})\left( 1- e^{-R (\overline{\cT}) }\right).
$$

\item[{\rm(b)}] $\displaystyle\int_{(0, \infty)} \Pi (dr) re^{-r\psi^{-1}(\mu_0)} \int_{\{ 
h(\overline{\cT}) > \epsilon;\overline{\cT} \neq \overline{\rho} 
\}} 
\Delta^r_{\mu_0 , \lambda}( d\overline{\cT}) \left( 1- e^{-R (\overline{\cT}) }\right)$

 $$\build{\longrightarrow}_{\lambda \to\infty}^{\; } 
\int_{\{ h(\overline{\cT}) > \epsilon \}}  M_1 (d\overline{\cT})\left( 1- e^{-R (\overline{\cT}) }\right).$$
\end{enumerate}
\end{lemma}

\vspace{2mm}

\noindent
{\bf End of the proof of the theorem:} Before proving the lemma, let us complete the proof of the theorem. Let $K$ be a non-negative continuous function on
$\ell_1 (\bN) \times \bT$ such that $K(\sigma , \overline{\cT})=0$ for every $\overline{\cT} \in \bT$ and every 
$\sigma$ such that $|| \sigma ||_1 \geq r_0$, where $r_0$ is a fixed positive
number. First deduce from (\ref{loclim}) that a.s. 
for every $\sigma \in S_{\mu_0}$ 
$$ \lim_{\lambda \to\infty} \uparrow h( \overline{F}_{\sigma} (\lambda)) = h(\overline{F}_{\sigma}  ) .$$
Fix $\epsilon >0$. Since $\cF$ is locally compact, there is only a finite number of $F_{\sigma}$'s such that 
$h(\overline{F}_{\sigma}  ) >\epsilon $ and $|| \sigma ||_1 \leq r_0$. Thus for any $i\in \{ 1,2\}$ a.s. 
$$  
\lim_{\lambda \to\infty} \sum_{\sigma \in S_{\mu_0 , \lambda } \cap S^i_{\mu_0}} 
K (\sigma , \overline{F}_{\sigma}  (\lambda)) \un_{ \{ h(\overline{F}_{\sigma} (\lambda) ) >\epsilon \}}= 
\sum_{\sigma \in S^i_{\mu_0}} 
K (\sigma , \overline{F}_{\sigma}) \un_{\{ h(\overline{F}_{\sigma}  ) >\epsilon \}} . $$

Since we have supposed (\ref{technicassu}),  
$$ \cP_1 (\lambda)=\{ (\sigma ,\overline{F}_{\sigma }(\lambda) ) , \; \sigma \in 
S^1_{\mu_0} \cap  S_{\mu_0, \lambda }
 \}\; , \; 
\cP_2 (\lambda)=\{ (\sigma ,\overline{F}_{\sigma }(\lambda) ) , \; \sigma \in 
S^2_{\mu_0} \cap  S_{\mu_0, \lambda } \}$$
and 
$\cP_3 (\lambda)=\{ (\sigma ,\overline{F}_{\sigma }(\lambda) ) \; , \; \sigma \in 
 {\rm Br}(T)\cup \{ \rho\}\} $ are distributed as specified in Remark
 (\ref{loifsigmulam}) with $T=\cF_{\mu_0}$. Then 
deduce from Remark \ref{loifsigmulam} (i) and (ii) and from Lemma \ref{techdelt} that 
\begin{eqnarray*}
&&\hspace{-0.5cm}\bE \left[\left. \exp \left( - \sum_{\sigma \in S^i_{\mu_0}} 
K (\sigma , \overline{F}_{\sigma}) \un_{\{ h(\overline{F}_{\sigma}  ) >\epsilon \}} \right) \right| \cF_{\mu_0} , 
S^i_{\mu_0} \cap  S_{\mu_0, \lambda }
\right] \\
&&=\exp \left( - \int \ell_{\cF_{\mu_0 } } (d\sigma )
\int_{\{ h(\overline{\cT}) > \epsilon \}}  M_i (d\overline{\cT})\left( 1- e^{-K (\sigma , \overline{\cT}) }\right)\right),  
\end{eqnarray*}
for $i\in \{ 1,2\}$. Now, let $\epsilon $ go to $0$: It implies that conditionally on 
$\cF_{\mu_0}$ the sets of points $\{ (\sigma , \overline{F_{\sigma}} ) \, ; \, \sigma \in S^i_{\mu_0}\}$, $i\in \{ 1,2\}$ 
are two independent Poisson point processes with resp. intensities 
$ \ell_{\cF_{\mu_0 } }\otimes  M_i$ , $i\in \{ 1,2\}$ .

Recall that (\ref{weakdelta}) asserts that for any $r>0$, the probability measure $\Delta^r_{\mu_0 , \lambda}$ on $\bT$ weakly 
converges to $P^r_{\mu_0}$. This observation combined with Remark \ref{loifsigmulam} imply that conditionally on 
$\cF_{\mu_0}$ for every $\sigma \in {\rm Br} (\cF_{\mu_0}) $, $\overline{F}_{\sigma}$ is distributed according to 
$$ \int \eta_{\mu_0 , l} (dr)  P^r_{\mu_0}  (d\overline{\cT}) , $$
(with $l= {\rm n} (\sigma , \cF_{\mu_0}) -1)$, and  that $\overline{F}_{\rho}$ is distributed according to  $P^a_{\mu_0}$. 
Then, Remark \ref{approracine} implies the first point of the theorem; this, combined with the previous 
observations, implies that 
conditionally on $\cF_{\mu_0}$, $\cP_1$ , $\cP_2$ and $\cP_3$ are distributed as specified in the theorem; then, their 
conditional independence 
is an easy consequence of the conditional independence of 
 $\cP_1(\lambda) $ , $\cP_2(\lambda) $ and $\cP_3(\lambda)$ stated in Remark \ref{loifsigmulam}.  This completes the proof of 
the theorem. \cqfd

\vspace{3mm}

\noindent
{\bf Proof of Lemma \ref{techdelt}:} Recall (\ref{deltashift}) that makes the connection between the distribution 
$\Delta_{\lambda}$ and $\Delta_{\mu_0 , \lambda}$. By replacing $\psi$ by $\psi_{\mu_0}$, the first point of the 
lemma is then equivalent to 
the following limit
\begin{equation}
\label{actuallim}
\psi^{-1}(\lambda)
\int_{\{ h(\overline{\cT}) > \epsilon \}} 
\Delta_{ \lambda}( d\overline{\cT}) \left( 1- e^{-R (\overline{\cT}) }\right)
\build{\longrightarrow}_{\lambda \to\infty}^{\; } 
\int_{\{ h(\overline{\cT}) > \epsilon \}}  \Theta (d\overline{\cT})\left( 1- e^{-R (\overline{\cT}) }\right) .
\end{equation}

Recall from Section \ref{excursionmes} the notation
$\cT_i$ , $i\in J$, for the subtrees of $\cF$ grafted at $\{ \rho\}$. Set 
$\cT_i (\lambda)=\cF_{\lambda} \cap \cT_i $, $i\in J$ and 
$$J(\lambda)=\{ i\in J \; :\; \cT_i (\lambda) \neq \{ \rho\}\} . $$
Since $\# J(\lambda)$ is a Poisson variable with parameter $a \psi^{-1}(\lambda)$ and since conditionally on 
$ J(\lambda)$, the trees $\overline{\cT}_i (\lambda)$, $i\in J(\lambda)$ are independent with the same distribution 
$\Delta_{\lambda}$, we get
\begin{eqnarray}
\label{expression}
&&\hspace{-0.5cm}\bE \left[ \exp \left( - \sum_{i\in J(\lambda)} R(\overline{\cT}_i (\lambda)) \un_{\{  h(\overline{\cT}_i (\lambda)) >
\epsilon \}}  \right) \right]\\
&&=\exp \left( -a \psi^{-1}(\lambda)  \int_{\{ h(\overline{\cT}) > \epsilon \}} 
\Delta_{ \lambda}( d\overline{\cT}) \left( 1- e^{-R (\overline{\cT}) }\right)\right) \; .
\end{eqnarray}
Now observe that a.s. for any $i\in J$ 
\begin{equation}
\label{pshaut}
\lim_{\lambda \to\infty} \delta (\overline{\cT}_i , \overline{\cT}_i (\lambda))=0 \quad {\rm and } \quad 
\lim_{\lambda \to\infty} \uparrow h(\overline{\cT}_i (\lambda))= h(\overline{\cT}_i ) \; . 
\end{equation}
Since $\cF$ is locally compact, there are only finitely many $\cT_i$'s such that $h(\overline{\cT}_i)>\epsilon$. Thus 
(\ref{pshaut}) implies 
\begin{equation}
\label{limisomme}
\lim_{\lambda \to\infty}
\sum_{i\in J(\lambda)} R(\overline{\cT}_i (\lambda)) \un_{\{  h(\overline{\cT}_i (\lambda)) >
\epsilon \}}  =
\sum_{i\in J} R(\overline{\cT}_i ) \un_{\{  h(\overline{\cT}_i ) >
\epsilon \}}, 
\end{equation}
and (a) follows from Proposition \ref{excuracine} and (\ref{expression}).

\vspace{3mm}

It remains to prove (b): An elementary computation based on (\ref{heightharris}) with $\varphi = \varphi_{\mu_0, \lambda}$ 
implies that 
$$ \Delta^r_{\mu_0 , \lambda} \left(  h(\overline{\cT}) > \epsilon \right) =1-
\exp \left( -r(\psi^{-1} (\lambda) - 
\psi^{-1} (\mu_0) )(1-e^{-v_{\mu_0 , \lambda} (\epsilon) }) \right)  $$
where $v_{\mu_0 , \lambda}$  satisfies the following equation
$$\int_{  (\psi^{-1} (\lambda) - 
\psi^{-1} (\mu_0) )(1-e^{-v_{\mu_0 , \lambda} (\epsilon) })  }^{\psi^{-1} (\lambda) - 
\psi^{-1} (\mu_0)} \frac{du}{\psi_{\mu_0}(u)} =\epsilon \; .$$
Thus, 
\begin{equation}
\label{heightmu_0}
\lim_{\lambda \to\infty} \uparrow \Delta^r_{\mu_0 , \lambda} \left(  h(\overline{\cT}) > \epsilon \right) = 
1-e^{-rv_{\mu_0 } (\epsilon) } =  P^r_{\mu_0} \left(  h(\overline{\cT}) > \epsilon \right) \; , 
\end{equation}
where $v_{\mu_0 }$ satisfies (\ref{intequamu_0}). By (\ref{weakdelta}), 
for any $r>0$ we get 
\begin{equation}
\label{partweak}
\lim_{\lambda \to\infty}
\int_{\{ h(\overline{\cT}) > \epsilon \; ; \; \overline{\cT} \neq 
\overline{\{\rho\}} \}}
 \Delta^r_{\mu_0 , \lambda}( d\overline{\cT}) 
\left( 1- e^{-R (\overline{\cT}) }\right) 
=\int_{\{ h(\overline{\cT}) > \epsilon \; ; \; \overline{\cT} 
\neq \overline{\rho} \}}
 P^r_{\mu_0 } ( d\overline{\cT} ) 
\left( 1- e^{-R (\overline{\cT}) }\right).
\end{equation}
Now by (\ref{heightmu_0}) 
$$ \int_{\{ h(\overline{\cT}) > \epsilon \; ; \; \overline{\cT} \neq \overline{\rho} \}}
 \Delta^r_{\mu_0 , \lambda}( d\overline{\cT}) \left( 1- e^{-R (\overline{\cT}) }\right) 
\leq  P^r_{\mu_0} \left(  h(\overline{\cT}) > \epsilon \right) =1-e^{-rv_{\mu_0 } (\epsilon) } . $$
Now note that 
$$ \int_{(0, \infty)} \Pi (dr) re^{-r\psi^{-1}(\mu_0)} (1-e^{-rv_{\mu_0 } (\epsilon) }) \; < \; \infty \; , $$
which implies (b) by (\ref{partweak}) and the dominated convergence theorem. This completes the proof of the lemma. 
\cqfd

\bibliographystyle{plain}

\begin{thebibliography}{10}

\bibitem{Ab92}
R.~Abraham.
\newblock Un arbre al{\'e}atoire infini associ{\'e} {\`a} l'excursion
  {B}rownienne.
\newblock In {\em S{\'e}m. de Proba.}, volume XXVI, pages 374--397. Springer,
  Berlin, 1992.

\bibitem{AlPit98}
D.~Aldous and J.~Pitman.
\newblock Tree-valued {M}arkov chains derived from {G}alton-{W}atson processes.
\newblock {\em Ann. Inst. H. Poincar\'e.}, 34:637--686, 1998.

\bibitem{Al1}
D.~J. Aldous.
\newblock The continuum random tree {I}.
\newblock {\em Ann. Probab.}, 19:1--28, 1991.

\bibitem{AthNey}
K.~Athreya and P.~Ney.
\newblock {\em Branching process}.
\newblock Number 196 in Grundlehren der Mathematischen Wissenschaften.
  Springer, 1972.

\bibitem{Bi2}
N.~H. Bingham.
\newblock Continuous branching processes and spectral positivity.
\newblock {\em Stochastic Process. Appl.}, 4:217--242, 1976.

\bibitem{BuBu}
Y.~Burago, D.~Burago and S.~Ivanov.
\newblock {\em A Course in Metric Geometry}, volume~33.
\newblock AMS, Boston, 2001.

\bibitem{Ho}
Hobson D.G.
\newblock Marked excursions and random trees.
\newblock In Springer, editor, {\em Lecture Notes Math.}, volume 1729 of {\em
  S\'eminaire de Probabilit\'es XXXIV}, pages 289--301, 2000.

\bibitem{Dress84}
A.~Dress.
\newblock Trees, tight extensions of metric spaces, and the cohomological
  dimension of certain groups: A note on combinatorial properties of metric
  spaces.
\newblock {\em Adv. Math.}, 53:321--402, 1984.

\bibitem{DMT96}
A.~Dress, V.~Moulton, and W.~Terhalle.
\newblock {T}-theory: an overview.
\newblock {\em European J. Combin.}, 17:161--175, 1996.

\bibitem{DT96}
A.~Dress and W.~Terhalle.
\newblock The real tree.
\newblock {\em Adv. Math.}, 120:283--301, 1996.

\bibitem{DuLG}
T.~Duquesne and J-F. Le~Gall.
\newblock {\em Random Trees, {L}\'evy Processes and Spatial Branching
  Processes}.
\newblock Ast\'erisque no 281, 2002.

\bibitem{DuLG2}
T.~Duquesne and J-F. Le~Gall.
\newblock Probabilistic and fractal aspects of {L}\'evy trees.
\newblock {\em To appear in Probab. Theorey and Rel. Fields}, 2004.

\bibitem{Ev00}
S.~Evans.
\newblock Snakes and spiders: {B}rownian motion on real trees.
\newblock {\em Probab. Theory Related Fields}, 117(3):361--386, 2000.

\bibitem{EvPitWin}
S.N. Evans, J.~Pitman, and A.~Winter.
\newblock Rayleigh processes, real trees, and root growth with re-grafting.
\newblock {\em To appear in Probab. Th. Rel. Fields}, 2005.

\bibitem{EvWin}
S.N. Evans and A.~Winter.
\newblock Subtree prune and re-graft: a reversible real tree valued {M}arkov
  process.
\newblock {\em preprint}, 2005.

\bibitem{Fe}
W.~Feller.
\newblock {\em An Introduction to Probability Theory and Its Applications, Vol.
  II, sec. ed.}
\newblock Wiley, New York., 1971.

\bibitem{GeiKau}
J.~Geiger and Kauffmann L.
\newblock The shape of large {G}alton-{W}atson trees with possibly infinite
  variance.
\newblock {\em Rand. Struct. Alg.}, 25(3):311--335, 2004.

\bibitem{Gro}
M.~Gromov.
\newblock {\em Metric Structures for {R}iemannian and non-{R}iemannian Spaces.}
\newblock Progress in Mathematics. Birkhäuser, 1999.

\bibitem{LyoHam}
B.M. Hambly and T.J. Lyons.
\newblock Uniqueness for the signature of a path of bounded variation and
  continuous analogues for the free group.
\newblock {\em Preprint}, 2004.

\bibitem{Ji}
M.~Jirina.
\newblock Stochastic branching processes with continous state-space.
\newblock {\em Czech. Math. J.}, 8:292--313, 1958.

\bibitem{La2}
J.~Lamperti.
\newblock Continuous-state branching processes.
\newblock {\em Bull. Amer. Math. Soc.}, 73:382--386, 1967.

\bibitem{La1}
J.~Lamperti.
\newblock The limit of a sequence of branching processes.
\newblock {\em Z. Wahrsch. Verw. Gebiete}, 7:271--288, 1967.

\bibitem{La3}
J.~Lamperti.
\newblock Limiting distributions of branching processes.
\newblock In {\em Fifth Berkeley Symp.}, volume II, Part 2, pages 225--241,
  1967.

\bibitem{LGLJ1}
J-F. Le~Gall and Y.~Le~Jan.
\newblock Branching processes in {L}\'evy processes: the exploration process.
\newblock {\em Ann. Probab.}, 26-1:213--252, 1998.

\bibitem{LGLJ2}
J-F. Le~Gall and Y.~Le~Jan.
\newblock Branching processes in {L}\'evy processes: {L}aplace functionals of
  snakes and superprocesses.
\newblock {\em Ann. Probab.}, 26:1407--1432, 1999.

\bibitem{LJ91}
Y.~Le~Jan.
\newblock Superprocesses and projective limits of branching {M}arkov processes.
\newblock {\em Ann. Inst. H.Poincar{\'e }}, 27:91--106, 1991.

\bibitem{LyPe}
R.~Lyons and Y~Peres.
\newblock {\em Probability on Trees and Networks}.
\newblock Cambridge University Press, in progress. Current version published on
  the web at http://php.indiana.edu/{$\sim$}rdlyons, 2004.

\bibitem{Ne}
J.~Neveu.
\newblock Arbres et processus de {G}alton-{W}atson.
\newblock {\em Ann. Inst. H. Poincar\'e}, 26:199--207, 1986.

\bibitem{Pau88}
F.~Paulin.
\newblock Topologie de {G}romov \'equivariante, structures hyperboliques et
  arbres r\'eels.
\newblock {\em Invent. Math.}, 94(1):53--80, 1988.

\bibitem{Pau89}
F.~Paulin.
\newblock The {G}romov topology on real-trees.
\newblock {\em Topology Appl.}, 32(3):197--221, 1989.

\bibitem{Pit02}
J.~Pitman.
\newblock {\em Combinatorial Stochastic Processes}.
\newblock Lecture Notes for St. Flour Course. Springer, 2002.

\bibitem{PitWink}
J.~Pitman and M.~Winkel.
\newblock Growth of the {B}rownian forest.
\newblock {\em to appear in Ann. Probab.}, 2005.

\bibitem{Sal}
P.~Salminen.
\newblock Cutting markovian trees.
\newblock {\em Ann.Acad. Scient. Fenn. A. I. Math.}, 17:123--137, 1992.

\end{thebibliography}

\end{document}